\documentclass[11pt,letterpaper]{article}
\usepackage[margin=1in]{geometry}
\usepackage[utf8]{inputenc} 
\usepackage[T1]{fontenc}    
\usepackage{dsfont}
\usepackage[maxfloats=100]{morefloats}[2015/07/22]
\usepackage{multirow}
\usepackage{wrapfig}
\usepackage[T1]{fontenc}
\usepackage{lmodern}

\usepackage{booktabs}       
\usepackage{amsfonts}       
\usepackage{nicefrac}       
\usepackage{microtype}      
\usepackage{enumerate}
\usepackage{lipsum}
\usepackage{mathtools}
\usepackage{cuted}
\usepackage{float}
\usepackage[dvipsnames,table,xcdraw]{xcolor}
\usepackage{bbm}
\usepackage{varioref}
\usepackage{hyperref}
                                
\usepackage{amssymb} 
\usepackage{rotating}
\usepackage{graphicx}
\usepackage{subcaption}
\usepackage[normalem]{ulem}

\usepackage{amsthm}
\DeclareUnicodeCharacter{00A0}{~}
\usepackage{algorithm}
\usepackage{algorithmic}

\newtheorem{assumption}{Assumption}
\newtheorem{theorem}{Theorem}
\newtheorem{lemma}{Lemma}
\newtheorem{proposition}{Proposition}
\newtheorem{definition}{Definition}
\newtheorem{corollary}{Corollary}
\theoremstyle{plain}
\newtheorem{remark}{Remark}
\RequirePackage[capitalize,nameinlink]{cleveref}[0.19]

\usepackage[textsize=small]{todonotes}

\usepackage{tcolorbox}
\usepackage{adjustbox}
\usepackage{hyperref}
\usepackage{url}
\usepackage{pdflscape}
\usepackage{graphicx} 
\usepackage{multirow}
\usepackage{multirow}
\usepackage{longtable}
\usepackage{longtable}
\usepackage{lipsum}
\usepackage{amsfonts}
\usepackage{sidecap}
\usepackage{float}  
\usepackage{rotating}
\usepackage{graphicx}
\usepackage{xcolor}
\usepackage{epstopdf}
\usepackage{wrapfig}
\usepackage{bm}

\usepackage{rotating} 
\usepackage{algorithmic}
\ifpdf
  \DeclareGraphicsExtensions{.eps,.pdf,.png,.jpg}
\else
  \DeclareGraphicsExtensions{.eps}
\fi

\newcommand{\ch}[1]{{\color{black}#1}}

\usepackage{makecell}
\newcommand{\newch}[1]{{\color{black}#1}}
\newcommand{\newchr}[1]{{\color{black}#1}}

\newcommand{\ser}[1]{\color{black}{#1}}
\newcommand{\fyr}[1]{\color{black}{#1}}

\newcommand{\ssg}[1]{\color{black}{#1}}

\crefname{section}{section}{sections}
\crefname{subsection}{subsection}{subsections}
\Crefname{section}{Section}{Sections}
\Crefname{subsection}{Subsection}{Subsections}

\Crefname{figure}{Figure}{Figures}

\crefformat{equation}{\textup{#2(#1)#3}}
\crefrangeformat{equation}{\textup{#3(#1)#4--#5(#2)#6}}
\crefmultiformat{equation}{\textup{#2(#1)#3}}{ and \textup{#2(#1)#3}}
{, \textup{#2(#1)#3}}{, and \textup{#2(#1)#3}}
\crefrangemultiformat{equation}{\textup{#3(#1)#4--#5(#2)#6}}%
{ and \textup{#3(#1)#4--#5(#2)#6}}{, \textup{#3(#1)#4--#5(#2)#6}}{, and \textup{#3(#1)#4--#5(#2)#6}}

\Crefformat{equation}{#2Equation~\textup{(#1)}#3}
\Crefrangeformat{equation}{Equations~\textup{#3(#1)#4--#5(#2)#6}}
\Crefmultiformat{equation}{Equations~\textup{#2(#1)#3}}{ and \textup{#2(#1)#3}}
{, \textup{#2(#1)#3}}{, and \textup{#2(#1)#3}}
\Crefrangemultiformat{equation}{Equations~\textup{#3(#1)#4--#5(#2)#6}}%
{ and \textup{#3(#1)#4--#5(#2)#6}}{, \textup{#3(#1)#4--#5(#2)#6}}{, and \textup{#3(#1)#4--#5(#2)#6}}

\crefdefaultlabelformat{#2\textup{#1}#3}

\title{\LARGE \bf
On iteratively regularized  first-order methods for simple bilevel optimization
}
\author{Sepideh Samadi 
\and Daniel Burbano 
\and Farzad Yousefian 
\thanks{The authors are affiliated with Rutgers University, Piscataway, NJ 08854, USA. They are contactable at sepideh.samadi@rutgers.edu, daniel.burbano@rutgers.edu, and farzad.yousefian@rutgers.edu, respectively.  A substantially preliminary version of this work appeared in the proceedings of the 2024 American Control Conference \cite{samadi2024achieving}.
This work was funded in part by the NSF under CAREER grant ECCS-$2323159$, in part by the ONR under grant N$00014$-$22$-$1$-$2757$, and in part by the DOE under grant DE-SC$0023303$.}
}

\begin{document}
\sloppy
\maketitle
\thispagestyle{empty}
\pagestyle{plain}
\maketitle

\begin{abstract}
We consider simple bilevel optimization \ser{(SBO)} problems where the goal is to compute among the optimal solutions of a composite convex optimization problem, one that minimizes a secondary objective function. Our main contribution is threefold. \ser{(i) When the upper-level objective is composite \ssg{and} strongly convex, we propose IR-ISTA$_\text{s}$, an iteratively regularized proximal gradient method with a prescribed update rule for the regularization parameter. We establish asymptotic convergence of the iterates to the unique optimal solution and \fyr{derive} simultaneous sublinear rates for suitably defined infeasibility and suboptimality error metrics.} \ser{(ii) For the same setting, we propose IR-VFISTA$_\text{s}$, an iteratively regularized accelerated proximal gradient method, establish its asymptotic convergence, and,} \fyr{under weak sharp minimality of the lower-level problem, derive faster simultaneous rates than IR-ISTA$_\text{s}$. These appear to be the best-known convergence rate guarantees for SBO problems with a strongly convex upper-level objective and improve upon the rates previously established for methods requiring convexity of both levels. (iii) When the upper-level objective is smooth and nonconvex, we propose IPR-VFISTA$_\text{nc}$, an inexactly projected iteratively regularized accelerated gradient method, and establish an asymptotic stationarity guarantee. Moreover, under lower-level weak sharp minimality, we derive nonasymptotic error bounds for both the stationarity and infeasibility error metrics.} We present preliminary numerical experiments on three ill-posed linear inverse problems \ser{and an optimal classifier-selection problem.}
\end{abstract}

%
%
%




%

\section{Introduction}\label{sec:1}
In this paper, we consider a class of constrained optimization problems, called simple bilevel optimization (SBO), of the form
 \begin{align}
\min \ & \bar{f}(x) := f(x) + \omega_f (x), \quad \text{s.t.} \quad x\in {X^*_{ \bar{h}}}:=   \arg\min_{x \in \mathbb{R}^n} \bar{h}(x) := h(x) + \omega_h (x), \label{prob:uni_centr}
\end{align}
where the upper- and lower-level objectives have a composite structure.  Here,  $f: \mathbb{R}^n \rightarrow  \mathbb{R} $ is a smooth
(possibly nonconvex) function, $ h: \mathbb{R}^n \rightarrow  \mathbb{R} $ is a smooth convex function, and $\omega_f , \omega_h  : \mathbb{R}^n \rightarrow  (-\infty, \infty ]  $ are extended-valued nonsmooth convex functions  that may represent structural constraints or regularization terms.  SBO problems naturally arise in optimal solution selection, a fundamental approach for addressing ill-posed optimization problems in image processing, machine learning, and signal processing~\cite{friedlander2008exact}. Beyond ill-posed problems, optimal solution selection is crucial in training over-parameterized models~\cite{samadi2024achieving}, portfolio optimization~\cite{beck2014first}, and stability analysis in multi-agent systems~\cite{qiu2024iteratively, kaushik2023incremental,ebrahimi2024distributed}. A key challenge in addressing this class of problems is that the standard constraint qualification conditions, e.g., \ssg{the} Slater condition, often fail to hold~\cite{friedlander2008exact}. This shortcoming has recently motivated the need for the design and analysis of iterative methods for addressing this class of problems.
 \subsection{Related work}\ser{
Table~\ref{table:lit} summarizes the most relevant methods for the SBO problem in \eqref{prob:uni_centr}, including their problem settings and asymptotic and nonasymptotic guarantees. Let $\bar h^*$ and $\bar f^*$ denote the optimal values of the lower- and upper-level objectives, respectively. The table reports bounds on the infeasibility $|\bar h(x)-\bar h^*|$, the suboptimality $|\bar f(x)-\bar f^*|$, and, when available, the distance to the SBO solution set.
The study of SBO problems originates from Tikhonov regularization for ill-posed problems~\cite{tikhonov1963solution}. Early methods established mainly asymptotic convergence or rates for the lower-level infeasibility~\cite{solodov2007explicit,solodov2007bundle,beck2014first,helou2017,sabach2017first,amini2019iterative}. \ssg{Some of} \fyr{the} first simultaneous nonasymptotic guarantees for the upper- and lower-level objectives were obtained through iteratively regularized methods for optimization problems with VI constraints~\cite{kaushik2021method}. Related distributed extensions were developed in~\cite{yousefian2021bilevel,kaushik2023incremental,qiu2024iteratively}.
More recent studies have pursued broader problem classes and improved rates. These include level-set methods~\cite{doron2023methodology}, bilevel subgradient schemes~\cite{merchav2023convex}, iteratively regularized methods for simple bilevel VIs~\cite{samadi2025improved}, dynamic proximal-gradient schemes~\cite{merchav2024dynamic}, accelerated gradient methods~\cite{cao2024accelerated}, and conditional-gradient methods\ser{~\cite{jiang2023conditional,giang2025projection}}. \ssg{More recently, \cite{giang2026conditional} developed an anytime conditional-gradient method for stochastic and finite-sum SBO with smooth objectives and established asymptotic and
nonasymptotic convergence guarantees.}
In our preliminary work~\cite{samadi2024achieving}, we proposed a regularized proximal-gradient method for SBO problems with a composite lower-level objective and established simultaneous guarantees for both levels. Under weak sharp minimality and an appropriately chosen regularization parameter, the method attains rates matching the optimal complexity of accelerated single-level convex optimization~\cite{beck2009fast}. \fyr{The present work extends these results to broader upper- and lower-level settings and develops a suite of iteratively regularized methods for both convex and nonconvex regimes, establishing asymptotic convergence guarantees and fast simultaneous nonasymptotic rates.} 
}
\begin{sidewaystable}
    \centering
    \tiny  
    \caption{\small {Most relevant existing works for addressing the simple bilevel problem in \eqref{prob:uni_centr}.}}
    \label{table:lit}
 
    \begin{tabular}{|l|c|c|c|c|c|c|c|c|c|c|}
      \hline
     \multirow{3}{*}{Reference} & \multirow{2}{*}{w.sh.} & \multicolumn{2}{c}{U.L.P} & \multicolumn{2}{|c|}{L.L.P}& \multicolumn{5}{|c|}{Convergence rates} \\  
&& \multicolumn{2}{|c|}{}& \multicolumn{2}{|c|}{} & \multicolumn{5}{|c|}{} \\  
      \cline{3-11}
      & on $X_{\text{l}}^*$& \multirow{2}{*}{ Obj.f.} & \multirow{2}{*}{Map.} & \multirow{2}{*}{ Obj.f.} & \multirow{2}{*}{Map.}& \multicolumn{2}{|c|} {Suboptimality} & \multicolumn{2}{|c|}{Infeasibility} & \multirow{2}{*}{$\operatorname{dist} (w_k, X_{\text{u}}^*)$}\\
      \cline{7-10}
      &&&&&&L.B.&U.B.&L.B.&U.B.& \\
      \hline
Solodov \cite{solodov2007explicit}& $\times$ &LS.C. &- &LS.C.& - &- &- &-& - & Asym.\\
 \hline
Solodov \cite{solodov2007bundle} & $\times$  &NS.C.&- & NS.C.& -&- & -&- &- &Asym. \\
 \hline
Beck et al. \cite{beck2014first}  &  $\times$  & D.SC.& - & LS.C. &-&  - &- & -  &$\tfrac{1}{\sqrt{K}}$&Asym.\\
 \hline
Helou et al. \cite{helou2017}& $\times$  &C.&-& LS.C.&-& -&- & - &-&Asym.\\
\hline
  Sabach et al. \cite{sabach2017first}  & $\times$ & LS.SC.& -& Comp.C.&- &- & - & 0 & $\tfrac{1}{{K}}$&-\\
  \hline
Amini et al. \cite{amini2019iterative} &$\times$&ND.SC.&-&ND.C.&-& -   &- & - &$\tfrac{1}{K^b}$&Asym.\\
\hline
  \multirow{2}{*} { Kaushik et al. \cite{kaushik2021method} }&  {$\times$} &  LC.C. &-&  - &C.& -& $\tfrac{1}{\sqrt[4]{K}}$&- &$\tfrac{1}{\sqrt[4]{K}}$ &-\\
  \cline{2-11}
  &{$\times$} & LS.SC. & - &-& C.& -& - &- &- &Asym.\\
    \hline
     Doron et al. \cite{doron2023methodology}  & $\times$ & LS.C.NL. & -& Comp.C.&-&  -&$\tfrac{1}{\sqrt{K}}$& 0 &$\tfrac{1}{{K}}$&-\\
    \hline
\multirow{2}{*}{Merchav et al. \cite{merchav2023convex}}  & $\times$ &Comp.C. &-& Comp.C. &-& -  & $\tfrac{1}{K^{1-c}}$&   0 &$\tfrac{1}{K^{c}}$, Asym. & Asym.\\
\cline{2-11}
& $\times$& Comp.SC.&- &Comp.C. &-& -&   $\rho_1^{K^{1-c}}$&0 & $\tfrac{1}{K^c}$& -  \\
\hline
\multirow{4}{*}{\ser{Jiang et al. \cite{jiang2023conditional}}}
& \ser{$\times$}
& \ser{LS.C.} & \ser{-}
& \ser{LS.C.} & \ser{-}
& \ser{-}
& \ser{$\tfrac{1}{K}$}
& \ser{$0$}
& \ser{$\tfrac{1}{K}$}
& \ser{-}
\\
\cline{2-11}
& \ser{$\times$}
& \ser{LS.NC.} & \ser{-}
& \ser{LS.C.} & \ser{-}
& \ser{-}
& \ser{$\tfrac{1}{\sqrt{K}}$}
& \ser{$0$}
& \ser{$\tfrac{1}{\sqrt{K}}$}
& \ser{-}
\\
\cline{2-11}
& \ser{$m\geq 1$}
& \ser{LS.C.} & \ser{-}
& \ser{LS.C.} & \ser{-}
& \ser{$-\tfrac{1}{K^{1/m}}$}
& \ser{$\tfrac{1}{K^{1/m}}$}
& \ser{$0$}
& \ser{$\tfrac{1}{K}$}
& \ser{-}
\\
\cline{2-11}
& \ser{$m\geq 1$}
& \ser{LS.NC.} & \ser{-}
& \ser{LS.C.} & \ser{-}
& \ser{$-\tfrac{1}{K^{1/(m+1)}}$}
& \ser{$\tfrac{1}{K^{1/(m+1)}}$}
& \ser{$0$}
& \ser{$\tfrac{1}{K^{m/(m+1)}}$}
& \ser{-}
\\
\hline
\multirow{7}{*}{Samadi et al. \cite{samadi2025improved}} & $\times$ & - & LS.C.&-& LS.C.& Asym.& $\tfrac{1}{\sqrt{K}}$&0 &$\tfrac{1}{\sqrt{K}}$&-\\
                                     & $m \geq 1$& - &LS.C.&-&LS.C.&$-\sqrt[2m]{\tfrac{1}{{K}}}$&$\tfrac{1}{\sqrt{K}}$&0 &$\tfrac{1}{\sqrt{K}}$&-\\
 \cline{2-11}
  &$\times$&- &LS.SC.& -&LS.C.&Asym.&$\tfrac{1}{K}$&0&$\tfrac{1}{K}$ &Asym.\\
    &$m\geq 1$&- &LS.SC.& -&LS.C.&$- \sqrt[m]{\tfrac{1}{K}}$&$\tfrac{1}{K}$&0&$\tfrac{1}{K}$ &$ \sqrt[m]{\tfrac{1}{K}}$\\
     &$\times$&- &LS.SC.& -&LS.SC.&-&$\tfrac{1}{K^p}$&0&$\tfrac{1}{K^{p+1}}+\tfrac{1}{K}$ &-\\
    &$m\geq 1$&- &LS.SC.& -&LS.C.&$-\sqrt[m]{\tfrac{1}{K^{p+1}}+\tfrac{1}{K}}$&$\tfrac{1}{K^p}$&0&$\tfrac{1}{K^{p+1}}+\tfrac{1}{K}$ &$\sqrt[m]{\tfrac{1}{K^{p+1}}+\tfrac{1}{K}}$\\
   \cline{2-11}
      &$m \geq 1$&LS.NC. &-& -&LS.C.&0&$\tfrac{1}{\sqrt{K}}$&0&$\tfrac{1}{K\sqrt{K}}$ &-\\
 \hline
 \multirow{2}{*}{Samadi et al. \cite{samadi2024achieving} } & $\times$ & LS.C. & - &Comp.C.&-& - &$\tfrac{1}{K}$& 0 &$\tfrac{1}{K}$& -\\
 & $m = 1$ & LS.C. & -&Comp.C.&-&$-\tfrac{1}{K^2}$& $\tfrac{1}{K^2}$ &0& $\tfrac{1}{K^2}$ &$\tfrac{1}{K^2}$ \\
     \hline
    \multirow{2}{*}{\ser{Merchav et al. ~\cite{merchav2024dynamic}}}  & $\times$ & Comp.C. & -& Comp.C.&-& - &$\tfrac{1}{K}$& 0 & $\tfrac{1}{K}$&-\\
   & $m=2$ & Comp.C. & -& Comp.C.&-& $\ser{-\tfrac{1}{K}}$ &$\ser{\tfrac{1}{K^{2-d}}}$& 0 & $\ser{\tfrac{1}{K^d}}$&Asym.\\
     \hline
   \multirow{2}{*}{Cao et al. \cite{cao2024accelerated}}  & $\times$ & LS.C. &  -&LS.C.&-& - &$\tfrac{1}{K^2}$&- &  $\tfrac{1}{K}$ &-\\
  &$m \geq 1$ & LS.C. &-  &LS.C.&-& $-\tfrac{1}{K^{2/2m-1}}-{\tfrac{1}{K^{2/m}}}$ &$\tfrac{1}{K^{2m/2m-1}}+\tfrac{1}{K^2}$&- &$\tfrac{1}{K^{2m/2m-1}}+\tfrac{1}{K^2}$ &-\\
     \hline
     \multirow{4}{*}{ Giang et al. \cite{giang2025projection}}  & $\times$ & LS.C. & -&  LS.C.&- & Asym.& Asym.&  Asym. & Asym.&-\\
     & $\times$ & LS.C. & -&  LS.C.&- & -& $\tfrac{1}{K^{1-q}}$ & - & $\tfrac{1}{K^{q}}$ &-\\
      & $m=2$ & LS.C. & -&  LS.C.&- & -& $\tfrac{1}{K^{\min\{q,1-q\}}}$ & - & $\tfrac{1}{K^{\min\{2q,1\}}}$ &-\\
      & $m=2$ & LS.SC. & -&  LS.C.&- & -& $\tfrac{1}{K^{\min\{\tfrac{2}{3}, 1-\tfrac{2q}{3},2-2q\}}}$ & - &$\tfrac{1}{K^{\min\{2q,1\}}}$  &-\\
      \hline
    \ssg{\multirow{4}{*}{Giang-Tran et al. \cite{giang2026conditional}}}
& \ssg{$\times$}
& \ssg{LS.C. } & \ssg{-}
& \ssg{LS.C.} & \ssg{-}
& \ssg{-}
& \ssg{$\tfrac{1}{K^{1/2}}$}
& \ssg{0}
& \ssg{$\tfrac{1}{K^{1/2}}$}
& \ssg{-}
\\
\cline{2-11}
& \ssg{$\times$}
& \ssg{LS.C. } & \ssg{-}
& \ssg{LS.C.} & \ssg{-}
& \ssg{Asym.}
& \ssg{Asym.}
& \ssg{0}
& \ssg{Asym.}
& \ssg{-}
\\
\cline{2-11}
& \ssg{$\times$}
& \ssg{LS.NC. } & \ssg{-}
& \ssg{LS.C.} & \ssg{-}
& \ssg{-}
& \ssg{$\tfrac{1}{K^{1/4}}$}
& \ssg{0}
& \ssg{$\tfrac{1}{K^{1/4}}$}
& \ssg{-}
\\
\cline{2-11}
& \ssg{$\times$}
& \ssg{LS.NC. } & \ssg{-}
& \ssg{LS.C.} & \ssg{-}
& \ssg{Asym.}
& \ssg{Asym.}
& \ssg{0}
& \ssg{Asym.}
& \ssg{-}
\\



\hline
    \end{tabular}
\par
\raggedright
{\tiny
\noindent{\ser{Notation:} w.sh. denotes weak sharp \ssg{minimality} of order $m$; \ssg{U.L.P. and L.L.P. denote upper- and lower-level problems, respectively, while U.L. and L.L. denote upper- and lower-level, respectively;}
LS., NS., and S. denote L-smooth, nonsmooth, and smooth functions, respectively; SC., C., and NC. denote strongly convex, convex, and nonconvex functions, respectively; LC. denotes Lipschitz continuity; D. and ND. denote differentiable and nondifferentiable cases, respectively; Asym. denotes asymptotic.  
The notation $w_K$ represents the output of the method, with parameter ranges: $b \in (0,0.5)$, $c \in (0.5,1)$, $\rho_1, \rho_2, \rho_3 \in (0,1)$, $d \in (1,2)$, $q \in (0,1)$, and $p \geq 1$. Comp. denotes composite functions. 
 NL. denotes norm-like functions, Obj.f. refers to the objective function, and Map. represents the mapping $F$ in $\mbox{VI}(X,F)$. The sets $X_{\text{u}}^*$ and $X_{\text{l}}^*$ denote the optimal solutions of SBO and \ser{the} lower-level problem, respectively. \ssg{For nonconvex cases, the error metric used is the stationarity measure specified in the corresponding reference.} 
 \ser{In \cite{jiang2023conditional}, we let $\epsilon_f=\epsilon_h = \epsilon$ and also \ssg{$K$ refers to the main CG-BiO iterations, excluding the initial preprocessing iterations.} \ssg{For \cite{giang2026conditional}, we report  the results for the finite-sum setting.}
 }
 }
 \par}
\end{sidewaystable}

 \subsection{Contributions}
Our main contributions are presented in the following and are also concisely summarized in Table~\ref{table:contributions_all}. \\
  {(i)} {\it An iteratively regularized proximal method with new guarantees for composite SBO with a strongly convex upper-level objective.}  We propose \ser{Iteratively Regularized ISTA} \ser{(IR-ISTA$_\text{s}$)} for addressing SBO problems with a composite strongly convex upper-level objective function and a composite convex lower-level objective function.
Under a diminishing regularization update rule, we show that the generated iterate converges asymptotically to the unique solution of the SBO problem. Further, we establish simultaneous sublinear convergence rates for infeasibility and the upper bound \ssg{on} suboptimality.  Under a weak sharp minimality assumption, we derive \fyr{a nonasymptotic lower bound on the suboptimality gap}. We also extend the rate analysis to the setting with a \ser{carefully designed} constant regularization parameter, where we refer to the method as R-ISTA$_\text{s}$.  Under a weak sharp minimality assumption for the lower-level problem, 
R-ISTA$_\text{s}$ attains a linear convergence rate. \fyr{These results provide new convergence guarantees for this class of composite SBO problems}. Importantly, when compared with existing methods in Table~\ref{table:lit}, IR-ISTA$_\text{s}$ is among the first IR schemes that is equipped with both asymptotic and (simultaneous) \ser{nonasymptotic} convergence guarantees for resolving SBO problems. \\
  {{(ii)}} {\it \ser{An iteratively} regularized accelerated proximal method for composite SBO with a strongly convex upper-level objective.} \ssg{To obtain faster convergence rates than those in (i) in certain settings,} \ser{we propose \ser{Iteratively Regularized Variant of FISTA} \ser{(IR-VFISTA$\text{s}$)}. \ssg{We also establish the asymptotic convergence of the proposed method.} 
  \fyr{When the lower-level problem admits a} \ssg{weak sharp minimality} \ssg{of} order $1 < m\leq 2$, \fyr{we derive simultaneous infeasibility and suboptimality convergence rates. These rates appear to be the best-known convergence rate guarantees for SBO problems with a strongly convex upper-level objective and improve upon the rates previously established for methods requiring convexity of both levels.} With a carefully designed constant regularization parameter, referred to as R-VFISTA$\text{s}$, we establish quadratically decaying sublinear convergence rates for both infeasibility and suboptimality error metrics.} When weak sharp minimality holds,  a linear convergence rate with an improved dependence on the condition number is achieved. 
\\
  {(iii)} \ser{{\it New asymptotic and nonasymptotic}  guarantees for composite SBO problems with a smooth nonconvex upper-level objective}. When the upper-level objective is smooth and nonconvex, we propose the Inexactly Projected Regularized VFISTA method ($\text{IPR-VFISTA}_{\mathrm{nc}}$). \ser{We establish \fyr{an} asymptotic stationarity guarantee for this class of SBO problems without requiring the \ssg{weak sharp minimality} of the lower-level problem, together with an asymptotic lower-level infeasibility guarantee. Under \fyr{lower-level weak sharp minimality}, we further derive  nonasymptotic guarantees for both the stationarity and infeasibility error metrics.} \fyr{The total iteration complexity of IPR-VFISTA$_{\mathrm{nc}}$ matches that of existing methods, while providing a sharper lower-level infeasibility guarantee.}

\begin{table*}[t]
\centering
\tiny
\caption{\small Summary of the main contributions in this work \ser{addressing \eqref{prob:uni_centr}.}}
\label{table:contributions_all}

\setlength{\tabcolsep}{3.8pt}
\renewcommand{\arraystretch}{1.15}

\begin{tabular}{|c|c|c|c|c|c|c|}
\hline

\multicolumn{7}{|c|}{Composite strongly convex U.L. and composite convex L.L.} \\ \hline

\multirow{2}{*}{%
  \parbox[c]{1.8cm}{\centering Our method}%
}
& \multirow{1}{*}{M.A.}
& \multicolumn{2}{c|}{$\bar{f}(w_K)-\bar{f}^*$}
& \multicolumn{2}{c|}{$\bar{h}(w_K)-\bar{h}^*$}
& \multirow{2}{*}{${\|w_K-x^*\|}^2_2$}
\\ \cline{3-6}

& w.sh.
& L.B.
& U.B.
& L.B.
& U.B.
&
\\ \hline

\multirow{2}{*}{
\begin{tabular}[c]{@{}c@{}}
IR-ISTA$_{\text{s}}$\\[-1pt]
\tiny \ser{Thm.~\ref{proposition:F_eta_IR-ISTA-S}},\\[-1pt] \ser{Cor.~\ref{theorem:ir-ista-s}}
\end{tabular}}
& $\times$
& Asym.
& $\tfrac{1}{K}$
& 0
& $\tfrac{1}{K}$
& Asym.
\\ \cline{2-7}

& $m \geq 1$
& $-\sqrt[m]{\tfrac{1}{K}}$
& $\tfrac{1}{K}$
& 0
& $\tfrac{1}{K}$
& $\sqrt[m]{\tfrac{1}{K}}$
\\ \hline

\multirow{3}{*}{
\begin{tabular}[c]{@{}c@{}}
R-ISTA$_{\text{s}}$\\[-1pt]
\tiny \ser{Cor.~\ref{cor:r-ista-s}}
\end{tabular}}
& $\times$
&  -
& $\tfrac{1}{K^{q}}$
& 0
& $\tfrac{1}{K^{q+1}}+\tfrac{1}{K}$
&  -
\\ \cline{2-7}

& $m \geq 1$
& $-\sqrt[m]{\tfrac{1}{K^{q+1}}+\tfrac{1}{K}}$
& $\tfrac{1}{K^{q}}$
& 0
& $\tfrac{1}{K^{q+1}}+\tfrac{1}{K}$
& $\tfrac{1}{K^q}+\sqrt[m]{\tfrac{1}{K^{q+1}}+\tfrac{1}{K}}$
\\ \cline{2-7}

& $m = 1$
& $-\tfrac{1}{\eta}(1-\tfrac{\eta}{\kappa})^K$
& $\tfrac{1}{\eta}(1-\tfrac{\eta}{\kappa})^K$
& 0
& $(1-\tfrac{\eta}{\kappa})^K$
& $\tfrac{1}{\eta \mu_f}(1-\tfrac{\eta}{\kappa})^K$
\\ \hline

\multirow{2}{*}{
\begin{tabular}[c]{@{}c@{}}
\ser{IR-VFISTA$_{\text{s}}$}\\[-1pt]
\ser{\tiny Thm.~\ref{theorem:IRVFISTA_var_rates}}
\end{tabular}}
& \ser{$\times$}
& \ser{Asym.}
& \ser{Asym.}
& \ser{0}
& \ser{Asym.}
& \ser{Asym.}
\\ \cline{2-7}

& \ser{$1 <m \leq 2$}
& \ser{$-\tfrac{1}{K^{\tfrac{a}{m-1}}}$}
&  \ser{$\tfrac{1}{K^{1+\tfrac{a(3-m)}{2(m-1)}}}$}
& \ser{0}
& \ser{$\tfrac{1}{K^{\tfrac{am}{m-1}}}$}
& \ser{$\tfrac{1}{K^{\tfrac{a}{m-1}}}$}
\\ \hline

\multirow{3}{*}{
\begin{tabular}[c]{@{}c@{}}
$\text{R-VFISTA}_{\text{s}}$\\[-1pt]
\tiny \fyr{Cor.~\ref{Cor:R-VFISTA}}
\end{tabular}}
& $\times$
& -
& $\tfrac{1}{K^{p-1}}$
& 0
& $\tfrac{1}{K^2}$
& -
\\ \cline{2-7}

& $m \geq 1$
& $-\sqrt[m]{\tfrac{1}{K^2}}$
& $\tfrac{1}{K^{p-1}}$
& 0
& $\tfrac{1}{K^2}$
& $\tfrac{1}{K^{p-1}}+\sqrt[m]{\tfrac{1}{K^2}}$
\\ \cline{2-7}

& $m = 1$
& $-\tfrac{1}{\eta}(1-\tfrac{1}{\sqrt{\kappa_{\eta}}})^K$
& $\tfrac{1}{\eta}(1-\tfrac{1}{\sqrt{\kappa_{\eta}}})^K$
& 0
& $(1-\tfrac{1}{\sqrt{\kappa_{\eta}}})^K$
& $\tfrac{1}{\eta \mu_f}(1-\tfrac{1}{\sqrt{\kappa_{\eta}}})^K$
\\ \hline

\end{tabular}

\begin{tabular}{|c|c|c|c|c|c|c|c|}
\hline

\multicolumn{8}{|c|}{\ser{Smooth nonconvex U.L.   and composite convex L.L.}} \\ \hline

\multirow{2}{*}{Our method}
& \multicolumn{2}{c|}{\ser{M.A.}}
& \multicolumn{2}{c|}{$\|G_{1/\hat{\gamma}}(\hat{w}_K^*)\|^2$}
& \multicolumn{2}{c|}{$\mathrm{dist}(w_K,X_{\bar h}^*)$}
& \multirow{2}{*}{\ser{Total complexity of $\epsilon$-stationarity}}
\\ \cline{2-3} \cline{4-7}

& \multirow{1}{*}{\ser{w.sh.}}
& \ser{Com. L.L. cons.}
& \ser{L.B.}
& \ser{U.B.}
& \ser{L.B.}
& \ser{U.B.}
&
\\ \hline

\multirow{2}{*}{
\begin{tabular}[c]{@{}c@{}}
$\text{IPR-VFISTA}_{\text{nc}}$\\[-1pt]
\tiny \ser{Thm.~\ref{theorem:asym_non_cons},} \\[-1pt]\ser{Thm.~\ref{theorem:nonc_all}}
\end{tabular}}
& \ser{$\times$}
& \ser{$\checkmark$}
& \ser{0}
& \ser{Asym.}
& \ser{0}
& \ser{Asym.}
& \ser{Asym.}
\\ \cline{2-8}


& \ser{$1 \ser{<} m \ser{<} 2$}
& \ser{$\checkmark$}
&  \multirow{2}{*}{\ser{0}}
&  \multirow{2}{*}{\ser{$\tfrac{1}{K}$}}
&  \multirow{2}{*}{\ser{0}}
&  \multirow{2}{*}{\ser{$\tfrac{1}{K^{m\ser{+\delta}}}$}} 
& \multirow{2}{*}{\ser{$\mathcal{O}(\epsilon^{-(m+1+\delta)})$}}
\\ \cline{2-3}

& \ser{$m=2$}
& \ser{$\times$}
& 
& 
& 
& 
& 
\\ \hline

\end{tabular}


\par
\raggedright
{\tiny
\noindent\tiny {Notation:  w.sh. denotes $\alpha$-weak sharp minima of order $m$;
$p > 2$; $q > 0$; \ser{$0<a<2$; $0<b<1$;} $w_K$ is output of the method;
L.B. and U.B. stand for lower bound and upper bound, respectively; Asym. denotes asymptotic convergence; M.A. denotes main assumption on; We ignore logarithmic numerical factors; We consider both $ m \geq 1 $ and $ m = 1 $, the rates for $ m = 1 $ are included within the analysis for $ m \geq 1 $. Additionally, we derive improved rates by explicitly setting $ m = 1 $ in a separate case in which we assume that $\eta$ falls below a threshold. \ser{Let $\delta>0$ be an arbitrarily small constant.} \ser{Com. L.L. cons.  denotes that
$\omega_h=\mathbb{I}_X$ for a nonempty compact convex set $X$.}}
\par}
\end{table*}

\ser{\subsection{Outline of the paper}
Section~\ref{Sec: 2} presents the preliminaries. Sections~\ref{sec: 4} and~\ref{sec:5} analyze composite SBO problems with a composite convex lower-level objective and, respectively, a composite strongly convex or smooth nonconvex upper-level objective. Section~\ref{sec: 6} presents the numerical results, and Section~\ref{sec:7} concludes the paper.}
\subsection{Notation}
\ser{For $x,y\in\mathbb{R}^n$, $\langle x,y\rangle$ and $x^\top$ denote the inner product and transpose. The $\ell_p$-norm is denoted by $\|\bullet\|_p$, with $\|\bullet\|:=\|\bullet\|_2$. For $g:\mathbb{R}^n\to(-\infty,\infty]$, $\operatorname{prox}_{g}[x]$ denotes its proximal map, defined in Definition~\ref{def1}. For $f:\mathbb{R}^n\to\mathbb{R}$, $\nabla f(x)$ denotes its gradient. For convex $f$, $\tilde{\nabla}f(x)$ is a subgradient at $x$ if $f(y)\geq f(x)+\langle\tilde{\nabla}f(x),y-x\rangle$ for all $y\in\operatorname{dom}(f)$, and $\partial f(x)$ is its subdifferential.}
\ser{For a set $X$, $\Pi_X[x]$ denotes the Euclidean projection of $x$ onto $X$; if $X$ is nonempty, closed, and convex, $\operatorname{dist}(x,X):=\|x-\Pi_X[x]\|$. For \eqref{prob:uni_centr}, $X^*$ and $\bar f^*$ denote its \fyr{optimal} solution set and optimal \fyr{objective} value. Also, $X_{\bar h}^*:=\operatorname{argmin}_{x\in\mathbb{R}^n}\bar h(x)$ and $\bar h^*:=\inf_{x\in\mathbb{R}^n}\bar h(x)$. The relative interior and interior of $C$ are denoted by $\operatorname{ri}(C)$ and $\operatorname{int}(C)$.}
\ser{Finally, $(\cdot)^\dagger$, $\mathbb{I}_S(x)$, $\hat C_{\bar f}:=\inf_{x\in\mathbb{R}^n}\bar f(x)$, and $[\bullet]_+:=\max\{0,\bullet\}$ denote the Moore--Penrose pseudoinverse, indicator function of $S$, infimum of $\bar f$, and positive-part operator, respectively.}

 \section{Preliminaries}\label{Sec: 2}
 We present some definitions and preliminary
results.
\begin{definition}[{\cite[Definition 6.1]{beck2017first}}]\em \label{def1}
Given a function $g: \mathbb{R}^n \rightarrow (-\infty, \infty]$, its proximal map is given as $\text{prox}_{g}[x] := {\text{argmin}}_{u \in \mathbb{R}^n}\ \{ g(u) + \tfrac{1}{2}{\|u-x\|}^2 \},$ for all $ x\in \mathbb{R}^n$.
\end{definition}
\begin{lemma}[{\cite[Theorem 6.39]{beck2017first}}]\label{lem:RAPM1}\em 
Given a proper, closed, and convex function $g: \mathbb{R}^n \rightarrow (-\infty, \infty]$,   $z:= \text{prox}_{\gamma g}[u]$ if and only if for any $u\in \mathbb{R}^n$ and $\gamma > 0$, we have $(u-z) \in \gamma \partial g(z)$. 
\end{lemma}
\begin{definition}\label{def:uni_funcs}\em
Consider the problem in~\eqref{prob:uni_centr} where $f, h: \mathbb{R}^n \rightarrow  \mathbb{R} $   and $\omega_f, \omega_h : \mathbb{R}^n \rightarrow  (-\infty, \infty ]  $ are  given functions. For ${\eta}, \gamma >0$ and any $ x\in\mathbb{R}^n$, we define 
\begin{equation*}
\begin{aligned}
    g_{\eta}(x) &:= h(x) + {\eta} f(x), 
    & \quad \omega_{{\eta}}(x) &:= \omega_h(x) + {\eta} \omega_f(x), \\
    \bar{g}_{\eta}(x) &:= g_{\eta}(x) + \omega_{{\eta}}(x), 
    & \quad \hbox{ and }\quad q_{{\eta}}(x) &:= \text{prox}_{\gamma {\omega_{{\eta}}}}\left(x - \gamma \left(\nabla h (x) + {\eta} \nabla f(x)\right)\right).
\end{aligned}
\end{equation*}

\end{definition}

\begin{definition}[{\cite[Definition 2.13]{beck2017first}}]\em
A proper function $f:\mathbb{R}^n \to (-\infty, \infty]$ is called coercive, if $\lim_{\|x\| \to \infty} f(x) = \infty$.
\end{definition}
\begin{definition}[Weak sharp minima{~\cite[Definition 1.1]{studniarski1999weak}}] \label{def:weaksharp} \em
Consider the problem $\min_{x \in \mathbb{R}^n} f(x)$, where $f: \mathbb{R}^n \to (-\infty, \infty]$. Let $X^*_f := \arg\min_{x \in \mathbb{R}^n} f(x)$ be a nonempty set. \ser{The set $X_f^*$ satisfies the \ssg{weak sharp minimality} of order $m \geq 1$,}
if there exists a constant $\alpha > 0$ such that  
$
f(x) - \inf_{x \in \mathbb{R}^n} f(x) \geq \alpha \, \text{dist}^m(x, X^*_f)$, for all $x \in \mathbb{R}^n.
$
\end{definition}
Under some non-degeneracy conditions, the optimal solution set of linear programs and linear \ssg{complementarity} problems admits weak sharp minima of order $m=1$~\cite{studniarski1999weak}. Further, quadratic programs under some assumptions admit the weak sharp minimality of order $m=1$~\cite[Section 3]{burke1993weak}. In nonlinear programming, this condition is also referred to as the H\"older continuity property of the solution
set, e.g., see \cite{ioffe1994sensitivity,bonnans1995quadratic, jiang2022holderian, bolte2017error}. The weak sharp minimality for problems with a unique optimal solution is also studied in \cite{auslender1984stability, studniarski1986necessary, ward1994characterizations}. \fyr{In the following lemma, we present two examples that satisfy the weak sharp minimality condition (see also \cite{jiang2022holderian,bolte2017error}).} \ssg{We provide the proof for the next lemma in  Appendix~\ref{subsection:proof_weaksharp_examples}.}
\begin{lemma}\label{lemma:weaksharp_examples}\em
{The following problems satisfy the weak sharp minimality condition in
Definition~\ref{def:weaksharp} with $m=2$ under the stated conditions.
}

\ser{(i) [constrained least-squares problem]} Let $\mathbf{A}\in\mathbb{R}^{n\times n}$ be nonzero,
symmetric, and positive semidefinite, and let $b\in\mathbb{R}^n$.
Consider the problem $\min_{\|x\|\leq 1} h(x):=\tfrac{1}{2}\|\mathbf{A}x-b\|^2.$
Define $\mathbf{Q}:=\mathbf{A}^\top\mathbf{A}$ and
$\mathbf{q}:=\mathbf{A}^\top b$, let $\lambda_{\min}$ denote the
smallest eigenvalue of $\mathbf{Q}$, and let $h^*$ denote the optimal
value of the problem. Assume that either
$\min_{x\in\mathbb{R}^n}h(x)<h^*$, or
$\min_{x\in\mathbb{R}^n}h(x)=h^*$ and either
$\lambda_{\min}>0$, or $\lambda_{\min}=0$ and
$\|\mathbf{Q}^{\dagger}\mathbf{q}\|<1$. 

\ser{(ii)} \fyr{[constrained logistic empirical risk problem] Let $\{(u_i,v_i)\}_{i=1}^N$ be a dataset with
$u_i\in\mathbb{R}^n$ and $v_i\in\{-1,+1\}$. Consider
$\min_{x\in X} f(x):=
\frac{1}{N}\sum_{i=1}^N
\log\bigl(1+\exp(v_i \fyr{\langle u_i, x\rangle})\bigr)$,
where
$\mathcal R:=\operatorname{span}\{u_1,\ldots,u_N\}\neq\{0\}$ and
$\mathcal N:=\mathcal R^\perp$.
Assume that $X=X_{\mathcal{R}}+X_{\mathcal{N}}$, where $X_{\mathcal{R}}\subset\mathcal{R}$ and $X_{\mathcal{N}}\subset\mathcal{N}$ are compact convex sets; This holds in common settings, such as box constraints, the $\ell_1$-ball, and the probability simplex.}
\end{lemma}
\begin{lemma} [{\cite[Theorem 27.2]{rockafellar1997convex}}]\label{lemma:convexity_X_h}\em
Let $f: \mathbb{R}^n \to (-\infty, \infty]$ be a proper closed convex function. Then, the optimal solution set of $\min_{x \in \mathbb{R}^n}  f(x)$ is convex.
\end{lemma}
\begin{lemma}\label{lemma:general_boundedness}\em 
  Let $b \in \mathbb{R}$ and $c,d > 0$ be given. Consider the \ssg{nonnegative}  sequence $\{r_k\}$ satisfying \ch{the recursion} $r_{k+1}
   \leq  b  +    \tfrac{c}{k^d}   r_{k}$, for any $k\geq 1$. 
   Let $\hat{r} := \fyr{\max\{\max_{1 \leq k \leq \lceil \sqrt[d]{2c} \rceil}r_k,2b \} }$. Then, $r_k \leq \ser{\hat{r}}$
  for all $k \geq 1$.
   \end{lemma}
   \begin{proof} 
First, we use mathematical induction to show that $r_k \leq \ser{\hat{r}}
$, for any $k$ such that $k^d\geq 2c$.  Let $k=\lceil\sqrt[d]{2c}\rceil$. By the definition of \ser{$\hat{r}$}, we have \ser{$r_k \leq \hat{r}$}.
Thus, the base case holds true. Now assume that $r_k$ is bounded by \ser{$\hat{r}$},
for some $k$ such that $\ssg{k^d \geq 2c}$. We aim to show that $r_{k+1}$ is bounded by $\ser{\hat{r}}$.
\ser{Regarding} the inductive hypothesis,
\ser{$r_{k+1}\leq b+\tfrac{c}{k^d}r_k \leq b +  \tfrac{c \hat{r}}{k^d} \leq  b +  \tfrac{c \hat{r}}{2c} \leq \tfrac{\hat{r}}{2} + \tfrac{\hat{r}}{2} =\hat{r}. $}
\ser{Hence, we have}  $r_{k+1} \leq \ser{\hat{r}}$. Since the finite number of initial terms of $\{ r_k\}$ for $ 1 \leq k \leq \lceil\sqrt[d]{2c}\rceil$ are captured by the definition of $\hat{r}$, it follows that $r_k \leq \ser{\hat{r}}
$ for all $k \geq 1$. 
 \end{proof}
 \ser{
\begin{lemma}\label{lem:eta_difference_bound}\em
Let $\eta_k=\tfrac{\eta_0}{(k+1)^a}$ for some $\eta_0>0$ and $a>0$, \fyr{for $k\geq 0$}. Then, for all $k\geq 0$, we have $\eta_k-\eta_{k+1}
\leq \tfrac{a\eta_0}{(k+1)^{a+1}}.$
\end{lemma}
\begin{proof}
Let $\phi(s)=s^{-a}$ for $s>0$. 
By the Mean Value Theorem applied to the function $\phi$ on $[k+1,k+2]$, there exists $\hat{k}\in(k+1,k+2)$ such that $\phi(k+2)-\phi(k+1)
=
\phi'(\hat{k})((k+2)-(k+1))
=
\phi'(\hat{k}).$
Therefore, $\tfrac{1}{(k+1)^a}-\tfrac{1}{(k+2)^a}
=
\phi(k+1)-\phi(k+2)
=
-\phi'(\hat{k})
=
\tfrac{a}{\hat{k}^{a+1}}.$
Also, since $k+1 < \hat{k}$ and $0<a$, we have
$
\tfrac{1}{\hat{k}^{a+1}}\leq \tfrac{1}{(k+1)^{a+1}}.
$
Therefore, we obtain $\eta_k-\eta_{k+1}\leq\tfrac{a\eta_0}{(k+1)^{a+1}}.$
\end{proof}
}
\ser{
\begin{lemma}[Chung's lemma~\cite{chung1954stochastic}]\label{lem:chung}\em
Let $\{u_k\}_{k\geq0}\subseteq\mathbb{R}_+$ satisfy $u_{k+1}
\leq
(1-\tfrac{c}{(k+1)^p})u_k
+
\tfrac{d}{(k+1)^{p+q}},
\ k\geq0,$
where $p\in(0,1]$ and $c,d,q>0$. If either $0<p<1$, or
$p=1$ and $c>q$, then $u_k
=
\mathcal{O}\left((k+1)^{-q}\right).$
\end{lemma}
}

 \section{Composite SBO with a strongly convex upper-level objective}\label{sec: 4}
In this section, we consider addressing the problem \ch{in}~\eqref{prob:uni_centr} under the following assumption. 
\begin{assumption}\label{assump:strongly_main}\em Consider  the problem in~\eqref{prob:uni_centr}. Let the following hold.\\
\noindent \ser{(i)} $ f: \mathbb{R}^n \rightarrow \mathbb{R}$ is  $L_f$-smooth and ${\mu_f}$-strongly convex.\\
\noindent \ser{(ii)}  $h: \mathbb{R}^n \rightarrow \mathbb{R}$ is  $L_h$-smooth and convex.\\
\noindent \ser{(iii)}  ${\omega_f}$ and ${\omega_h}: \mathbb{R}^n \rightarrow (-\infty, \infty]$ are proper, closed, and convex.\\
\noindent \ser{(iv)}  The set {$X^*_{\bar{h}}$} is nonempty.\\
\noindent \ser{(v)} $\hat{C}_{\bar{f}} = \inf_{x \in \mathbb{R}^n} \bar{f}(x) > -\infty$.\\
\noindent \ser{(vi)}  $x^* \in \text{int}(\text{dom}(\omega_f))$, where $x^*$ is the unique optimal solution to \eqref{prob:uni_centr}.\\
\noindent \ser{(vii)}  $\mathrm{ri}(\mathrm{dom}(\omega_f)) \cap X^*_{\bar{h}} \neq \emptyset$.
\end{assumption}
\begin{remark}\label{remark:nonempty_X}
\fyr{Assumption~\ref{assump:strongly_main}~(v) holds, for example, when $\bar{f}$ is coercive. This situation commonly arises when the upper-level objective acts as a regularizer (see Section~\ref{sec: 6}). We next present several sufficient conditions for the existence of an optimal solution to problem~\eqref{prob:uni_centr}. (I) By \cite[Prop. 2.1.1]{bertsekas2003convex}, $X^*_{\bar{h}}$ is nonempty and compact whenever one of the following conditions holds: (a) $\operatorname{dom}(\bar{h})$ is bounded; (b) there exists $\lambda\in\mathbb{R}$ such that $\{x\mid \bar{h}(x)\leq \lambda\}$ is nonempty and bounded; or (c) $\bar{h}$ is coercive. Consequently, under Assumption~\ref{assump:strongly_main}~\ssg{(i)-(iv)}, $X^*$  is nonempty (cf.~\cite[Thm. 2.12]{beck2017first}). (II) If $\omega_h$ is an indicator function of a compact set $C$, then under Assumption~\ref{assump:strongly_main}~(\ser{ii}) and~(iii), $X^*_{\bar{h}}$ is nonempty by~\cite[Thm. 2.12]{beck2017first}. Since $X^*_{\bar{h}} \subseteq C$, the set  $X^*_{\bar{h}}$ is compact. Thus, $X^*$ is nonempty by Assumption~\ref{assump:strongly_main}~(i) and~(iii) together with \cite[Thm.~2.12]{beck2017first}.
(III) If $\omega_f\equiv 0$, $\omega_h$ is the indicator function of a closed set, and both $h$ and $f$ are coercive functions, then $X^*$ is nonempty by~\cite[Thm. 2.14]{beck2017first}.}  
\end{remark}
\begin{remark}[Uniqueness of the optimal solution] The nonemptiness of $X^*_{\bar{h}}$, the convexity of this set (cf. Lemma~\ref{lemma:convexity_X_h}), and the strong convexity of $\bar{f}$ guarantee that the problem in \eqref{prob:uni_centr} admits a unique optimal solution (see~\cite[ Prop. 2.1.2]{bertsekas2003convex}).
\end{remark}
\begin{definition}\label{def:uni_L}\em
\fyr{Given a sequence $\{\eta_k\}$ of positive scalars, for each $k\ge0$, we define
$L_{\eta_k} := L_h + \eta_k L_f$.}
\end{definition}

\subsection{\texorpdfstring{The \ssg{$\text{IR-ISTA}_{\text{s}}$} method}{The IR-ISTA-s method}}
We propose Algorithm~\ref{alg:IR-ISTA-s} to address the problem in~\eqref{prob:uni_centr} under Assumption~\ref{assump:strongly_main}. \fyr{The proposed method, referred to as the} Iteratively Regularized  Iterative Shrinkage-Thresholding Algorithm ($\text{IR-ISTA}_\text{s}$)  \fyr{extends} the classical ISTA method~\cite[Section 10.5]{beck2009fast} \fyr{for} single-level composite \fyr{optimization. Specifically,} $\text{IR-ISTA}_\text{s}$  \fyr{is characterized by two features:} (i) an iterative regularization technique, whereby at each iteration $k$, $x_k$ is updated using the proximal operator applied to the regularized function $\bar{h}(\fyr{\cdot}) + \eta_k \bar{ f}(\fyr{\cdot})$; \fyr{and} (ii) a weighted averaging sequence \fyr{with geometrically varying weights $\{\theta_k\}$}.
\begin{algorithm}   
 	\caption{Iteratively Regularized ISTA (IR-ISTA$_\text{s}$)}
 	\begin{algorithmic}[1]
 		\STATE {{\bf input:} $\bar{x}_0=x_0 \in \mathbb{R}^n$, nonincreasing sequence $\{\eta_k\}$ for $k\geq 0$, stepsize $\gamma > 0$ such that $\gamma \leq  \tfrac{1}{L_h + {\eta_0} L_f}$, $\Gamma_0 = 0$, $\theta_0=\tfrac{1}{(1 - {\eta_0\gamma{\mu_f}  })}$, and  {$K\geq 1$.}}
 		\FOR {$k = 0,1, \dots, {K-1}$}
 		
 		\STATE   {$x_{k+1} = \text{prox}_{\gamma {\omega_{\eta_k}}}\left( x_{k} -\gamma \left(\nabla h( x_{k}) + {\eta_k} {\nabla} f( x_{k})\right)\right)$}	

       \STATE $\bar{x}_{k+1} = \tfrac{\Gamma_k \bar{x}_{k} + \eta_k\theta_{k} x_{k+1}}{\Gamma_{k+1}}$
        \STATE $\Gamma_{k+1} = \Gamma_{k} + \eta_k\theta_{k}$ \hbox{and} $\theta_{k+1}=\tfrac{\theta_k}{(1 - {\eta_{k+1}\gamma{\mu_f}  })}$
 		\ENDFOR
	\STATE {{\bf return:} {$\bar{x}_K$}}
 	\end{algorithmic}
  \label{alg:IR-ISTA-s}
 \end{algorithm}
In the next lemma, we study some properties of the sequence  $\{\theta_k\}$.
 \begin{lemma}\label{lemma:defining_theta}\em
 \fyr{Let $\{\theta_k\}$ be the sequence} generated by the recursive update rule in Algorithm~\ref{alg:IR-ISTA-s}. Then, the following \fyr{properties} hold.
 
\ser{\noindent (i)} $\theta_{k}=1/\textstyle\prod_{t=0}^{k}(1 - {\eta_t\gamma{\mu_f}  })$, for all $k \geq 0$. Further, $\theta_k >1$ for all $k\geq 0$.\\
\ser{\noindent (ii)} [diminishing regularization] Let $\eta_k = \tfrac{\eta_{0,u}}{\eta_{0,l}+k}$, where $\eta_{0,u} =(\gamma \mu_f)^{-1} $, $\eta_{0,l} = {2L_{f}}/{\mu_f}$, and $\gamma \leq 0.5/L_h$. Then, we have $\gamma \leq  \tfrac{1}{L_h + {\eta_0} L_f}$, 
$\theta_k =\tfrac{\eta_{0,l}+k}{\eta_{0,l}-1}
$ for all $k\geq 0$, $ \textstyle\sum_{k=0}^{K-1}\theta_k \eta_k = \tfrac{K}{\gamma (2L_f - \mu_f)},$ \fyr{and}
 $ \textstyle\sum_{j=0}^{K-1}\eta_j^2 \theta_j = \tfrac{\eta_{0,u}^2}{(\eta_{0,l}-1)} (\fyr{\frac{1}{\eta_{0,l}}}+\ln (\tfrac{K+\eta_{0,l}-1}{\eta_{0,l}})).
$ 

\ser{\noindent (iii)} [constant regularization]  Let $\eta={(p+1)\ln(K)}/({\gamma \mu_f K})$, for some $p > 0$ and $K> 1$ such that $\tfrac{K}{\ln(K)} \geq 2(p+1)\tfrac{L_f}{\mu_f}$. Suppose $\gamma \leq \fyr{\tfrac{1}{2L_h}}$. \fyr{Then,} $\gamma \leq  \tfrac{1}{L_h + {\eta} L_f}$ and 
${ (\eta\textstyle\sum_{j=0}^{K-1}\theta_j ) }^{-1} \leq   \tfrac{{\gamma \mu_f }}{{(p+1)\ln(K)}K^{p}}.$
  \end{lemma}
 \begin{proof}
 {{(i)}} The equation follows directly from the update rule $\theta_{k+1}=\tfrac{\theta_k}{(1 - {\eta_{k+1}\gamma{\mu_f}  })}$ for all $k \geq 1$, where $\theta_0 = \tfrac{1}{(1-\eta_0 \gamma \mu_f)}$. To show that $\theta_k >1 $, note that from the condition $\gamma \leq  \tfrac{1}{L_h + {\eta_0} L_f}$, we have $\eta_0\gamma \mu_f <1$, and thus $\eta_k\gamma \mu_f <1$ for all $k\geq 0$.

\noindent (ii) From the update rule of $\eta_k$, we have $\eta_0 = \fyr{1/(2\gamma L_f)}$. Together with $\gamma \leq \fyr{1/(2L_h)}$, we obtain $\gamma(L_h + {\eta_0} L_f) \leq 0.5+0.5 =1$. \ch{Thus, $\gamma \leq \tfrac{1}{\eta_0 L_f + L_h}$.} Next, we show $\theta_k =\tfrac{\eta_{0,l}+k}{\eta_{0,l}-1}
$. Using $\eta_k: = \tfrac{\eta_{0,u}}{\eta_{0,l}+k}$, where $\eta_{0,u} =(\gamma \mu_f)^{-1} $ and $\eta_{0,l} = \tfrac{2L_{f}}{\mu_f}$, for any $k\geq 0$ we have $\textstyle\prod_{t=0}^k(1 - {\eta_t\gamma{\mu_f}  }) = \textstyle\prod_{t=0}^k(1 - \tfrac{1}{t+\eta_{0,l}}) = \tfrac{\eta_{0,l}-1}{\eta_{0,l}+k}$. From the equation in (i), we obtain $\theta_k =\tfrac{\eta_{0,l}+k}{\eta_{0,l}-1}
$ for all $k\geq 0$. Next, we show that $ \textstyle\sum_{k=0}^{\infty}\theta_k \eta_k = \infty$. We have
 \begin{align}\label{eq:sum_thetaeta}
 \textstyle\sum_{k=0}^{K-1}\theta_k \eta_k = \textstyle\sum_{k=0}^{K-1} \tfrac{\left(\eta_{0,l}+k\right)}{\left(\eta_{0,l}-1\right)}\tfrac{\eta_{0,u}}{\left(\eta_{0,l}+k\right)}= \textstyle\sum_{k=0}^{K-1}\tfrac{\eta_{0,u}}{\eta_{0,l}-1} = \tfrac{K\eta_{0,u}}{\eta_{0,l}-1} = \tfrac{K}{\gamma (2L_f - \mu_f)}
 \end{align} 
 \ser{and $\textstyle\sum_{k=0}^{K-1}\theta_k \eta_k^2  = \tfrac{\eta_{0,u}^2}{(\eta_{0,l}-1)}\textstyle\sum_{k=0}^{K-1}\tfrac{1}{({\eta_{0,l}+k})} \leq \tfrac{\eta_{0,u}^2}{(\eta_{0,l}-1)} (\fyr{\frac{1}{\eta_{0,l}}}+\ln (\tfrac{K+\eta_{0,l}-1}{\eta_{0,l}})),$}
where the last inequality is implied by invoking \cite[Lemma 9]{yousefian2017smoothing}. 
 
 \noindent (iii) From the condition $\tfrac{K}{\ln(K)} \geq 2(p+1)\tfrac{L_f}{\mu_f}$ and the value of $\eta$, we have $\gamma \eta L_f \leq 0.5$. Thus, from $\gamma \leq 0.5/L_h$, it follows that $\gamma \leq  \tfrac{1}{L_h + {\eta} L_f}$. From the equation in (i), for a constant regularization parameter $\eta$, we have  $\theta_k = {1}/\left({1-\eta \gamma \mu_f}\right)^ {k+1}$ for all $k\geq 0$. Note that trivially, we can write $  \textstyle\sum_{j=0}^{K-1}\theta_j >  \theta_{K-1} =  (1 - {\eta\gamma{\mu_f}  })^{-K} $. Thus,
$
{ (\textstyle\sum_{j=0}^{K-1}\theta_j ) }^{-1} \leq (1 - {\eta\gamma{\mu_f}  })^{K}.
$
From $\eta=\tfrac{(p+1)\ln(K)}{\gamma \mu_f K}$, we have $(1-\eta\gamma\mu_f)^K = (1-\tfrac{(p+1)\ln(K)}{K})^K$. 
\fyr{Invoking the identity $1-x \leq \exp(-x) $ for $x \in \mathbb{R}$, we} obtain 
 \begin{align*}
(\textstyle\sum_{j=0}^{K-1}\theta_j)^{-1} \leq (1-\eta\gamma\mu_f)^K = (1-\tfrac{(p+1)\ln(K)}{K})^K
 \leq  \exp(-(p+1)\ln(K)) = \tfrac{1}{K^{p+1}}.
\end{align*}
Thus, from the value of $\eta$, we obtain \ser{$
1/{ (\eta\textstyle\sum_{j=0}^{K-1}\theta_j ) } \leq   {{(\gamma \mu_f )}}/{({(p+1)\ln(K)}K^{p})}.
$}
 \end{proof}
\fyr{Next, we show that the sequence $\{\bar{x}_{k}\}$ in Algorithm~\ref{alg:IR-ISTA-s} admits a weighted-average representation}.
\begin{lemma}\label{lemma:weighted}\em

\fyr{Consider} Algorithm~\ref{alg:IR-ISTA-s}. \fyr{For} any $K \geq 1$, we have
\begin{align*}
\ssg{\bar{x}_{K} = {(\textstyle\sum_{k=0}^{K-1}\theta_k \eta_k x_{k+1})}/{(\textstyle\sum_{j=0}^{K-1}\theta_j \eta_j)}.}
\end{align*}

\end{lemma}
\begin{proof}
 \fyr{We proceed by induction on} $K \geq 1$. Let $K=1$. From Algorithm~\ref{alg:IR-ISTA-s}, we have  
$
\bar{x}_{1} = \tfrac{\Gamma_0 \bar{x}_{0} + \eta_0\theta_{0} x_{1}}{\Gamma_{1}}.
$
Since $\Gamma_0 = 0$ and $\Gamma_1 = \eta_0\theta_0$, this simplifies to  
$
\bar{x}_{1} = \tfrac{\eta_0\theta_{0} x_{1}}{\eta_0 \theta_0} =  ({\textstyle\sum_{k=0}^{0}\theta_k \eta_k x_{k+1}})/{\textstyle\sum_{j=0}^{0}\theta_j \eta_j}.
$
Thus, the base case \fyr{holds}. Suppose  
$
\bar{x}_{K} = ({\textstyle\sum_{k=0}^{K-1}\theta_k \eta_k x_{k+1}})/{\textstyle\sum_{j=0}^{K-1}\theta_j \eta_j}$, for some $K \geq 1$. \fyr{Using this induction hypothesis in the update}  $
\bar{x}_{K+1} = \tfrac{\Gamma_K \bar{x}_{K} + \eta_K\theta_{K} x_{K+1}}{\Gamma_{K+1}}
$ \fyr{in} Algorithm~\ref{alg:IR-ISTA-s} \fyr{and} recalling that $\Gamma_{K+1} ={\textstyle\sum_{j=0}^{K}\theta_j \eta_j}  $, we obtain \ser{$\bar{x}_{K+1} 
= {(\textstyle\sum_{k=0}^{K}\theta_k \eta_k x_{k+1})}/{(\textstyle\sum_{j=0}^{K}\theta_j \eta_j)}.$} \fyr{This completes the proof.}
\end{proof}
In the following, we derive a preliminary result that will be utilized in the analysis.
\begin{lemma}\label{lemma:lower_f_strongly}\em Consider the problem in~\eqref{prob:uni_centr} under Assumption~\ref{assump:strongly_main}. Let $x^*$ be the unique optimal solution to this problem. Then, the following \fyr{properties} hold.

\ser{\noindent (i)} For any $x \in \mathbb{R}^n$ and \fyr{any} $\tilde{\nabla} \bar{f}(x^*) \in \partial \bar{f}(x^*)$, we have
\begin{align}\label{eq:lower_f_strongly}
&-\|\tilde{\nabla} \bar{f}(x^*)\|\,  \operatorname{dist}(x,X^*_{\bar{h}})+\tfrac{\mu_f}{2} {\|x - x^*\|}^2 \leq \bar{f}(x) - \bar{f}^*.
\end{align}
\ser{\noindent (ii)} If $ X^*_{\bar h}$ is $\alpha$-weak sharp minima of order $m \geq 1$. Then, for any $x \in \mathbb{R}^n$,
\begin{align}\label{prob:distance-weak-strongly}
&{\|{x} - x^*\|}^2 \leq    \tfrac{2}{\mu_f} \left(\bar{f}({x}) -\bar{f}^* 
+\|\tilde{\nabla} \bar{f}(x^*)\| \sqrt[m]{\alpha^{-1} (\bar{h}({x}) - \bar{h}^*)}\right).
\end{align}
\end{lemma}
\begin{proof}
{{(i)}} Under Assumption~\ref{assump:strongly_main}~(vi), the set $\partial \bar{f}(x^*)$ is nonempty and bounded~\cite[Theorems 3.14 and 3.18]{beck2017first}.
Let us define $\hat{x} := \Pi_{X^*_{\bar{h}}}[x] \in X^*_{\bar{h}}$. By invoking the optimality condition on the problem \ch{in}~\eqref{prob:uni_centr} and  \fyr{using} Lemma~\ref{lemma:convexity_X_h}, we obtain 
$
\langle \tilde{\nabla} \bar{f}(x^*), \hat{x} - x^* \rangle \geq 0$,
where $\tilde{\nabla} \bar{f}(x^*) \in \partial \bar{f}(x^*)$\ch{~\cite[Corollary 3.68]{beck2017first}}. Then, by using the strong convexity property of $\bar{f}$ and \ser{the} Cauchy-Schwarz inequality, we obtain the desired lower bound for the suboptimality as 
\begin{align*}
\ser{\bar{f}(x) - \bar{f}^*}
&\ser{\geq
\langle \tilde{\nabla} \bar{f}(x^*),x-\hat{x}\rangle
+\langle \tilde{\nabla} \bar{f}(x^*),\hat{x}-x^*\rangle
+\tfrac{\mu_f \|x-x^*\|^2}{2}}
\\
&\ssg{\geq}
\ser{-\|\tilde{\nabla} \bar{f}(x^*)\|\,
\operatorname{dist}(x,X^*_{\bar{h}})
+\tfrac{\mu_f \|x-x^*\|^2}{2}.}
\end{align*}

\noindent {{(ii)}} \ser{Rearranging \eqref{eq:lower_f_strongly} gives
$\tfrac{\mu_f\|x-x^*\|^2}{2}
\leq \bar f(x)-\bar f^*
+\|\tilde{\nabla}\bar f(x^*)\|\,\operatorname{dist}(x,X^*_{\bar h})$.
The result then follows from the weak sharp minimality of $X^*_{\bar h}$ in Definition~\ref{def:weaksharp}.}
\end{proof}
Next, we \fyr{establish an auxiliary result used in the convergence analysis.}
\ser{\begin{lemma}\label{lem:RAPM2}\em
Consider Definition~\ref{def:uni_funcs} and suppose  Assumption~\ref{assump:strongly_main} holds. Let $\{\eta_k\}$ be a positive nonincreasing sequence. 

\noindent  (i) \fyr{Let $\gamma_k=\frac{1}{L_{\eta_k}}$ for any $k\geq0$. Then, for any $x,y\in\mathbb{R}^n$ and any $k\ge 0$, we have} 
\begin{align} \label{eq:lemma_part_i}
   \bar g_{\eta_k}(x)-\bar g_{\eta_k}(q_{\eta_k}(y))
&\geq
\tfrac{L_{\eta_k}}{2}\|q_{\eta_k}(y)-x\|^2
-\tfrac{L_{\eta_k}-\eta_k\mu_f}{2}\|y-x\|^2.
\end{align} 

\noindent  (ii) Let $\fyr{0<\gamma \le {1}/{L_{\eta_0}}}$. Then, for any $x,y\in\mathbb{R}^n$ and any \fyr{$k\ge 0$}, we have 
\begin{align}\label{lem1boundF_stronger}
\bar{g}_{\eta_k}(x)-\bar{g}_{\eta_k}(q_{\eta_k}(y))
\ge \tfrac{1}{2\gamma}\|x-q_{\eta_k}(y)\|^2
- (\tfrac{1}{2\gamma}-\tfrac{\eta_k\mu_f}{2})\|x-y\|^2.
\end{align}
\end{lemma}
\begin{proof} (i) \fyr{Fix $k\geq 0$. Let $q_{\eta_k}(y)$ be} defined as in Definition~\ref{def:uni_funcs}, with \fyr{$\gamma:=\frac{1}{L_{\eta_k}}$}. By Assumption~\ref{assump:strongly_main} and
\fyr{Definition~\ref{def:uni_L}}, the function
$g_{\eta_k}$ is $L_{\eta_k}$-smooth. Therefore, for any
$y\in\mathbb{R}^n$, we have
\begin{align}\label{eq:smooth_part_i}
g_{\eta_k}(q_{\eta_k}(y))
&\leq
g_{\eta_k}(y)
+\langle q_{\eta_k}(y)-y,\nabla g_{\eta_k}(y)\rangle
+\tfrac{L_{\eta_k}}{2}\|q_{\eta_k}(y)-y\|^2.
\end{align}
Also, since $g_{\eta_k}(x)=\eta_k f(x)+h(x)$ and $f$ is
$\mu_f$-strongly convex, the function $g_{\eta_k}$ is
$\eta_k\mu_f$-strongly convex. Thus, for any
$x,y\in\mathbb{R}^n$, we have
\begin{align}\label{eq:convex_part_i}
g_{\eta_k}(y)
&\leq
g_{\eta_k}(x)
-\langle x-y,\nabla g_{\eta_k}(y)\rangle
-\tfrac{\eta_k\mu_f}{2}\|x-y\|^2.
\end{align}
Summing \eqref{eq:smooth_part_i} and
\eqref{eq:convex_part_i}, we obtain
\begin{align}\label{eqn:smooth_f_eta_1_part_i}
&
\langle x-q_{\eta_k}(y),\nabla g_{\eta_k}(y)\rangle
-\tfrac{L_{\eta_k}}{2}\|q_{\eta_k}(y)-y\|^2
+\tfrac{\eta_k\mu_f}{2}\|x-y\|^2
\leq
g_{\eta_k}(x)-g_{\eta_k}(q_{\eta_k}(y)).
\end{align}
\ssg{Applying Lemma~\ref{lem:RAPM1} and  Definition~\ref{def:uni_funcs},} \ssg{yields} $\tfrac{y-q_{\eta_k}(y)}{\gamma_k}
-\nabla g_{\eta_k}(y)
\in
\partial\omega_{\eta_k}(q_{\eta_k}(y)).$
Therefore, by the definition of the subgradient of
$\omega_{\eta_k}$, we obtain
\begin{align}\label{eqn:smooth_f_eta_2_part_i}
\langle
x-q_{\eta_k}(y),
\fyr{L_{\eta_k}(y-q_{\eta_k}(y))} 
-\nabla g_{\eta_k}(y)
\rangle
\leq
\omega_{\eta_k}(x)
-\omega_{\eta_k}(q_{\eta_k}(y)).
\end{align}
Summing \eqref{eqn:smooth_f_eta_1_part_i} and
\eqref{eqn:smooth_f_eta_2_part_i}, \fyr{we obtain}
\begin{align*} 
&
\fyr{L_{\eta_k}}
\langle x-q_{\eta_k}(y),y-q_{\eta_k}(y)\rangle
-\fyr{\tfrac{L_{\eta_k} \ssg{\|q_{\eta_k}(y)-y\|^2}}{2}}
+\tfrac{\eta_k\mu_f \ssg{\|x-y\|^2}}{2}
\leq
\bar g_{\eta_k}(x)
-\bar g_{\eta_k}(q_{\eta_k}(y)).
\end{align*}
Substituting
$\langle x-q_{\eta_k}(y),y-q_{\eta_k}(y)\rangle
=
\tfrac{1}{2}\|x-q_{\eta_k}(y)\|^2
+\tfrac{1}{2}\|y-q_{\eta_k}(y)\|^2
-\tfrac{1}{2}\|x-y\|^2$, \fyr{we obtain \eqref{eq:lemma_part_i}}.

\noindent (ii) Following the same argument as in part~(i),  using the constant stepsize $\gamma \leq 1/L_{\eta_0}$, and noting that the nonincreasing property of
$\{\eta_k\}$ implies that $g_{\eta_k}$ is $L_{\eta_0}$-smooth for all $k$, we obtain \eqref{lem1boundF_stronger}.
\end{proof}}

\begin{lemma}\label{lemma:3.9}\em
Consider the problem in~\eqref{prob:uni_centr} \fyr{under} Assumption~\ref{assump:strongly_main}. Let $x^*$ \fyr{denote its} unique optimal solution, and let $\{x_k\}$ be \fyr{the} sequence generated by Algorithm~\ref{alg:IR-ISTA-s}. \ser{Let $\{\eta_k\}$ be a positive nonincreasing sequence, and let $\gamma>0$ satisfy $\gamma \le \tfrac{1}{L_{\eta_0}}$.} Then, for any $k \geq 0$, 
\ser{\begin{align}\label{eq:ir-ista-barh}
{\eta_k}(\bar{f}(x_{k+1})-\bar{f}^*)+ \bar{h}(x_{k+1}) - \bar{h}^* &\leq
(\tfrac{1}{2\gamma}-\tfrac{\eta_k\mu_f}{2})\|x_k-x^*\|^2 - \tfrac{\|x_{k+1}-x^*\|^2}{2\gamma}.
\end{align}}
Further, for any $K\geq 1$, 
 \begin{align}\label{eq:beforedroppping}
& \textstyle\sum_{k=0}^{K-1}\theta_k \eta_k (\bar{f}(x_{k+1})-\bar{f}^*)\leq  \tfrac{1}{2\gamma }( {\| x_0 - x^*\|}^2 -{\theta_{K-1}}{\|x_K -x^*\|}^2). 
\end{align}
\end{lemma}
\begin{proof} 
\ser{Note that, from Definition~\ref{def:uni_funcs} and Algorithm~\ref{alg:IR-ISTA-s}, we have   $q_{\eta_k} (x_k) = x_{k+1} $. \fyr{Setting $y=x_k$ and $x= x^*$} in   \eqref{lem1boundF_stronger}, we obtain the result in \eqref{eq:ir-ista-barh}.}
Next, we show the inequality \eqref{eq:beforedroppping}. \fyr{Dropping} the nonnegative term $ \bar{h}(x_{k+1}) - \bar{h}^* $ from the left-hand side of \eqref{eq:ir-ista-barh}, 
we obtain
 \begin{align}\label{eq:irs-ista-s-f}
& {\eta_k} (\bar{f}(x_{k+1})-\bar{f}^*)
\leq  (\tfrac{1}{2\gamma })( \left(1 - {\eta_k\gamma{\mu_f}  }\right){\|x_k -x^*\|}^2  -  {\|x_{k+1} - x^*\|}^2).
\end{align}
 Now, consider \eqref{eq:irs-ista-s-f} for $k=0$. \fyr{Multiplying}  \ser{both sides} by $\theta_0$, we have
 \begin{align}\label{eq:irs-ista-s-f-k=0-sum}
&  \theta_0 \eta_0(\bar{f}(x_{1})-\bar{f}^*)\leq   ( \tfrac{1}{2\gamma })( {\|x_0 -x^*\|}^2  -  \theta_0{\|x_{1} - x^*\|}^2).
\end{align}
\noindent Considering \eqref{eq:irs-ista-s-f}  and multiplying  both sides by $\theta_k=1/\textstyle\prod_{t=0}^k(1 - {\eta_t\gamma{\mu_f}  })$, and recalling the nonincreasing property of $\{ \eta_k\}$, we have for any $k\geq 1$\fyr{,} 
 \begin{align*}
& \theta_k \eta_k (\bar{f}(x_{k+1})-\bar{f}^*)\leq   (\tfrac{1}{2\gamma })({ \theta_{k-1}}{\|x_k -x^*\|}^2  - {\theta_k} {\|x_{k+1} - x^*\|}^2),
\end{align*}
and  \ser{summing}  both sides over $k = 1, 2, \ldots, K-1$, yields
 \begin{align}\label{eq:irs-ista-s-f-k>0-sum}
&\textstyle\sum_{k=1}^{K-1}\theta_k \eta_k (\bar{f}(x_{k+1})-\bar{f}^*)\leq  (\tfrac{1}{2\gamma } ) ({\theta_0} {\| x_1 - x^*\|}^2 - {\theta_{K-1}}{\|x_K -x^*\|}^2).
\end{align}
\noindent Finally, \fyr{summing} \eqref{eq:irs-ista-s-f-k=0-sum} and \eqref{eq:irs-ista-s-f-k>0-sum}, we obtain \eqref{eq:beforedroppping}.
\end{proof}
In the following result, we establish suboptimality and infeasibility error bounds and provide an asymptotic convergence guarantee for \ch{$\text{IR-ISTA}_\text{s}$.}
\begin{theorem}\label{proposition:F_eta_IR-ISTA-S}\em
Consider the problem \ch{in}~\eqref{prob:uni_centr} under Assumption \ref{assump:strongly_main}. Let $x^*$ \fyr{denote its  unique optimal solution}.  
\ch{Let $\{\bar{x}_k\}$ be a sequence generated by Algorithm~\ref{alg:IR-ISTA-s}} where $\gamma \leq  \tfrac{1}{L_h + {\eta_0} L_f}$. \ser{Let $\{\eta_k\}$ be a nonincreasing sequence.} 

\ser{\noindent (i)} [suboptimality  bounds] For any $K \geq 1$, 
 \begin{align*}
 -\|\tilde{\nabla} \bar{f}(x^*)\|  \operatorname{dist}(\bar{x}_K,X^*_{\bar{h}})+\tfrac{\mu_f {\|\bar{x}_K - x^*\|}^2}{2}  & \leq \bar{f}(\bar{x}_{K})-\bar{f}^* \leq  \| x_0 - x^*\|^2\ser{/(2\gamma{\textstyle\sum_{j=0}^{K-1}\theta_j \eta_j})}. 
\end{align*}
\ser{\noindent (ii)} [infeasibility  bounds] 
For any $K \geq 1$,
\begin{align*}
& 0 \leq \bar{h}(\bar{x}_{K}) - \bar{h}^* \leq  {(\tfrac{ \eta_0  {\|x_0 - x^*\|}^2 }{2\gamma} 
+ {(\bar{f}^* - \hat{C}_{\bar{f}})} \textstyle\sum_{j=0}^{K-1}\eta_j^2 \theta_j  )}/ {\textstyle\sum_{j=0}^{K-1} \eta_j\theta_j}.
\end{align*}
\ser{\noindent (iii)} [asymptotic convergence] \fyr{If $\liminf_{k \to \infty}{\textstyle\sum_{j=0}^{k-1}\theta_j \eta_j} >0$, then $\{\bar{x}_k\}$ is bounded. Suppose $\lim_{k\to \infty}\eta_k =0$,} \ser{$\textstyle\sum_{j=0}^{\infty} \eta_j\theta_j = \infty$, and} $\lim_{k\to \infty}(\textstyle\sum_{j=0}^{k}\eta_j^2 \theta_j /\textstyle\sum_{j=0}^{k} \eta_j\theta_j) = 0$  (e.g., see Lemma~\ref{lemma:defining_theta}~(ii)).  Then, $\{\bar{x}_k\}$ has a limit point and 
 $\lim_{k \to \infty} \bar{x}_k =x^*$\ser{.}
\end{theorem}
\begin{proof}
\noindent {{(i)}} The lower bound holds due to \eqref{eq:lower_f_strongly}. \fyr{To establish the upper bound, we proceed as follows.}  
\fyr{Dropping} the nonnegative term ${\theta_{K-1}}{\|x_K -x^*\|}^2$ from the right-hand side of \eqref{eq:beforedroppping}\fyr{,} 
dividing \ser{both sides} by $\textstyle\sum_{j=0}^{K-1}\theta_j \eta_j$,  \fyr{and applying} Jensen's inequality  \fyr{together with} invoking Lemma~\ref{lemma:weighted}, we obtain the result.

\noindent {{(ii)}} Note that $\bar{h}^*:= \inf_{x\in \mathbb{R}^n} \bar{h}(x)$ implies that the lower bound for infeasibility is zero. Then, \fyr{considering} \eqref{eq:ir-ista-barh}, 
 we have
 \begin{align}\label{eq:ir0ista-h-first}
 \bar{h}(x_{k+1}) - \bar{h}^* \leq  \tfrac{ 
    \left(1 - {\eta_k \gamma\mu_f}\right) {\|x_k - x^*\|}^2 
    - {\|x_{k+1} - x^*\|}^2 
}{2\gamma}  
+ \eta_k  (\bar{f}^* - \hat{C}_{\bar{f}}).
\end{align}
Consider the preceding inequality for $k=0$. \fyr{Multiplying} \ser{both sides} by $\eta_0\theta_0 = \eta_0/(1 - {\eta_0\gamma{\mu_f}  })$, and recalling that $\eta_1 \leq \eta_0$, we have
 \begin{align}\label{eq:ir-ista-h-k=0}
& \eta_0\theta_0( \bar{h}(x_{1}) - \bar{h}^*  )  \leq  \tfrac{ 
    \eta_0{\|x_0 - x^*\|}^2
    - \eta_1\theta_0{\|x_{1} - x^*\|}^2 
 }{2\gamma}  
+ \eta_0^2 \theta_0  (\bar{f}^* - \hat{C}_{\bar{f}}).
\end{align}
\noindent \fyr{Multiplying} \ser{both sides} of \eqref{eq:ir0ista-h-first} by $\eta_k\theta_k=\eta_k/\textstyle\prod_{t=0}^k(1 - {\eta_t{\gamma\mu_f}  })$\ser{,} using $\eta_{k+1} \leq \eta_k$, \ser{and then} summing \ser{both sides} \fyr{of the resulting inequality}
over $k=1, \ldots, K-1$,  
\begin{align}\label{eq:suetakmu_f}
  \textstyle\sum_{k=1}^{K-1} \eta_k\theta_k(\bar{h}(x_{k+1}) - \bar{h}^*) & \leq \tfrac{\left( 
    \eta_1\theta_{0} {\|x_1 - x^*\|}^2
    - \eta_{K}\theta_{K-1}{\|x_{K} - x^*\|}^2
\right) }{2\gamma}  
\\\notag & +  \textstyle\sum_{k=1}^{K-1}\eta_k^2 \theta_k (\bar{f}^* - \hat{C}_{\bar{f}}).
\end{align}
\fyr{Summing} \eqref{eq:ir-ista-h-k=0} and \eqref{eq:suetakmu_f}, and dropping 
\ch{$ \eta_{K}\theta_{K-1}{\|x_{K} - x^*\|}^2 \geq 0$} from the right-hand side, we obtain
\ser{$\textstyle\sum_{k=0}^{K-1} \eta_k\theta_k(\bar{h}(x_{k+1}) - \bar{h}^*)\leq  \tfrac{ 
    \eta_0 {\|x_0 - x^*\|}^2
}{2\gamma}  
+  \textstyle\sum_{k=0}^{K-1}\eta_k^2 \theta_k  (\bar{f}^* - \hat{C}_{\bar{f}}).$}
\fyr{Dividing} \ser{both sides} of the preceding inequality by $\textstyle\sum_{j=0}^{K-1} \eta_j\theta_j$ and \fyr{applying}  Jensen's inequality and Lemma~\ref{lemma:weighted}, we obtain the result.

\noindent {{(iii)}} 
\fyr{First, we show the boundedness of $\{\bar{x}_k\}$. From (i), together with $\operatorname{dist}(\bar{x}_K,X^*_{\bar{h}}) \leq \|\bar{x}_K-x^*\|$, we have 
 $ \tfrac{\mu_f}{2} \|\bar{x}_K - x^*\|^2  \leq  \|\tilde{\nabla} \bar{f}(x^*)\|\|\bar{x}_K-x^*\|+ \| x_0 - x^*\|^2\ser{/(2\gamma{\textstyle\sum_{j=0}^{K-1}\theta_j \eta_j})}.$ This quadratic inequality, together with assumption $\liminf_{k \to \infty}{\textstyle\sum_{j=0}^{k-1}\theta_j \eta_j} >0$, \ssg{implies} that $\{\bar{x}_k\}$ is bounded. Next, we show the asymptotic convergence claim.} \ser{Theorem~\ref{proposition:F_eta_IR-ISTA-S}~(ii), together with
$\lim_{k\to \infty}\textstyle\sum_{j=0}^{k}\eta_j^2\theta_j/
\textstyle\sum_{j=0}^{k}\eta_j\theta_j=0$ and
$\textstyle\sum_{j=0}^{\infty}\eta_j\theta_j=\infty$, \fyr{imply} that
$\lim_{k\to\infty}\bar h(\bar{x}_k)=\bar h^*$.} \fyr{Moreover, since
$\{\bar{x}_k\}$ is bounded, it has at least one convergent subsequence. Let $\{\bar{x}_{k_j}\}$ denote an arbitrary convergent subsequence with
$\lim_{j\to\infty}\bar{x}_{k_{j}}=\hat{x}$. By Assumption~\ref{assump:strongly_main}, $\bar h$ is lower
semicontinuous. Thus,  $\bar h(\hat{x})\leq
\liminf_{j\to\infty} \bar h(\bar{x}_{k_{j}})=\bar h^*$. We also have $\bar h^*\leq \bar h(\hat{x})$, implying that $\hat{x} \in X^*_{\bar{h}}$. Since $\operatorname{dist}(\bar{x}_K,X^*_{\bar{h}}) \leq \|\bar{x}_K-\hat{x}\|$, part (i) implies that $ \tfrac{\mu_f {\|\bar{x}_K - x^*\|}^2}{2}   \leq    {\| x_0 - x^*\|}^2\ser{/{(2\gamma{\textstyle\sum_{j=0}^{K-1}\theta_j \eta_j})}} +\|\tilde{\nabla} \bar{f}(x^*)\| \|\bar{x}_K-\hat{x}\|.$ Taking the limit on both sides along $\{\bar{x}_{k_{j}}\}$, we obtain $x^* = \hat{x}$. Since $\{\bar{x}_{k_j}\}$ is an arbitrary convergent subsequence, we conclude that $\lim_{k\to\infty}\bar{x}_k=x^*$.}
 \end{proof}
In the following result, we provide both the asymptotic guarantee and  nonasymptotic rate statements under a prescribed diminishing regularization sequence.
\begin{corollary} \label{theorem:ir-ista-s}\em
Consider the problem \ch{in}~\eqref{prob:uni_centr} under Assumption~\ref{assump:strongly_main}.  \ch{Let $x^*$ \fyr{denote its unique optimal solution,  and let} $\{\bar{x}_k \}$ be generated by Algorithm~\ref{alg:IR-ISTA-s}.} \ch{Let $\gamma \leq \fyr{\tfrac{1}{2L_h}}$ and $\eta_k: = \tfrac{\eta_{0,u}}{\eta_{0,l}+k}$, where $\eta_{0,u} =(\gamma \mu_f)^{-1} $ and $\eta_{0,l} = \tfrac{2L_{f}}{\mu_f}$.}    \fyr{Then,} $\lim_{{k} \to \infty} \bar{x}_{k}  = x^*$. Additionally, the following results hold for any $K\geq 1$.\\
\ser{\noindent (i)} {{[suboptimality  \fyr{upper bound and a lower bound estimate}]}}   For any $\tilde{\nabla} \bar{f}(x^*) \in \partial \bar{f}(x^*)$, we have \ser{$-\|\tilde{\nabla} \bar{f}(x^*)\|\, \operatorname{dist}(\bar{x}_K,X^*_{\bar{h}})+\tfrac{\mu_f}{2} {\|\bar{x}_K - x^* \|}^2 \leq \bar{f}(\bar{x}_{K})-\bar{f}^* \leq  \tfrac{ u_1}{K}
$}\ssg{, where  $u_1 := \fyr{\tfrac{1}{2}} {\| x_0 - x^*\|}^2(2L_f - \mu_f)$.}\\
\ser{\noindent (ii)} {{[infeasibility  bounds]}}  We have
 $0\leq  \bar{h}(\bar{x}_{K}) - \bar{h}^*  \leq \tfrac{u_{2,K}}{K}$, where \\ $
u_{2,K}  :=  \gamma (2L_f - \mu_f) (\tfrac{ \eta_0   {\| x_0 - x^*\|}^2 }{2\gamma}
+  (\tfrac{(\bar{f}^* - \hat{C}_{\bar{f}} )\eta_{0,u}^2}{\eta_{0,l}-1} (\fyr{\frac{1}{\eta_{0,l}}}+\ln (\tfrac{K+\eta_{0,l}-1}{\eta_{0,l}}))  ) ).$\\
\ser{\noindent (iii)} {{[distance to the unique optimal solution]}}  Suppose $X^*_{\bar{h}}$ is $\alpha$-weak sharp minima of order $m \geq 1$. Then, ${\|\bar{x}_{K} - x^*\|}^2  \leq   \tfrac{2 u_1}{\mu_f K}
+ \tfrac{2 \|\tilde{\nabla} \bar{f}(x^*)\|}{\mu_f} \sqrt[m]{\tfrac{u_{2,K}}{\alpha K}}$.\\
\ser{\noindent (iv)} {{[suboptimality lower bound]}} 
If $X^*_{\bar{h}}$ is $\alpha$-weak sharp minima of order $m \geq 1$, then for any $\tilde{\nabla} \bar{f}(x^*) \in \partial \bar{f}(x^*)$, 
$ -\|\tilde{\nabla} \bar{f}(x^*)\| \sqrt[m]{\tfrac{u_{2,K}}{\alpha K}} \leq \bar{f}(\bar{x}_K) - \bar{f}^*.$
\end{corollary}
\begin{proof} 
\noindent The asymptotic convergence result \fyr{follows from} Lemma~\ref{lemma:defining_theta} (ii) \fyr{together with} Theorem~\ref{proposition:F_eta_IR-ISTA-S}~(iii).

\noindent (i) Consider Theorem~\ref{proposition:F_eta_IR-ISTA-S}~(i). The result follows by invoking Lemma~\ref{lemma:defining_theta}~(ii).  

\noindent (ii) Consider Theorem~\ref{proposition:F_eta_IR-ISTA-S}~(ii). The result follows by invoking Lemma~\ref{lemma:defining_theta}~(ii).

 \noindent (iii) Consider \eqref{prob:distance-weak-strongly}. The result \fyr{follows from} the bounds in {(i)} and {(ii)}.
 
\noindent (iv)  Consider \eqref{eq:lower_f_strongly} \fyr{where we set $x=\bar{x}_K$}. Dropping the nonnegative term $\tfrac{\mu_f {\ssg{\|\fyr{\bar{x}_K} - x^* \|}^2}}{2}$ from the left-hand side and using Definition~\ref{def:weaksharp}, we obtain the result.
\end{proof}
\ch{Next, we provide rate statements under 
a constant regularization parameter.}
\begin{corollary} \label{cor:r-ista-s}\em
Consider the problem \ch{in}~\eqref{prob:uni_centr} under Assumption~\ref{assump:strongly_main}.  Let $x^*$ \fyr{denote its unique optimal solution, and let} $\bar{x}_K$ be generated by Algorithm~\ref{alg:IR-ISTA-s}.

 \noindent (1) Let \fyr{$\gamma \leq \tfrac{1}{2L_h}$} and  $\eta\ser{=}\tfrac{(p+1)\ln(K)}{\gamma \mu_f K}$ for some arbitrary $ p > 0$ and $K>1$ such that $\tfrac{K}{\ln(K)} \geq 2(p+1)\tfrac{L_f}{\mu_f}$.  Then, the following \fyr{statements}  hold.\\
 \noindent {{(1.i) [suboptimality  bounds]}} Let  $u_3 := \tfrac{ {\| x_0 - x^*\|}^2 \mu_f } { \fyr{2(p+1)}}$. Then, for any $\tilde{\nabla} \bar{f}(x^*) \in \partial \bar{f}(x^*)$, \ser{we have $-\|\tilde{\nabla} \bar{f}(x^*)\|  \operatorname{dist}(\bar{x}_K,X^*_{\bar{h}})+\tfrac{\mu_f}{2} {\|\bar{x}_K - x^* \|}^2 \leq \bar{f}(\bar{x}_{K})-\bar{f}^* \leq \tfrac{u_3}{\ln(K) K^{p}}\ser{.}$}

\noindent {{(1.ii) [infeasibility  bounds]}} Let $u_4 := \tfrac{   {\| x_0 - x^*\|}^2 }{2\gamma}$ and $u_5 := \tfrac{(p+1){\left(\bar{f}^* - \hat{C}_{\bar{f}}\right)}}{\gamma \mu_f}$. We have\\ \ser{$0 \leq \bar{h}(\bar{x}_{K}) - \bar{h}^* \leq  \tfrac{u_4}{K^{(p+1)}}  +  \tfrac{u_5\ln(K)}{ K}.$}

\noindent {{\ser{(1.iii)} [distance to optimal solution]}}  Suppose  $X^*_{\bar{h}}$ is $\alpha$-weak sharp minima of order $m \geq 1$. Then, ${\|\bar{x}_{K} - x^*\|}^2 \leq \tfrac{2 u_3}{\mu_f   \ln(K) K^{p} } + \tfrac{2 \|\tilde{\nabla} \bar{f}(x^*)\|}{\mu_f}   \sqrt[m]{\tfrac{u_4}{\alpha K^{(p+1)}}  +  \tfrac{u_5\ln(K)}{\alpha K}}.$

\noindent {{\ser{(1.iv)} [suboptimality lower bound]}}
Furthermore, if $X^*_{\bar{h}}$ is $\alpha$-weak sharp minima of order $m \geq 1$, then  
$ -\|\tilde{\nabla} \bar{f}(x^*)\| \sqrt[m]{\tfrac{u_4}{\alpha K^{(p+1)}}  +  \tfrac{u_5\ln(K)}{\alpha K}} \leq \bar{f}(\bar{x}_K) - \bar{f}^*.$

\noindent (2)  \ch{Assume that $X^*_{\bar{h}}$ is $\alpha$-weak sharp minima of order $m=1$ and 
$\eta \leq \tfrac{\alpha}{2\|\tilde{\nabla} \bar{f}(x^*)\|} $, for some $\tilde{\nabla} \bar{f}(x^*) \in \partial \bar{f}(x^*)$.  Then, the following results hold for any $K\geq 1$.}

\noindent {{(2.i)} [infeasibility bounds]} We have
$ 0 \leq  {\bar{h}(\bar{x}_{K}) - \bar{h}^*} \leq  \tfrac{  {\| x_0 - x^*\|}^2}{\gamma }(1-\eta\gamma\mu_f)^K.$ 

\noindent {{(2.ii) }[distance to lower-level solution set]}
 $   \operatorname{dist}(\bar{x}_K,X^*_{\bar{h}})   \leq  \tfrac{  {\| x_0 - x^*\|}^2}{\alpha\gamma }(1-\eta\gamma\mu_f)^K.$

\noindent {{(2.iii)} [suboptimality bounds]} \ser{$-\tfrac{{\| x_0-x^* \|}^2 (1-\gamma\eta\mu_f)^K}{2\eta\gamma} \leq \bar{f} (\bar{x}_K) - \bar{f}^*   \leq  \tfrac{  {\| x_0 - x^*\|}^2 (1-\gamma\eta\mu_f)^K}{2\eta\gamma}.$}\\
\noindent {{(2.iv)} [distance to  optimal solution]} We have
$
{\|\bar{x}_{K} - x^*\|}^2 \leq \tfrac{ 2{\| x_0 - x^*\|}^2 }{\eta \gamma\mu_f} (1-\eta\gamma\mu_f)^K.$
\end{corollary}
\begin{proof}
\noindent {{(1.i)}} The lower bound is obtained from \eqref{eq:lower_f_strongly}. The upper bound follows
\fyr{from} Theorem~\ref{proposition:F_eta_IR-ISTA-S}~(i) for a constant $\eta$ \fyr{together with} Lemma~\ref{lemma:defining_theta}~(iii).

\noindent {{(1.ii)}}
 Consider \ch{the relation in} Theorem~\ref{proposition:F_eta_IR-ISTA-S}~(ii) with the constant regularization. \ch{Then, by invoking  Lemma~\ref{lemma:defining_theta}~(iii),  we obtain  the result.}

 \noindent{{\ser{(1.iii)}}}
Consider \eqref{prob:distance-weak-strongly}. Then, we obtain the results by applying (1.i) and (1.ii).

\noindent {{\ser{(1.iv)}}} \fyr{Considering} \eqref{eq:lower_f_strongly}, dropping the nonnegative term from the left-hand side, and using Definition~\ref{def:weaksharp}, we derive the result.

\noindent{{(2.i)}} Note that $\bar{h}^*:= \inf_{x\in \mathbb{R}^n} \bar{h}(x)$ implies that the lower bound for infeasibility is zero. Now,
consider \eqref{eq:lower_f_strongly}. \fyr{Dropping} the nonnegative term $\tfrac{\mu_f}{2} {\|\bar{x}_{k+1} - x^*\|}^2$ from the left-hand side, we obtain
\begin{align}\label{eq:lower-bound-f-ir-ista}
&-\|\tilde{\nabla} \bar{f}(x^*)\|  \operatorname{dist}(\bar{x}_{k+1},X^*_{\bar{h}}) \leq     \bar{f}(\bar{x}_{k+1}) -\bar{f}^*.
\end{align}
\fyr{Considering} \eqref{eq:ir-ista-barh} with constant $\eta$, multiplying \ser{both sides} by \ser{$2\gamma$,} 
we obtain
\begin{align*}
& {2\gamma}\eta\left(\bar{f}(x_{k+1})-\bar{f}^*\right)+ {2\gamma}\left(\bar{h}(x_{k+1}) - \bar{h}^*\right)\leq    (1 - \eta \gamma\mu_f ) {\|x_{\ser{k}} -x^*\|}^2  -  {\|x_{k+1} - x^*\|}^2.
\end{align*}
 \fyr{Invoking} \eqref{eq:lower-bound-f-ir-ista} in the preceding inequality
 and using Definition~\ref{def:weaksharp}, we arrive at
   \begin{align*}
& {2\gamma} (1 - \tfrac{\eta\|\tilde{\nabla} \bar{f}(x^*)\|}{\alpha}  ) ({\bar{h}(x_{k+1}) - \bar{h}^*}) \leq    (1 - \eta \gamma\mu_f ) {\|x_{\ser{k}} -x^*\|}^2  -  {\|x_{k+1} - x^*\|}^2.
\end{align*}
\ser{Multiplying} \ser{both sides}  by $\theta_k = 1/(1-\eta\gamma\mu_f)^{k+1}$, \ser{then} summing \ser{both sides} over $k=0, 1, \ldots, K-1$, for $K\geq 1$, and dropping the nonnegative term $ \theta_{K-1} {{\|x_{\ser{K}} - x^*\|}^2}$ \ser{from} the left-hand side, \ser{yield} $ {2\gamma}  (1 -  {\eta\|\tilde{\nabla} \bar{f}(x^*)\|}/{\alpha}  )\textstyle\sum_{k=0}^{K-1}\theta_k  \left({\bar{h}(x_{k+1}) - \bar{h}^*}\right) \leq { {\|x_0 -x^*\|}^2}.$
\ser{Then} using  Jensen's inequality and Lemma~\ref{lemma:weighted}, we obtain
\begin{align*}
& (\tfrac{2\gamma}{\alpha})  (\alpha - {\eta\|\tilde{\nabla} \bar{f}(x^*)\|}   )(\textstyle\sum_{k=0}^{K-1}\theta_k)\left({\bar{h}(\bar{x}_{K}) - \bar{h}^*}\right) \leq { {\|x_0 -x^*\|}^2}.
\end{align*}
From $\eta \leq \tfrac{\alpha}{2\|\tilde{\nabla} \bar{f}(x^*)\|}$, we have $\alpha \leq  2(\alpha-\eta\|\tilde{\nabla} \bar{f}(x^*)\|)$. Then, in view of $\textstyle\sum_{k=0}^{K-1}\theta_k \geq \theta_{K-1} = 1/(1-\eta\gamma\mu_f)^K$ (cf. Lemma~\ref{lemma:defining_theta}~(i)), and multiplying \ser{both sides} by $\alpha > 0$, we \ch{arrive at the result.}\\
\noindent{{(2.ii)}} The result follows from  (2.i) and Definition~\ref{def:weaksharp}. 

\noindent{{(2.iii)}}
\fyr{Considering} the \ssg{result} in Theorem~\ref{proposition:F_eta_IR-ISTA-S}~(i) for a constant $\eta$ and using $\textstyle\sum_{j=0}^{K-1}\theta_j \geq \theta_{K-1} = 1/(1-\gamma\eta\mu_f)^K$ once again,  we have
 \begin{align}\label{eq:f-upper-after}
& \bar{f}(\bar{x}_{K})-\bar{f}^* \leq \tfrac{  {\| x_0 - x^*\|}^2}{2\eta\gamma}(1-\gamma\eta\mu_f)^K. 
\end{align}
To obtain the lower bound, consider  \eqref{eq:lower-bound-f-ir-ista}. Then, \fyr{using} (2.ii) and $\|\tilde{\nabla} \bar{f}(x^*)\| \leq \tfrac{\alpha}{2 \eta}$, we obtain the lower bound. \\
\noindent{{(2.iv)}}
 Consider  \eqref{eq:lower_f_strongly}. \ser{The} result follows by
 using (2.ii), \eqref{eq:f-upper-after}, and  $\|\tilde{\nabla} \bar{f}(x^*)\| \leq \tfrac{\alpha}{2\eta}$.
\end{proof}

\begin{remark}  Notably, it appears that the asymptotic convergence guarantee in Theorem~\ref{proposition:F_eta_IR-ISTA-S},  exemplified by Corollary~\ref{theorem:ir-ista-s}, is established for the first time in the literature for \ch{addressing} the problem \ch{in}~\eqref{prob:uni_centr}. Additionally, the nonasymptotic convergence rate statements in Corollary~\ref{theorem:ir-ista-s} and Corollary~\ref{cor:r-ista-s} for addressing the SBO problem are novel and are concisely presented in Table~\ref{table:contributions_all}. 
\end{remark}

\subsection{\texorpdfstring{\ser{\ssg{The} \ssg{IR-VFISTA$_{\text{s}}$} method}}{IR-VFISTA-s method}}\label{subsection:3.2}
\ser{\fyr{Here, we aim to leverage acceleration techniques to improve the guarantees for IR-ISTA$_\text{s}$. To this end, we propose} Algorithm~\ref{alg:IR-R-VF_ne} to address the problem \ch{in}~\eqref{prob:uni_centr} under Assumption~\ref{assump:strongly_main}. $\text{IR-VFISTA}_\text{s}$ is an iteratively regularized variant of VFISTA~\cite[Section 10.7.7]{beck2009fast} \fyr{for} single-level composite strongly convex optimization. \fyr{At iteration $k$, a standard proximal step is performed using the regularized mapping $\nabla h  + {\ser{\eta}_k} {\nabla} f$. Further, we use a momentum parameter $\tfrac{\sqrt{\kappa_k} -1}{\sqrt{\kappa_{k+1}}+1}$ where $\kappa_k$ is characterized in terms of $\eta_k$ as $\kappa_k := \tfrac{L_h + {\ser{\eta_k}} L_f}{\ser{\eta_k}{\mu_f}}$. Notably, the momentum parameter utilizes both $\eta_k$ and $\eta_{k+1}$ distinctly.}}
\begin{algorithm}   
 	\caption{\ser{Iteratively} Regularized Variant of  FISTA ($\ser{\text{IR-VFISTA}_\text{s}}$)}
 	\begin{algorithmic}[1]\label{alg:IR-R-VF_ne}
 		\STATE {{\bf input:} $y_0=x_0 \in \mathbb{R}^n$, \ser{positive nonincreasing sequence $\{\eta_k\}$},   $\kappa_k := \tfrac{L_h + {\ser{\eta_k}} L_f}{\ser{\eta_k}{\mu_f}} $,  \ch{$\ser{\gamma_k}  := 1/ \left(L_h + {\ser{\eta_k}} L_f\right)$},  and {$K\geq 1$.}}
 		\FOR {$k = 0, 1, 2, \dots, {K-1}$}
 		\STATE   {$x_{k+1} := \text{prox}_{\ser{\gamma_k} {\omega_{\ser{\eta_k}}}}\left( y_{k} -\gamma_k \left(\nabla h( y_{k}) + {\ser{\eta}_k} {\nabla} f( y_{k})\right)\right)$}	
        \STATE {$y_{k+1} := x_{k+1} + \fyr{\left(\tfrac{\sqrt{\kappa_k} -1}{\sqrt{\kappa_{k+1}}+1}\right)}(x_{k+1}-x_{k})$}
        
 		\ENDFOR
	\STATE {{\bf return:} {${x}_K$}}
 	\end{algorithmic}
 \end{algorithm} 
\ser{
\begin{proposition}\em\label{propositon:IRVFISTA_bound_scaled}
\fyr{Consider Algorithm~\ref{alg:IR-R-VF_ne} under Assumption~\ref{assump:strongly_main}}. \fyr{For $k\geq 0$, define $t_k := \sqrt{\kappa_k} $, and for $k\geq 1$, define the function}\begin{align}\label{eq:IRVFISTA_var_dk}
d_k(x)
&:=\bar g_{\eta_k}(x_k)-\bar g_{\eta_k}(x)+
\tfrac{\eta_k\mu_f}{2}\|
t_{k-1}x_k-
\left(
x+(t_{k-1}-1)x_{k-1}\right)\|^2,
\end{align}
where $d_0(x)
:=
\bar{g}_{\fyr{\eta_0}}(x_0)- \bar{g}_{\eta_0}(x)
+
\tfrac{\eta_0\mu_f}{2}\|x_0-x\|^2.$ 
Then, for all $x\in\mathbb{R}^n$ and  \fyr{$k\geq 0$},
\begin{align}\label{eq:IRVFISTA_var_d_recursion_pro}
d_{k+1}(x)
&\leq
(1-\tfrac{1}{t_k})d_k(x)
+
(
\eta_{k+1}-\eta_k
)
(
\bar f(x_{k+1})-\bar f(x)
).
\end{align}
Moreover, if the regularization parameter is constant, i.e., $\eta_k=\fyr{\eta >0}$ for all $k \geq 0$,
\begin{align}\label{prop:R-VFISTA}
\bar g_{\eta}(x_K)-\bar g_{\eta}(x)
&\leq
(1-\tfrac{1}{\sqrt{\kappa}})^K
\left(\bar g_{\eta}(x_0)-\bar g_{\eta}(x)
+\tfrac{\eta\mu_f}{2}\|x_0-x\|^2
\right).
\end{align}
\end{proposition}
}
\ser{
\begin{proof}
\ssg{Considering $x_{k+1} := \text{prox}_{\ser{\gamma_k} {\omega_{\ser{\eta_k}}}}\left( y_{k} -\gamma_k \left(\nabla h( y_{k}) + {\ser{\eta}_k} {\nabla} f( y_{k})\right)\right)$ and applying  Lemma~\ref{lem:RAPM2} (i),} \fyr{for any $x \in\mathbb{R}^n$ and any $k\geq 0$ we have}
\begin{align}\label{eq:IRVFISTA_var_first}
\bar g_{\eta_k}
(
\tfrac{1}{t_k}x+
(1-\tfrac{1}{t_k})x_k
)
-\bar g_{\eta_k}(x_{k+1})
&\geq
\tfrac{L_{\eta_k}}{2t_k^2}
\|
t_kx_{k+1}
-
\left(
x+(t_k-1)x_k
\right)
\|^2
\\\notag
&-
\tfrac{L_{\eta_k}-\eta_k\mu_f}{2t_k^2}
\|
t_ky_k-
\left(
x+(t_k-1)x_k
\right)
\|^2.
\end{align}
Since $\bar g_{\eta_k}$ is $\eta_k\mu_f$-strongly convex, \fyr{we may write}
\begin{align}\label{eq:IRVFISTA_var_strong_convexity}
&\bar g_{\eta_k}
(
\tfrac{x}{t_k}+
(1-\tfrac{1}{t_k})x_k
)
\leq
\fyr{(\tfrac{1}{t_k})\bar g_{\eta_k}(x)}
+
(1-\tfrac{1}{t_k})
\bar g_{\eta_k}(x_k)
-
\tfrac{\eta_k\mu_f}{2t_k}
(1-\tfrac{1}{t_k})
\|x_k-x\|^2.
\end{align}
Combining \eqref{eq:IRVFISTA_var_first} and
\eqref{eq:IRVFISTA_var_strong_convexity} gives
\begin{align}\label{eq:IRVFISTA_var_combined}
&\bar g_{\eta_k}(x_{k+1})
-
\bar g_{\eta_k}(x)
+\tfrac{L_{\eta_k}}{2t_k^2}
\|t_kx_{k+1}-(x+(t_k-1)x_k)\|^2
\\\notag
&\leq
(1-\tfrac{1}{t_k})
\left(
\bar g_{\eta_k}(x_k)-\bar g_{\eta_k}(x)
\right)
+\tfrac{L_{\eta_k}-\eta_k\mu_f}{2t_k^2}
\|t_ky_k-\left(x+(t_k-1)x_k\right)
\|^2
\\\notag
&-\tfrac{\eta_k\mu_f}{2t_k}
(1-\tfrac{1}{t_k})
\|x_k-x\|^2.
\end{align}
Note that, $t_ky_k-
(
x+(t_k-1)x_k
)
=
x_k-x+t_k(y_k-x_k).$
For all $a,b\in\mathbb{R}^n$ and every $\beta\in[0,1)$, $\|a+b\|^2-\beta\|a\|^2
=
(1-\beta)
\|
a+\tfrac{b}{1-\beta}
\|^2
-
\tfrac{\beta}{1-\beta}\|b\|^2
\leq
(1-\beta)
\|
a+\tfrac{b}{1-\beta}
\|^2.$
Applying the preceding inequality
with
$a
=
x_k-x,
b=
t_k(y_k-x_k),
\beta=
\tfrac{\eta_k\mu_f(t_k-1)}{L_{\eta_k}-\eta_k\mu_f}$
and using $L_{\eta_k}=\eta_k\mu_ft_k^2$, we have $\beta
=
\tfrac{t_k-1}{t_k^2-1}
=
\tfrac{1}{t_k+1},
\
1-\beta
=
\tfrac{t_k}{t_k+1}.$
Therefore, \fyr{for any $k\geq0$,}
\begin{align}\label{eq:IRVFISTA_var_quadratic}
&\tfrac{L_{\eta_k}-\eta_k\mu_f}{2t_k^2}
\|
t_ky_k-
\left(
x+(t_k-1)x_k
\right)
\|^2
-
\tfrac{\eta_k\mu_f}{2t_k}
(1-\tfrac{1}{t_k})
\|x_k-x\|^2
\\\notag
&\leq
(1-\tfrac{1}{t_k})
\tfrac{\eta_k\mu_f}{2}
\|
x_k-x+(t_k+1)(y_k-x_k)
\|^2.
\end{align}
Rearranging the momentum update as $(t_k+1)(y_k-x_k)
=
(t_{k-1}-1)(x_k-x_{k-1})$, yields $x_k-x+(t_k+1)(y_k-x_k)
=
t_{k-1}x_k-
(
x+(t_{k-1}-1)x_{k-1}
).$
Substituting the preceding equality 
into
\eqref{eq:IRVFISTA_var_quadratic} yields, for $k\geq1$,
\begin{align}\label{eq:IRVFISTA_var_quadratic_final}
&\tfrac{L_{\eta_k}-\eta_k\mu_f}{2t_k^2}
\|
t_ky_k-
\left(
x+(t_k-1)x_k
\right)
\|^2
-
\tfrac{\eta_k\mu_f}{2t_k}
(1-\tfrac{1}{t_k})
\|x_k-x\|^2
\\\notag
&\leq
(1-\tfrac{1}{t_k})
\tfrac{\eta_k\mu_f}{2}
\|
t_{k-1}x_k-
\left(
x+(t_{k-1}-1)x_{k-1}
\right)
\|^2.
\end{align}
\fyr{From the preceding inequality, the definition of $d_{k}(x)$, and \eqref{eq:IRVFISTA_var_combined}, we obtain, for any $k\geq 1$,
\begin{align*} 
& \bar g_{\eta_k}(x_{k+1})
-
\bar g_{\eta_k}(x)
+
\tfrac{L_{\eta_k}}{2t_k^2}
\|
t_kx_{k+1}
-
(
x+(t_k-1)x_k
)
\|^2\leq
(1-\tfrac{1}{t_k})d_k(x).
\end{align*}
Using $\tfrac{L_{\eta_k}}{2t_k^2}
=
\tfrac{\eta_k\mu_f}{2}$, the preceding inequality implies that, for any $k\geq1$,
\begin{align*}
& d_{k+1}(x) -(\bar g_{\eta_{k+1}}(x_{k+1})-\bar g_{\eta_{k+1}}(x)) +(\bar g_{\eta_k}(x_{k+1})-\bar g_{\eta_k}(x))\leq
(1-\tfrac{1}{t_k})d_k(x).
\end{align*}
Recall that $\bar{g}_{\eta_k}(\cdot) = \bar h(\cdot) + \eta_k \bar f(\cdot)$. Using this and rearranging the preceding inequality, we arrive at \eqref{eq:IRVFISTA_var_d_recursion_pro} for any $k\geq 1$. Next, we show that \eqref{eq:IRVFISTA_var_d_recursion_pro} holds for $k=0$. Consider \eqref{eq:IRVFISTA_var_quadratic} evaluated at $k=0$. Since $x_0=y_0$, together with \eqref{eq:IRVFISTA_var_combined}, we obtain 
\begin{align}\label{eq:IRVFISTA_var_combined3}
&\bar g_{\eta_0}(x_{1}) - \bar g_{\eta_0}(x) + \tfrac{L_{\eta_0}}{2t_0^2} \| t_0x_{1} - ( x+(t_0- 1)x_0 ) \|^2
\\\notag
&\leq (1-\tfrac{1}{t_0}) \left( \bar g_{\eta_0}(x_0)- \bar g_{\eta_0}(x) \right)
+(1-\tfrac{1}{t_0})\tfrac{\eta_0\mu_f}{2} \| x_0-x\|^2.
\end{align}
Invoking the definition of $d_1(x)$ and $d_0(x)$, together with $\tfrac{L_{\eta_0}}{2t_0^2}
=
\tfrac{\eta_0\mu_f}{2}$, and rearranging the terms in \eqref{eq:IRVFISTA_var_combined3}, we obtain \eqref{eq:IRVFISTA_var_d_recursion_pro} for $k=0$. Next, we consider the constant regularization setting. If} $\eta_k = \eta$ for all $k \geq 0$, then considering \eqref{eq:IRVFISTA_var_dk} and \eqref{eq:IRVFISTA_var_d_recursion_pro}, we obtain $d_{k+1} \leq (1-\tfrac{1}{\sqrt{\kappa}}) d_k$. Applying this inequality recursively gives $d_K(x)
\leq(1-\tfrac{1}{\sqrt{\kappa}}
)^K
d_0(x).$ \fyr{Considering} \eqref{eq:IRVFISTA_var_dk} at $k= K$
and the fact that squared-norm term is nonnegative, we obtain \eqref{prop:R-VFISTA}.
\end{proof}
}
\ser{
\begin{lemma}\label{lemma:general_rho}\em
Let $\{a_k\}$, $\{b_k\}$, $\{d_k\}$, and $\{\theta_k\}$ be
sequences satisfying
\begin{align}\label{def:d_rho}
d_{k+1}
\leq (1-\theta_k)d_k - (\eta_k-\eta_{k+1})b_{k+1},
\qquad \fyr{\hbox{for all }} k\geq0,
\end{align}
where $d_k
:= a_k+\eta_kb_k,
\ k\geq0$. Suppose $a_k \geq 0$ and  $\{\eta_k\}$ is  positive and nonincreasing.  Suppose  there exists a constant $\rho>0$, such that
$0 <\rho\sqrt{\eta_k}\leq
\theta_k
\leq1,$ \fyr{for all $k\geq0$.}

\noindent(i) Suppose  $\{\eta_k\}$ satisfies $\lim_{k \to \infty}\eta_k=0,
\
\sum_{k=0}^{\infty}\sqrt{\eta_k} = \infty,$
and $\lim_{k\to \infty}\tfrac{\eta_k-\eta_{k+1}}
{\eta_{k+1}\sqrt{\eta_k}}
=0.$
Let $[-b_k]_+ \leq B$ for some \fyr{$B \geq 0$}, \fyr{i.e.,} $\{b_k\}$ is bounded below. \fyr{Assume that $a_k\to 0$ implies $[-b_k]_+\to 0$ as $k\to\infty$.} Then, $\lim_{k \to \infty}a_k=0
\ \text{and} \
\lim_{k\to\infty} b_k=0.$

\noindent (ii) Let $\eta_k
:=
\eta_0(k+1)^{-a}$
for some $a\in(0,2)$ and \fyr{$\eta_0>0$}.
\fyr{Let} the coupling condition $b_k
\geq
-C\sqrt[m]{a_k}$
\fyr{hold} for some \fyr{$m>1$ and $C\geq0$}. \fyr{Then,
$0 \leq a_k
\leq \mathcal{O}(k^{-\tfrac{am}{m-1}})$ and $-\mathcal{O}(k^{-\tfrac{a}{m-1}})
\leq b_k\leq\mathcal{O}(k^{-1-\tfrac{a(3-m)}{2(m-1)}}
).$}
\end{lemma}}
\begin{proof}
\ser{\fyr{Since the positive-part operator is monotone and subadditive, and since
$1-\theta_k\geq0$ and $\eta_k-\eta_{k+1}\geq0$, taking positive parts in
\eqref{def:d_rho} gives}
\begin{align}\label{eq:positive_d_recursion}
[d_{k+1}]_+
\leq
(1-\rho\sqrt{\eta_k})[d_k]_+
+
(\eta_k-\eta_{k+1})[-b_{k+1}]_+.
\end{align}
\noindent(i)
Using $[-b_{k+1}]_+\leq B$ in
\eqref{eq:positive_d_recursion} and dividing by $\eta_{k+1}>0$, \fyr{we obtain}
\begin{align}\label{eq:positive_d_recursion2}
\tfrac{[d_{k+1}]_+}{\eta_{k+1}}
\leq
(1-\rho\sqrt{\eta_k})
(
1+\tfrac{\eta_k-\eta_{k+1}}{\eta_{k+1}}
)
\tfrac{[d_k]_+}{\eta_k}
+
B(\tfrac{\eta_k-\eta_{k+1}}{\eta_{k+1}}).
\end{align}
\fyr{The assumption $\lim_{k\to \infty}\tfrac{\eta_k-\eta_{k+1}}
{\eta_{k+1}\sqrt{\eta_k}}
=0$ implies that for an arbitrary} $\varepsilon>0$, for all
sufficiently large $k$,
$(1-\rho\sqrt{\eta_k})
(1+\tfrac{\eta_k-\eta_{k+1}}{\eta_{k+1}})
\leq1-\tfrac{\rho\sqrt{\eta_k}}{2}$
and
$B(\tfrac{\eta_k-\eta_{k+1}}{\eta_{k+1}})
\leq\tfrac{\rho\sqrt{\eta_k}}{2}\varepsilon$.
Therefore, \fyr{subtracting $\varepsilon$ from both sides in \eqref{eq:positive_d_recursion2}, followed by  taking positive-parts on both sides, and utilizing the preceding relations, we obtain} $[
\tfrac{[d_{k+1}]_+}{\eta_{k+1}}-\varepsilon
]_+
\leq
(1-\tfrac{\rho\sqrt{\eta_k}}{2})
[
\tfrac{[d_k]_+}{\eta_k}-\varepsilon
]_+.$
Iterating from a sufficiently large $k_\varepsilon$ yields
\begin{align*}
[
\tfrac{[d_{k+1}]_+}{\eta_{k+1}}-\varepsilon
]_+
\leq
\textstyle\prod_{j=k_\varepsilon}^{k}
(1-\tfrac{\rho\sqrt{\eta_j}}{2})
[
\tfrac{[d_{k_\varepsilon}]_+}{\eta_{k_\varepsilon}}
-\varepsilon
]_+.
\end{align*}
Since $1-x\leq\exp(-x)$ and
$\textstyle\sum_{k=0}^{\infty}\sqrt{\eta_k}=\infty$, the product converges to
zero. Hence,
$\limsup_{k\to\infty}[d_k]_+/\eta_k\leq\varepsilon$.
Since $\varepsilon>0$ is arbitrary,
$\lim_{k\to\infty}[d_k]_+/\eta_k=0$.
Moreover,
$0\leq a_k\leq[d_k]_++\eta_k[-b_k]_+
\leq\eta_k([d_k]_+/\eta_k+B)$,
and hence $\lim_{k\to\infty}a_k=0$.
Also, $a_k+\eta_kb_k=d_k\leq[d_k]_+$ and $a_k\geq0$ imply
$b_k\leq[d_k]_+/\eta_k$, and thus
$\limsup_{k\to\infty}b_k\leq0$.
Since $\lim_{k\to\infty}a_k=0$, the assumption
$\lim_{k\to\infty}a_k=0
\Longrightarrow
\lim_{k\to\infty}[-b_k]_+=0$
gives $\lim_{k\to\infty}[-b_k]_+=0$, and hence
$\liminf_{k\to\infty}b_k\geq0$.
Therefore, $\lim_{k\to\infty}b_k=0$.

\noindent(ii)
Let 
$A_1:=\tfrac{2(m-1)}{m}
(\tfrac{2C^m}{m})^{{1}/{(m-1)}}$,
$A_2:=2C A_1^{1/m}$, and
$A_3:=\tfrac{2(m-1)}{m}
(\tfrac{4C^m}{m\rho})^{{1}/{(m-1)}}$.
The coupling condition and
$d_k=a_k+\eta_kb_k\leq[d_k]_+$
imply
$a_k-C\eta_k\sqrt[m]{a_k}\leq[d_k]_+$.
By Young's inequality, $C\eta_k\sqrt[m]{a_k}
\leq
\tfrac{a_k}{2}
+
\tfrac{m-1}{m}
(
\tfrac{2C^m\eta_k^m}{m}
)^{{1}/{(m-1)}}.$
Combining the preceding two inequalities and multiplying by $2$ yield $a_k
\leq
2[d_k]_+
+
A_1\eta_k^{{m}/{(m-1)}}.$}
\ser{Using the coupling condition, the preceding inequality, and
$\sqrt[m]{x+y}\leq\sqrt[m]{x}+\sqrt[m]{y}$ gives
\begin{align}\label{eq:negative_b_bound_general_m}
[-b_{k+1}]_+
\leq
\sqrt[m]{2C^m[d_{k+1}]_+}
+
C A_1^{1/m}\eta_{k+1}^{{1}/{(m-1)}}.
\end{align}
Furthermore, Young's inequality gives
\begin{align}\label{eq:young_d_bound_general_m}
&(\eta_k-\eta_{k+1})
\sqrt[m]{2C^m[d_{k+1}]_+}
\leq
\tfrac{\rho\sqrt{\eta_k}}{2}[d_{k+1}]_+
+
\tfrac{m-1}{m}
(
\tfrac{4C^m(\eta_k-\eta_{k+1})^m}
{m\rho\sqrt{\eta_k}}
)^{\tfrac{1}{m-1}}.
\end{align}
Substituting \eqref{eq:negative_b_bound_general_m} and
\eqref{eq:young_d_bound_general_m} into
\eqref{eq:positive_d_recursion}, and using
$\tfrac{1-\rho\sqrt{\eta_k}}
{1-\rho\sqrt{\eta_k}/2}
\leq1-\tfrac{\rho\sqrt{\eta_k}}{2}$
and
$\tfrac{1}{1-\rho\sqrt{\eta_k}/2}\leq2$, yield
\begin{align}\label{eq:d_recursion_reduced_general_m}
[d_{k+1}]_+
&\leq
(1-\tfrac{\rho\sqrt{\eta_k}}{2})[d_k]_+
+
A_2(\eta_k-\eta_{k+1})
\eta_{k+1}^{\tfrac{1}{m-1}}
+
A_3
(
\tfrac{(\eta_k-\eta_{k+1})^m}
{\sqrt{\eta_k}}
)^{\tfrac{1}{m-1}}.
\end{align}
}
\ser{By Lemma~\ref{lem:eta_difference_bound} and since
$\tfrac{m(a+1)-a/2}{m-1}>1+\tfrac{am}{m-1}$, the last two terms
in~\eqref{eq:d_recursion_reduced_general_m} are bounded by
$A_4(k+1)^{-1-\frac{am}{m-1}}$, where
$A_4:=aA_2\eta_0^{\frac{m}{m-1}}
+A_3a^{\frac{m}{m-1}}\eta_0^{\frac{2m-1}{2(m-1)}}$.
Thus, $[d_{k+1}]_+
\leq
(
1-\tfrac{\rho\sqrt{\eta_0}}
{2(k+1)^{a/2}}
)[d_k]_+
+
{A_4}/
{(k+1)^{1+\tfrac{am}{m-1}}}.$
Applying Lemma~\ref{lem:chung} with
$p=\tfrac{a}{2}$,
$q=1+\tfrac{a(m+1)}{2(m-1)}$, and contraction coefficient
$\tfrac{\rho\sqrt{\eta_0}}{2}$ yields
$[d_k]_+
=\mathcal{O}(k^{-1-\tfrac{a(m+1)}{2(m-1)}})$.
Combining this estimate with $a_k
\leq
2[d_k]_+
+
A_1\eta_k^{{m}/{(m-1)}}$ 
and
$a<2$ proves \fyr{the bounds for $a_k$}.
Finally,
$a_k+\eta_kb_k=d_k\leq[d_k]_+$ and $a_k\geq0$ imply
$b_k\leq[d_k]_+/\eta_k
=\mathcal{O}(k^{-1-\tfrac{a(3-m)}{2(m-1)}})$.
The coupling condition and \fyr{the earlier bounds on $a_k$} also give
$b_k\geq-\mathcal{O}(k^{-\tfrac{a}{m-1}})$. \fyr{This completes the proof.}}
\end{proof}
\ser{In the following theorem, we \fyr{establish} convergence \fyr{guarantees} for IR-VFISTA$_\text{s}$.}
\ser{\begin{theorem}\em\label{theorem:IRVFISTA_var_rates}
Consider \fyr{the problem in}~\eqref{prob:uni_centr} under
Assumption~\ref{assump:strongly_main}. Let $x^*$ \fyr{denote its} unique optimal solution. Consider Algorithm~\ref{alg:IR-R-VF_ne}. Suppose 
$\{\eta_k\}$ is positive and nonincreasing.  \fyr{Define $t_k := \sqrt{\kappa_k}$ for $k\geq0$}.

\noindent (i) {[asymptotic convergence]} \fyr{Assume that}
$\lim_{k\to \infty}\eta_k=0,\ 
\textstyle\sum_{k=0}^{\infty}\sqrt{\eta_k}=\infty,$
\allowbreak\ \ssg{and}\
$\lim_{k \to \infty}\tfrac{\eta_k-\eta_{k+1}}
{\eta_{k+1}\sqrt{\eta_k}}=0.$
Then, $\lim_{k\to\infty}x_k=x^*.$
\\
\noindent (ii) Suppose $X_{\bar h}^*$ is $\alpha$-weak sharp
of order $m\in(1,2]$. Let $\eta_k
=
\eta_0(k+1)^{-a},
\
a\in(0,2),
\
\eta_0>0.$ \fyr{Then,} for any $K\geq 1$, \fyr{the following holds.}
\\
\noindent (ii.1) [infeasibility bounds]   $0
\leq
\bar h(x_K)-\bar h(x^*)
\leq
\mathcal{O}(K^{-\tfrac{am}{m-1}}).$
\\
\noindent (ii.2) [suboptimality bounds]  $-\mathcal{O}(
K^{-\tfrac{a}{m-1}}
)
\leq
\bar f(x_K)-\bar f(x^*)
\leq
\mathcal{O}(K^{-1-\tfrac{a(3-m)}{2(m-1)}}).$
\\
\noindent (ii.3) [distance to the unique optimal solution] $\| x_K -x^*\|^2 \leq \mathcal{O}(K^{-\tfrac{a}{m-1}}) .$
\end{theorem}
\begin{proof}
Consider Lemma~\ref{lemma:general_rho}.  Let us define $a_k
:=\bar h(x_k)-\bar h(x^*)
+\tfrac{\eta_k\mu_f}{2}
\|t_{k-1}x_k
-(
x^*+(t_{k-1}-1)x_{k-1}
)
\|^2$, $b_k:=
\bar f(x_k)-\bar f(x^*),$ $\theta_k
:=
\tfrac{1}{t_k} = \tfrac{1}{\sqrt{\kappa_k}}
=
\sqrt{
\tfrac{\eta_k\mu_f}
{L_h+\eta_kL_f}
}.$ \fyr{Next, we verify the conditions} in Lemma~\ref{lemma:general_rho}. \fyr{Since} $x^*\in X_{\bar h}^*$, we have
$\bar h(x_k)-\bar h(x^*)\geq0$, \fyr{which verifies $ a_k \geq 0$. Let us set} $\rho :=\sqrt{
\tfrac{\mu_f}
{L_h+\eta_0L_f}
}.$ \fyr{Since} $\{\eta_k\}$ is nonincreasing, $\mu_f \leq L_f$,  and $L_h \geq 0$, we have
$\rho \sqrt{\eta_k} \leq \theta_k \leq 1$. Moreover, using \eqref{eq:IRVFISTA_var_dk} and 
$\bar g_{\eta_k}\fyr{(\cdot)}=\bar h\fyr{(\cdot)}+\eta_k\bar f\fyr{(\cdot)}$,  for all $k\geq1$,
we obtain $d_k(x^*)
=
a_k+\eta_k b_k.$
Furthermore, letting $x=x^*$ in
\eqref{eq:IRVFISTA_var_d_recursion_pro}, we obtain $d_{k+1}(x^*)
\leq
(1-\theta_k)d_k(x^*)
-
(\eta_k-\eta_{k+1})b_{k+1}.$
\fyr{Hence, the recursion and the basic conditions of} Lemma~\ref{lemma:general_rho} are satisfied. 

\noindent (i) 
As $\inf_{x\in\mathbb{R}^n}\bar f(x)>-\infty$, then $-b_k
=
\bar f(x^*)-\bar f(x_k)
\leq
\bar f(x^*)
-
\inf_{x\in\mathbb{R}^n}\bar f(x),
\ k\geq0.$
Therefore, $\{b_k\}$ is bounded below.
We next show that $\lim_{k\to\infty}a_k=0
\ \Longrightarrow\
\lim_{k\to\infty}[-b_k]_+=0.$
Suppose  $\lim_{k\to\infty}a_k=0.$
Then,  $0
\leq
\bar h(x_k)-\bar h(x^*)
\leq
a_k$,
yields $\lim_{k\to\infty}
(
\bar h(x_k)-\bar h(x^*)
)
=
0.$
For any
$\tilde{\nabla}\bar f(x^*)\in\partial\bar f(x^*)$,
by \eqref{eq:lower_f_strongly}, we have
\begin{align}\label{eq:IRVFISTA_var_error_bound_i}
-&
\|
\tilde{\nabla}\bar f(x^*)
\|
\operatorname{dist}
\left(
x_k,X_{\bar h}^*
\right)
+
\tfrac{\mu_f}{2}
\|x_k-x^*\|^2
\leq
b_k.
\end{align}
Suppose, by contradiction, that $\lim_{k\to\infty}[-b_k]_+\neq0.$
Then, there exist $\epsilon>0$ and a subsequence
$\{x_{k_i}\}$ such that $b_{k_i}
\leq
-\epsilon.$
Since $x^*\in X_{\bar h}^*$, we have
$\operatorname{dist}
(x_{k_i},X_{\bar h}^*)
\leq
\|x_{k_i}-x^*\|.$
Therefore, \eqref{eq:IRVFISTA_var_error_bound_i} and $b_{k_i}
\leq
-\epsilon$,
yield
\begin{align}\label{eq:IRVFISTA_var_error_bound_i_2}
\tfrac{\mu_f}{2}\|x_{k_i}-x^*\|^2
&\leq
b_{k_i}
+\|\tilde{\nabla}\bar f(x^*)\|\operatorname{dist}\left(x_{k_i},X_{\bar h}^*\right)
\leq \|\tilde{\nabla}\bar f(x^*)\|\|x_{k_i}-x^*\|.
\end{align}
Hence,
$\{x_{k_i}\}$ is bounded. Moreover, \fyr{the assumption}
$\lim_{k\to\infty}(\bar h(x_k)-\bar h(x^*))=0$ \fyr{implies} that
$\lim_{i\to\infty}\bar h(x_{k_i})=\bar h^*$. \fyr{Since $\{x_{k_i}\}$ is bounded, let $x_{k_{i_j}}$ denote an arbitrary convergent subsequence of $\{x_{k_i}\}$ with the limit point $\hat{x}$. From Assumption~\ref{assump:strongly_main}, $\bar{h}$ is lower semicontinuous. Thus,  $\bar h(\hat{x})\leq
\liminf_{j\to\infty} \bar h({x}_{k_{i_j}})=\bar h^*$. We also have $\bar h^*\leq \bar h(\hat{x})$, implying that $\hat{x} \in X^*_{\bar{h}}$. Since $\operatorname{dist}({x}_{k_i},X^*_{\bar{h}}) \leq \|{x}_{k_i}-\hat{x}\|$, \eqref{eq:IRVFISTA_var_error_bound_i_2} implies that $ \tfrac{\mu_f {\|x_{k_i} - x^*\|}^2}{2}   \leq    \|\tilde{\nabla}\bar f(x^*)\|\operatorname{dist}(x_{k_i},X_{\bar h}^*) \leq \|\tilde{\nabla}\bar f(x^*)\|\|x_{k_i}-\hat{x}\|.$ Taking the limit along $\{{x}_{k_{i_j}}\}$ yields $x^* = \hat{x}$. Since $\{{x}_{k_{i_j}}\}$ is an arbitrary convergent subsequence, we conclude that $\lim_{k\to\infty}{x}_{k_i}=x^*$, implying that} $\lim_{i\to\infty}\operatorname{dist}
(x_{k_i},X_{\bar h}^*)=0$.
Dropping the nonnegative quadratic term from
\eqref{eq:IRVFISTA_var_error_bound_i}, we obtain $b_{k_i}
\geq
-
\|
\tilde{\nabla}\bar f(x^*)
\|
\operatorname{dist}(
x_{k_i},X_{\bar h}^*).$
Consequently, $\liminf_{k_i\to\infty}b_{k_i}
\geq
0$
which contradicts $b_{k_i}
\leq
-\epsilon.$
Hence, \fyr{if $\lim_{k\to\infty}a_k=0$, then  $\lim_{k\to\infty}[-b_k]_+=0$.} 
Therefore, all the conditions of
Lemma~\ref{lemma:general_rho}~(i) are satisfied. \fyr{Applying that result, we} obtain $\lim_{k \to \infty} (\bar f(x_k)-\bar f(x^*))=
0$ and $\lim_{k \to \infty} (\bar h(x_k)-\bar h(x^*))=
0.$\\
We next show that $\{x_k\}$ is bounded. Since $\lim_{k\to\infty}b_k
=
0$, there
exists $k_0\geq0$ such that $b_k\leq1,
\ k\geq k_0.$
Using \eqref{eq:IRVFISTA_var_error_bound_i} and
$\operatorname{dist}
(
x_k,X_{\bar h}^*)
\leq
\|x_k-x^*\|,$
we obtain, for all $k\geq k_0$,
$\tfrac{\mu_f}{2}
\|x_k-x^*\|^2
\leq
1
+
\|
\tilde{\nabla}\bar f(x^*)
\|
\|x_k-x^*\|.$
Thus,  $\{x_k\}$ is bounded.
Since $\{x_k\}$ is bounded, it has at least one convergent
subsequence. Let $\{x_{k_t}\}$ be an arbitrary convergent
subsequence and let $\lim_{k_t \to \infty}x_{k_t}
=
\hat x.$
By the lower semicontinuity of $\bar h$, we obtain $\bar h(\hat x)
\leq
\liminf_{k_t\to\infty}\bar h(x_{k_t})
=
\bar h(x^*).$
Since $\bar h(x^*)$ is the minimum value of $\bar h$, it follows
that $\hat x\in X_{\bar h}^*$. Therefore,
$\lim_{k_t \to \infty}\operatorname{dist}
(x_{k_t},X_{\bar h}^*)
\leq
\lim_{k_t \to \infty}\|x_{k_t}-\hat x\|
=
0.$
Taking the limit along $\{x_{k_t}\}$ in
\eqref{eq:IRVFISTA_var_error_bound_i}, and using
$\lim_{k_t \to \infty}b_{k_t}=0$, we obtain
$\tfrac{\mu_f}{2}
\|\hat x-x^*\|^2
\leq
0.$
Thus, $\hat x=x^*$. Hence, every limit point of
$\{x_k\}$ is equal to $x^*$. Since $\{x_k\}$ is bounded, it
follows that $\lim_{k\to\infty}x_k
=
x^*.$

\noindent (ii) For any
$\tilde{\nabla}\bar f(x^*)
\in\partial\bar f(x^*)$, let us define $C
:=
\|
\tilde{\nabla}\bar f(x^*)
\|
\alpha^{-\tfrac{1}{m}}.$
Since $X_{\bar h}^*$ is $\alpha$-weak sharp of order $m$,
\eqref{prob:distance-weak-strongly} gives, for all $k\geq1$,
\begin{align}\label{eq:IRVFISTA_var_previous_lemma_ii}
\|x_k-x^*\|^2
\leq
\tfrac{2}{\mu_f}
\left(
b_k
+
C
\ser{\sqrt[m]{
\bar h(x_k)-\bar h(x^*)
}}
\right).
\end{align}
Since the left-hand side of
\eqref{eq:IRVFISTA_var_previous_lemma_ii} is nonnegative, 
\fyr{using} $0
\leq
\bar h(x_k)-\bar h(x^*)
\leq
a_k$, 
we  obtain
$b_k\geq-C\sqrt[m]{a_k}.$
Thus, the coupling condition in
Lemma~\ref{lemma:general_rho}~(ii) is satisfied. Therefore, all the conditions of
Lemma~\ref{lemma:general_rho}~(ii) are \fyr{verified}. Consequently,  we obtain the desired results in (ii.1) and (ii.2). Moreover, regarding \eqref{eq:IRVFISTA_var_previous_lemma_ii}, invoking  (ii.1) and (ii.2), and considering   $a \in (0,2)$ and $m \in (1,2]$, we obtain (ii.3).
\end{proof}}
\ser{\begin{remark}\label{remark:3.17}
As a consequence of Theorem~\ref{theorem:IRVFISTA_var_rates}, in the special case $m=2$, for any $\epsilon>0$, by choosing $a<2$ sufficiently close to $2$, we obtain $-\mathcal{O}(K^{-2+\epsilon}) \leq \bar f(x_K)-\bar f(x^*) \leq \mathcal{O}(K^{-2+\epsilon})$ and $0 \leq \bar h(x_K)-\bar h(x^*) \leq \mathcal{O}(K^{-4+\epsilon})$. Therefore, compared with the convergence rates of \cite{merchav2024dynamic} in the case $m=2$, as summarized in Table~\ref{table:lit}, IR-VFISTA \fyr{improves both the suboptimality and infeasibility convergence rates}.
In addition, Theorem~\ref{theorem:IRVFISTA_var_rates} provides convergence rates for the broader range $m\in(1,2]$, whereas
\cite{merchav2024dynamic} considers only the case $m=2$. For $m\in(1,2)$, by choosing $a<2$ sufficiently close to $2$, we obtain rates faster than those obtained in the case $m=2$. For example, when $m=1.5$, for any $\epsilon>0$, by choosing $a<2$ sufficiently close to $2$, we obtain $-\mathcal{O}(K^{-4+\epsilon}) \leq \bar f(x_K)-\bar f(x^*) \leq \mathcal{O}(K^{-4+\epsilon})$ and $0 \leq \bar h(x_K)-\bar h(x^*) \leq \mathcal{O}(K^{-6+\epsilon})$. \fyr{Compared with \cite{zhang2024functionally}, which considers SBO problems with convex upper- and lower-level objectives, IR-VFISTA${\mathrm{s}}$ exploits strong convexity of the upper-level objective and provides substantially sharper rates under weak sharp minimality of order 2. Moreover, IR-VFISTA${\mathrm{s}}$ is a single-loop anytime method that guarantees asymptotic convergence without requiring weak sharp minimality.}
\end{remark}}
In the \fyr{following}, we establish convergence rate statements for \fyr{a variant of $\text{IR-VFISTA}_\text{s}$ with a constant regularization parameter, which we refer to as $\text{R-VFISTA}_\text{s}$.} \ssg{We provide the proof for the next corollary in  Appendix~\ref{appendix:proof_cor18}.}
\begin{corollary}[error bounds for $\text{R-VFISTA}_\text{s}$]\label{Cor:R-VFISTA}\em
 Suppose Assumption~\ref{assump:strongly_main} holds. Let $x_K$ be  generated by \ser{Algorithm~\ref{alg:IR-R-VF_ne} with a constant regularized parameter $\eta$.} The results in the following two settings hold.
 
\noindent (1) Let \ch{$\eta = \tfrac{(L_h + \bar{\eta}L_f)}{\mu_f}(\tfrac{ (p+1) \newchr{\ln({K})}}{ \newchr{K}})^2$},  where $\bar{\eta} > 0$ and   $p > 2$   are arbitrary.  Then, for $K$  satisfying  \ch{$ \tfrac{ (L_h + \bar{\eta}L_f) (p+1)^2}{\mu_f \bar{\eta}} \leq \newch{(\tfrac{K}{\ln(K)})^2} $}, \ch{the following statements hold.}

\noindent {{(1.i)}}  [suboptimality  bounds]
\begin{align*}
 -\|\tilde{\nabla} \bar{f}(x^*)\|  \operatorname{dist}(x_K,X^*_{\bar{h}})+\tfrac{\mu_f}{2} {\|x_K -x^* \|}^2 &\leq \bar{f}(x_K) - \bar{f}^* 
\leq\tfrac{u_6}{{K}^{(p+1)}} 
 +\tfrac{\newchr{u_7}}{\newchr{{K}^{(p-1)}}\ln({K})},
 \end{align*}
 where $u_6 :=  \bar{f}(x_0) - \bar{f}^*   + \tfrac{{\mu_f}}{2} {\|x_0 - x^*\|}^2 $ and \ch{$u_7 := \tfrac{\mu_f (\bar{h}(x_0) -\bar{h}^*)}{{(L_h + \bar{\eta}L_f) (p+1) } }$}.
 
\noindent {{(1.ii)}}  [infeasibility  bounds] \ser{$0 \leq \bar{h}(x_K) -{\bar{h}^*}  
  \leq u_8 (\tfrac{\newch{\ln({K})}}{\newch{K}})^2+ \tfrac{\newch{u_9(\ln({K}))^2}}{\newchr{K}^{(p+3)}}    
  +\tfrac{u_{10}}{\newchr{{K}^{(p+1)}}},$ }
 where \ch{$u_8 := \tfrac{{( \bar{f}(\Pi_{X^*_{\bar{h}}}[x_0])-\hat{C}_{\bar{f}})}(L_h + \bar{\eta}L_f) (p+1)^2 }{\mu_f }$},  \ch{$u_9 := \tfrac{(L_h + \bar{\eta}L_f) (p+1)^2 (  \bar{f}(x_0) -\hat{C}_{\bar{f}}  +0.5{{\mu_f} \operatorname{dist}^2(x_0, X^*_{\bar{h}})} )  }{\mu_f }$},\\ and $u_{10} := \bar{h}(x_0) - {\bar{h}^*} $. 
 
\ser{\noindent {(1.iii)}} [distance to optimal solution]  If $X^*_{\bar{h}}$ is $\alpha$-weak sharp minima of order $m \geq 1$, 
\begin{align*}
\tfrac{\mu_f}{2} {\|x_{\ssg{K}} - x^*\|}^2  
&\leq\tfrac{u_6}{{K}^{(p+1)}} 
+\tfrac{\newchr{u_7}}{\newchr{{K}^{(p-1)}}\ln({K})} 
 +   \tfrac{\|\tilde{\nabla} \bar{f}(x^*)\|}{ \sqrt[m]{\alpha}} \sqrt[m]{u_8 (\tfrac{\newch{\ln({K})}}{\newch{K}})^2+ \tfrac{\newch{u_9(\ln({K}))^2}}{\newchr{K}^{(p+3)}}    
  +\tfrac{u_{10}}{\newchr{{K}^{(p+1)}}}}.
\end{align*}
\ser{\noindent {(1.iv)}} [suboptimality lower bound] Furthermore, if $X^*_{\bar{h}}$ is $\alpha$-weak sharp minima of order $m \geq 1$, then for any $K\geq 1$, we obtain
  \begin{align*}
  -\tfrac{\|\tilde{\nabla} \bar{f}(x^*)\|}{ \sqrt[m]{\alpha}} \sqrt[m]{u_8 (\tfrac{\newch{\ln({K})}}{\newch{K}})^2+ \tfrac{\newch{u_9(\ln({K}))^2}}{\newchr{K}^{(p+3)}}    
  +\tfrac{u_{10}}{\newchr{{K}^{(p+1)}}}}\leq\bar{f}(x_K) - \bar{f}^*.
  \end{align*} 
\noindent (2) Assume that $X^*_{\bar{h}}$ is $\alpha$-weak sharp minima of the order $m = 1$. Let  $\eta \leq \tfrac{\alpha}{2 \|\tilde{\nabla} \bar{f}(x^*)\|} $, for some $  \tilde{\nabla} \bar{f}(x^*) \in \partial \bar{ f}(x^*)$.  Let \ch{us define} $u_{11} :=  \bar{f}(x_0)-\bar{f}^* + \tfrac{\bar{h}(x_0)- \bar{h}^*}{\eta} + \tfrac{ \mu_f {\|x_0-x^*\|}^2}{2} $. Then, for any $ K \geq 1$, we have the following statements.\\
\noindent {{(2.i)}} [infeasibility  bounds] 
$
   0 \leq \bar{h}(x_K) -{\bar{h}^*}  
  \leq   2\eta u_{11} (1 - \tfrac{1}{\sqrt{\kappa}} )^K.$
  
\noindent {{(2.ii)}} [distance to the lower-level solution set] $
     \operatorname{dist}(x_K,X^*_{\bar{h}}) \leq  \tfrac{2\eta u_{11}}{\alpha}  (1 - \tfrac{1}{\sqrt{\kappa}} )^K.$
     
\noindent {{(2.iii)}} [suboptimality  bounds]
$
  -u_{11} (1 - \tfrac{1}{\sqrt{\kappa}} )^K \leq \bar{f}(x_K)-{\bar{f}^*} \leq  u_{11} (1 - \tfrac{1}{\sqrt{\kappa}} )^K.$
  
 \noindent {{(2.iv)}} [distance to optimal solution] 
$ {\|x_{\ssg{K}} - x^*\|}^2  \leq \tfrac{4u_{11}}{\mu_f}  (1 - \tfrac{1}{\sqrt{\kappa}} )^K.$
\end{corollary}
\begin{remark} In \fyr{Corollary}~\ref{Cor:R-VFISTA},  \fyr{choosing} $p=3$, we derive quadratically decaying sublinear convergence rates of the order $1/K^2$ for both infeasibility and suboptimality error metrics. These appear to be the fastest known rates in addressing the SBO problem in \eqref{prob:uni_centr} \ser{under Assumption~\ref{assump:strongly_main}}.  When weak sharp minimality holds,  a linear convergence rate is obtained. Notably, when compared to R-ISTA$_\text{s}$, the linear rate has an improved dependence on the condition number. 
\end{remark}

\section{Composite SBO  with a smooth nonconvex upper-level objective} \label{sec:5}
In this section, we consider the following problem.
\begin{align}
\min \ & f(x), \quad \text{s.t.} \quad \quad x \in \arg\min_{x \in \mathbb{R}^n} \bar{h}(x) := h(x) + \omega_h (x),
\label{prob:centr-nonconvex}
\end{align}
where $f(x)$ is not necessarily convex. The formal assumptions are stated below.
 \begin{assumption}\label{assump:nonconvex_main}\em Consider the problem \ch{in}~\eqref{prob:centr-nonconvex}. Let the following hold.\\
\noindent\ser{(i)} $ f: \mathbb{R}^n \rightarrow \mathbb{R}$ is   $L_f$-smooth and possibly nonconvex.\\
\noindent\ser{(ii)}  $h: \mathbb{R}^n \rightarrow \mathbb{R}$ is  $L_h$-smooth and convex.\\
\noindent\ser{(iii)}   $\omega_h: \mathbb{R}^n \rightarrow (-\infty, \infty]$ is proper, closed, and convex.\\
\fyr{\noindent \ser{(iv)}  The set {$X^*_{\bar{h}}$} is nonempty.}\\
\noindent\fyr{(v)} $\inf_{x \in \mathbb{R}^n} f(x) > -\infty$.
\end{assumption}
\begin{definition}\label{def:residual}\em
Let Assumption~\ref{assump:nonconvex_main} hold and $0<\hat{\gamma}<\tfrac{1}{L_f}$. \ser{For any $x\in\mathbb{R}^n$, define the residual mapping}
$G_{1/\hat{\gamma}}(x):=\tfrac{1}{\hat{\gamma}}(x-\Pi_{X^*_{\bar h}}[x-\hat{\gamma}\nabla f(x)])$ \ser{\cite[Definition 10.5]{beck2017first}. In view of \cite[Theorem 10.7]{beck2017first}, $x^*$ is a stationary point of the problem in~\eqref{prob:centr-nonconvex} if and only if $G_{{1}/{\hat{\gamma}}}(x^*)=0$.}
\end{definition}

 \subsection{The  \texorpdfstring{\ssg{$\text{IPR-VFISTA}_{\text{nc}}$}}{IPR-VFISTA-nc} method}
We propose \fyr{an} {Inexactly Projected Regularized Variant of FISTA} $(\text{IPR-VFISTA}_{\text{nc}})$, presented in Algorithm~\ref{alg:IPR-VFISTA}. \fyr{At each outer-loop iteration $k$, it employs
$\text{R-VFISTA}_\text{s}$, using the modified constant
regularization parameter $\eta_k$ for an adaptive number of inner-loop
iterations $J_k$. The parameters
$J_k$, $\eta_k$, and $\kappa_k$ are \fyr{iteratively} updated at each outer-loop iteration. To this end}, we reformulate the problem \ch{in}~\eqref{prob:centr-nonconvex} as
$
\min_{x \in X^*_{\bar{h}}} f(x).
$
One may consider the standard projected gradient method
$\hat{x}_{k+1} = \Pi_{X^*_{\bar{h}}}[\hat{x}_k - \hat{\gamma} \nabla f(\hat{x}_k)].$
However, since $X^*_{\bar{h}}$ is not explicitly known, \fyr{an exact projection cannot
be computed}. Motivated by the work in~\cite{samadi2025improved} for addressing nonconvex optimization with variational \ser{inequality} constraints, \fyr{we instead employ an
inexact projection.} Given 
$z_k := \hat{x}_k - \hat{\gamma} \nabla f(\hat{x}_k),$
consider the projection problem
\begin{align} \label{prob:proj}
\min_{x \in X^*_{\bar{h}}} f_{s,k}(x) := \tfrac{1}{2}{\|x - z_k\|}^2,
\end{align}
where we use the label $s$ to \ser{emphasize} the strong convexity of the objective. Let $x^*_{f_{\fyr{s,k}}}$ denote the unique optimal solution to the problem in \eqref{prob:proj}, i.e.,
$
x^*_{f_{\fyr{s,k}}} = \Pi_{X^*_{\bar{h}}}[z_k].
$
Since the problem \ch{in}~\eqref{prob:proj} satisfies Assumption~\ref{assump:strongly_main}, \ser{Algorithm~\ref{alg:IR-R-VF_ne}} can be employed to obtain this inexact projection.
However, inexact projections may lead to infeasibility of $\hat{x}_k$ for the problem \ch{in}~\eqref{prob:centr-nonconvex}. \fyr{To establish convergence
guarantees, we design the sequences $J_k$, $\eta_k$, and $\kappa_k$. In particular, the number of inner-loop iterations
$J_k$ is chosen so that $\lim_{k\to\infty}J_k=\infty$, with its growth determined
by properties of the lower-level problem, while $\eta_k$ and $\kappa_k$ are
defined accordingly. Finally,} we initialize $\text{R-VFISTA}_\text{s}$ at each outer-loop iteration $k$ with
$y_{k+1,0} = x_{k+1,0} = \Pi_{\mathcal{B}} \fyr{\left[x_{k,J_k}\right]},$
where \ser{$\mathcal{B} \subseteq \mathbb{R}^n$} is an arbitrary \fyr{compact set} and $x_{k,J_k}$ is the output of $\text{R-VFISTA}_\text{s}$ after $J_k$ inner iterations at iteration $k$. \fyr{This initialization ensures that the sequence of inner-loop starting points remains bounded.}
 \begin{algorithm}   
 	\caption{ Inexactly Projected Regularized  VFISTA ($\text{IPR-VFISTA}_{\text{nc}}$)}
 	\begin{algorithmic}[1]\label{alg:IPR-VFISTA}
 		\STATE \textbf{input:} Initial vectors $\hat{x}_0$, \fyr{$y_{0,0} = x_{0,0}$}, outer-loop stepsize \ser{$\hat{\gamma} \leq \tfrac{1}{2L_f}$}, \fyr{an increasing sequence \ssg{of integers $\{J_k\}$ with $J_k\geq 2$},} \fyr{$\bar{\eta}>0$,}  
a \ssg{nonempty} \fyr{compact set} \ssg{$\mathcal{B} \subseteq \mbox{dom}(\omega_h)$}, and total number of outer-loop iterations \ser{$K \geq 1$}.

 		\FOR {$	k = 0, 1, 2, \dots, {K-1}$}
 		\STATE {$z_k = \hat{x}_k - \hat{\gamma}\nabla f(\hat{x}_k)$}
 		\STATE{\fyr{$\eta_k :=16(L_h + \bar{\eta}) ( \tfrac{ \newch{\ln({{J_k}})}}{ \newch{J_k}})^2 $,  $\gamma_k := \tfrac{1}{L_h + \eta_k}$}}, $\kappa_k := \tfrac{L_h + \eta_k}{\eta_k}$
 		\FOR {$j = 0, 1, 2, \dots, {J_{k}-1}$}
 		\STATE   {$x_{k,j+1} := \mbox{prox}_{\ch{\gamma_k} \omega_h}\left( y_{k,j} -\ch{\gamma_k} \left(\nabla h( y_{k,j}) + \eta_k  (y_{k,j} - z_k) \right)\right)$}	
        \STATE {$y_{k,j+1} := x_{k,j+1} + (\tfrac{\sqrt{\kappa_k} -1}{\sqrt{\kappa_k}+1})(x_{k,j+1}-x_{k,j})$}
        
 		\ENDFOR
 		\STATE {$\hat{x}_{k+1} := x_{k,J_k}, y_{{k+1},0}:=x_{{k+1},0}:=\Pi_{\mathcal{B}}[x_{k,J_k}] $}
 		\ENDFOR
	\STATE {{\bf return:} {${\hat{x}}_K$}}
	
 	\end{algorithmic}
 \end{algorithm}
\fyr{Our goal in this section is to establish  asymptotic convergence and nonasymptotic rate guarantees. A key challenge is that the objective function in~\eqref{prob:proj} varies with the outer-loop iteration $k$, since $z_k$ is updated recursively. We first introduce several auxiliary quantities and establish their boundedness properties, which will play a central role in the subsequent analysis.}
 \begin{definition}\label{def:eanddelta}\em
 Let us define  $z_k := \hat{x}_k - \hat{\gamma} \nabla f(\hat{x}_k)$, 
 $\delta_k := \hat{x}_{k+1} - \Pi_{{X}^*_{\bar{h}}}[z_k]$, and $e_k := \hat{x}_k - \Pi_{{X}^*_{\bar{h}}}[\hat{x}_k]$, i.e., $\|e_k\| := \operatorname{dist}(\hat{x}_k , X^*_{\bar{h}})$, for $k \geq 0$.  
\end{definition}
\begin{definition}\label{def:T}\em

\ssg{For $k\geq 0$ and $i\in\{1,\cdots, 5\}$, let us define $T_{i,k}$ as follows.} \ch{$T_{1,k} := (L_h + \bar{\eta})\big({f_{s,k}}(x_{k,0})+\fyr{\tfrac{1}{2}}\operatorname{dist}^2(x_{k,0},X^*_{\bar{h}})\big)$}, 
$T_{2,k} := \bar{h}(x_{k,0})-\bar{h}^*$,  
\ch{$T_{3,k} := (L_h+\bar{\eta}){f_{s,k}}(\Pi_{X^*_{\bar{h}}}[x_{k,0}])$}, 
$T_{4,k} := \ser{2{f_{s,k}}(x_{k,0})+\|x_{k,0}-x^*_{f_{\fyr{s,k}}}\|^2}$ 
{, and} 
\ch{$T_{5,k} := \tfrac{T_{2,k}}{2(L_h+\bar{\eta})}$}.

\end{definition}

 \begin{definition}\label{def:basic}\em

Let $B_{X^*_{\bar{h}}}:=\sup_{x\in\mathcal{B}}\operatorname{dist}(x,X^*_{\bar{h}})$, \ssg{\allowbreak}
$B_{\mathcal{B}}:=\sup_{x\in\mathcal{B}}\|x\|$, \ssg{\allowbreak}
$B_{\bar{h}}:=\sup_{x\in\mathcal{B}}\bar{h}(x)-\bar{h}^*$, \ssg{\allowbreak}
\fyr{and $\hat{C}_f:=\inf_{x\in\mathbb{R}^n}f(x)$. When $X^*_{\bar{h}}$ is bounded, we let
$D_{X^*_{\bar{h}}}:=\sup_{x\in X^*_{\bar{h}}}\|x\|$ \ssg{\allowbreak}
and $B_f:=\sup_{x\in\mathcal{B}\cup X^*_{\bar{h}}}\|\nabla f(x)\|$. In settings where we show that $\{\hat{x}_k\}$ is bounded, we let
$\hat{B}_f:=\sup_{k\geq1}\|\nabla f(\hat{x}_k)\|$ \ssg{\allowbreak}
and $\hat{D}_f:=\sup_{k\geq1}f(\hat{x}_k)$.}

\end{definition} 
\begin{lemma}\label{lemma:bounded_T}\em
       \fyr{The following properties hold. (i) We have} $\ T_{1,k-1}  \leq
c_1 +c_2  {\| \hat{x}_{k-1}\|}^2,$ \fyr{$T_{2,k-1} \leq B_{\bar{h}},$} and $T_{3,k-1}  \leq
c_3 + c_2   {\| \hat{x}_{k-1}\|}^2,$
\ser{where {$c_1 := 2(L_h + \bar{\eta}) (1.5 B_{\mathcal{B}}^2 + \hat{\gamma}^2  B_f^2   +0.25 B_{{X}_{\bar{h}}^*}^2 )$, $c_2 :=3(L_h + \bar{\eta}) $}, and $c_3:= 2 (L_h + \bar{\eta})(1.5  D_{X_{\bar{h}}^*}^2  + \hat{\gamma}^2  B_f^2)$. {(ii)} \fyr{If} $\{ \hat{x}_k\}$ \fyr{is} a bounded \fyr{sequence, then,} for all $i \in \{1,2,3,4,5\}$, there exists $M_i > 0$ such that $T_{i,k} \leq M_i$, \fyr{for all $k\geq 0$}.}
 \end{lemma}
 \begin{proof}
     \ser{\noindent {(i)} First, we derive an  upper bound for $T_{1,k-1}$. 
We have 
\begin{align}\label{eq:bound_fs}
{f_{s,k-1}}(x_{k-1,0})&=  0.5 {\| x_{k-1,0}- \hat{x}_{k-1} + \hat{\gamma} \nabla f(\hat{x}_{k-1})\|}^2\\\notag 
& \leq
{\| x_{k-1,0}- \hat{x}_{k-1}\|}^2 + {\| \hat{\gamma} \nabla f(\hat{x}_{k-1})-\hat{\gamma} \nabla f(x_{k-1,0})+\hat{\gamma} \nabla f({x}_{k-1,0})\|}^2 \\\notag &\leq 
{\| x_{k-1,0}- \hat{x}_{k-1}\|}^2 + 2 \hat{\gamma}^2{ L_f^2 {\|  {x}_{k-1,0} - \hat{x}_{k-1}\|}^2} + 2 \hat{\gamma}^2 B_f^2   \\\notag & \leq
2(1+2\hat{\gamma}^2 L_f^2) {{\|x_{k-1,0} \|}^2} + 2 (1+2\hat{\gamma}^2 L_f^2) {\| \hat{x}_{k-1}\|}^2 + 2\hat{\gamma}^2  B_f^2 \\\notag  & \leq 3 B_{\mathcal{B}}^2 +3  {\| \hat{x}_{k-1}\|}^2 + 2\hat{\gamma}^2  B_f^2,
\end{align} 
\noindent where we used $2\hat{\gamma}^2 L_f^2\leq 0.5$ in the last inequality. Hence, $T_{1,k-1}  \leq
c_1 +c_2  {\| \hat{x}_{k-1}\|}^2.$
Moreover, $T_{2,k-1}=  \bar{h}(x_{k-1,0}) - \bar{h}^* \leq B_{\bar{h}}.$
Following the approach used in \ser{\eqref{eq:bound_fs}}, 
$
{f_{s,k-1}}(\Pi_{X^*_{\bar{h}}}[x_{{k-1},0}])\leq 3  D_{X_{\bar{h}}^*}^2 + 3 {\| \hat{x}_{k-1}\|}^2 + 2\hat{\gamma}^2  B_f^2
$ yields $T_{3,k-1}  \leq
c_3 + c_2   {\| \hat{x}_{k-1}\|}^2.$}\\
 \ser{ \noindent {(ii)} \fyr{Let ${\|\hat{x}_k\|}^2 \leq M$, for all $k\geq 0$ and for some $M>0$}. \fyr{Substituting} ${\|\hat{x}_k\|}^2 \leq M$ into the bounds in \fyr{part (i)},    
 we conclude that there exist constants $M_1 =c_1+c_2 M $ and $M_3 = c_3+c_2 M$, such that $T_{1,k} \leq M_1$ and $T_{3,k} \leq M_3$.
 Let $M_2 = B_{\bar{h}}$. From part~(i), we \fyr{obtain} 
 $T_{2,k} \leq  M_2$.  The boundedness of $\{T_{4,k}\}$ follows by invoking $x^*_{f_{s,k}} \in X^*_{\bar{h}}$, $x_{k,0} \in \mathcal{B}$ \fyr{and \eqref{eq:bound_fs}. Lastly,} $T_{5,k} \leq M_5 $ for $M_5=\tfrac{M_2}{2(L_h + \fyr{\bar{\eta}})}$.}
 \end{proof}

\subsection{\texorpdfstring{Asymptotic convergence \fyr{for \ssg{$\text{IPR-VFISTA}_{\text{nc}}$}}}{Asymptotic convergence for IPR-VFISTA-nc}}
\label{subsection:4.2}
\fyr{We first establish that $\|e_k\|$, that is the distance from $\hat{x}_k$ to the lower-level optimal solution set $X^*_{\bar{h}}$, converges to zero as the outer-level iteration index $k$ goes to infinity. We also show that the norm of the projection error, i.e., $\|\delta_k\|$, converges to zero.}
\begin{lemma}\label{lemma:limit_dist_nonconvex_cons}\em
 \ser{Let $\{\hat{x}_k\}$ be generated by Algorithm~\ref{alg:IPR-VFISTA}, \fyr{and let} Assumption~\ref{assump:nonconvex_main} hold. Let $\{J_k\}$ be an increasing sequence of positive integers such that $\lim_{k \to \infty}J_k = \infty$. Suppose $\omega_h = \mathbb{I}_X$, where $X\subset\mathbb{R}^n$ is nonempty, compact, and convex. \fyr{Then, $\lim_{k \to \infty} \|e_k\|=\lim_{k \to \infty} \|\delta_k\| = 0.$}}
\end{lemma}
\begin{proof}
\ser{\fyr{From Algorithm~\ref{alg:IPR-VFISTA}, $x_{k-1,J_{k-1}}$} is the output of \ser{Algorithm~\ref{alg:IR-R-VF_ne} \fyr{with regularization parameter} $\eta_k$ \fyr{at} outer iteration $k$} after $J_{k-1}$ iterations for solving the problem in \eqref{prob:proj}. \fyr{From} Algorithm~\ref{alg:IPR-VFISTA}, we have $\hat{x}_{k} = x_{k-1,J_{k-1}}$ for $k \geq 1$.  \fyr{Invoking Corollary}~\ref{Cor:R-VFISTA} (1.ii) with $p=3$, noting the nonnegativity of $f_{s,k}(x)$, and invoking Lemma~\ref{lemma:bounded_T}, we obtain
\begin{align}\label{eq:bound_lim_nonconvex_cons}
0 \le \bar{h}(\hat{x}_k) - \bar{h}^*
\le \tfrac{16M_1(\ln(J_{k-1}))^2}{J_{k-1}^6}
+ \tfrac{M_2}{J_{k-1}^4}
+ \tfrac{16M_3(\ln(J_{k-1}))^2}{J_{k-1}^2}.
\end{align} 
Since $X$ is compact, the sequence $\{\hat{x}_k\}$ generated by Algorithm~\ref{alg:IPR-VFISTA} is bounded. Moreover, since $\lim_{k\to\infty}J_k=\infty$, the bound in
\eqref{eq:bound_lim_nonconvex_cons} yields
$\lim_{k\to\infty}\bar h(\hat{x}_k)=\bar h^*$.
\fyr{Following the same argument as in the proof of Theorem~\ref{proposition:F_eta_IR-ISTA-S} (iii), and using the lower semicontinuity of $\bar h$, we obtain
$\lim_{k \to \infty} \operatorname{dist}(\hat{x}_k,X^*_{\bar h})=0$, implying that $\lim_{k \to \infty}\|e_k\|=0$. Next, we show that $\lim_{k \to \infty}\|\delta_k\|=0$. From Definition~\ref{def:eanddelta},}  $\delta_k = \hat{x}_{k+1} - \Pi_{X^*_{\bar{h}}}[z_k]$. On the other hand, $\hat{x}_{k+1}$ is the output of Algorithm~\ref{alg:IR-R-VF_ne} with \fyr{the regularization parameter $\eta_k:= 16(L_h + \bar{\eta}) ( \tfrac{ {\ln({{J_k}})}}{ {J_k}})^2$ and $p:=3$} for solving the projection problem in~\eqref{prob:proj}. Therefore, by \fyr{Corollary}~\ref{Cor:R-VFISTA}~(i), we obtain
\begin{align}\label{eqn:ekdeltakasymp}
\ssg{0 \leq} \fyr{\ssg{\tfrac{1}{2}}}\ssg{\|\delta_k\|}^{\fyr{\ssg{2}}}
&\ssg{=} \fyr{\ssg{\tfrac{1}{2}}}
\ssg{\|\hat{x}_{k+1} - \Pi_{X^*_{\bar{h}}}[z_k]\|}^{\fyr{\ssg{2}}}
\\ \notag
&\ssg{\leq} 
\fyr{\ssg{
\tfrac{u_{6,k}}{J_k^{4}}
+\tfrac{u_{7,k}}{J_k^{2}\ln(J_k)}
}}
\ssg{
+ \|\nabla f_{s,k}(\hat{x}_{k+1})\|
\,\mathrm{dist}(\hat{x}_{k+1}, X^*_{\bar{h}})
}\ssg{,}
\end{align}
\fyr{where $u_{6,k}:= f_{s,k}(x_{k,0}) - f^*_{s,k}  + \tfrac{{1}}{2} {\|x_{k,0} - \Pi_{X^*_{\bar{h}}}[z_k]\|}^2$ and $u_{7,k}:= \tfrac{{h}(x_{k,0}) -{h}^*}{{4(L_h + \bar{\eta})   } }$. \ser{Note that for the constrained lower-level problem, since $\omega_h=\mathbb{I}_X$, each proximal update reduces to a projection
onto the feasible set $X$. Therefore, the rate results in
Corollary~\ref{Cor:R-VFISTA} \ssg{apply} to the  function $h$.} We have $\mathrm{dist}(\hat{x}_{k+1}, X^*_{\bar{h}}) = \|e_{k+1}\|$ and thus, it converges to zero as $k\to \infty$. To show $\lim_{k \to \infty}\|\delta_k\| = 0$, it suffices to prove that $\{u_{6,k}\}$, $\{u_{7,k}\}$ and $\{\|\nabla f_{s,k}(\hat{x}_{k+1})\|\}$ are bounded. From Assumption~\ref{assump:nonconvex_main}, we have that $\bar{h}^*>-\infty$, which together with $x_{k,0} \in \mathcal{B}$,  imply that $\{u_{7,k}\}$ is bounded. From \eqref{prob:proj} and $z_k = \hat{x}_k - \hat{\gamma}\nabla f(\hat{x}_k)$, we have $\|\nabla f_{s,k}(\hat{x}_{k+1})\| =\|\hat{x}_{k+1} - \hat{x}_k + \hat{\gamma}\nabla f(\hat{x}_k)\| $. Using $\hat{x}_{k+1} = x_{k,J_k}\in X$ for all $k$, together with the boundedness of $X$, imply that $\|\nabla f_{s,k}(\hat{x}_{k+1})\|$ is bounded. Next, from \eqref{prob:proj} and that $f^*_{s,k}=0$, we have $
    u_{6,k} = \tfrac{1}{2}\|x_{k,0} - \hat{x}_k + \hat{\gamma}\nabla f(\hat{x}_k)\|^2   + \tfrac{{1}}{2} {\|x_{k,0} - \Pi_{X^*_{\bar{h}}}[z_k]\|}^2.$ Noting that $X^*_{\bar{h}} (\subseteq X)$ is a bounded set, $x_{k,0}$ and $\hat{x}_k $ both remain in $X$, and that f is continuously differentiable, we conclude that $\{u_{6,k}\}$ is bounded. Taking limits in \eqref{eqn:ekdeltakasymp} and recalling that
$\lim_{k\to\infty}J_k=\infty$, we conclude that
$\lim_{k\to\infty}\|\delta_k\|=0$.} 
} 
\end{proof}
\ser{\begin{proposition}\label{pro:bound_g_new}\em
\fyr{Let Assumption~\ref{assump:nonconvex_main} hold and let $\{\hat{x}_k \}$ be generated by Algorithm~\ref{alg:IPR-VFISTA}. Then,} \ssg{for all  $k\geq 0$}
\begin{align}\label{eq:bound_g_new_forasym}
\ser{\|G_{1/\hat{\gamma}}(\hat{x}_k)\|^2
\le \tfrac{16 (f(\hat{x}_k)-f(\hat{x}_{k+1}))}{\hat{\gamma}} 
+{ \tfrac{36}{L_f\hat{\gamma}^3}}\|\delta_k\|^2 +\tfrac{16}{\hat{\gamma}}\|e_k-\delta_k\|
\hat{B}_f
+\tfrac{16}{\hat{\gamma}^2}\|e_k\|^2.}
\end{align}
\end{proposition}}
\begin{proof}
\noindent  \ser{By Definitions~\ref{def:residual} and~\ref{def:eanddelta}, as in \fyr{Lemma~5.6 of~\cite{samadi2025improved}}, we obtain
\begin{align} \label{eq:new_g_bound}
    \tfrac{\hat{\gamma}^2}{2} \|G_{1/\hat{\gamma}}(\hat{x}_k)\|^{\fyr{2}} \leq {\|\hat{x}_{k+1} - \hat{x}_k \|}^2+ {\| \delta_k\|}^2.
\end{align} 
Next, by the projection theorem, for any $k\geq 0$, $\langle
\Pi_{X_{\bar h}^*}[z_k]-\Pi_{X_{\bar h}^*}[\hat x_k],
z_k-\Pi_{X_{\bar h}^*}[z_k]
\rangle\geq 0.$
\fyr{Substituting $\Pi_{X_{\bar h}^*}[z_k]=\hat{x}_{k+1}-\delta_k$, $\Pi_{X_{\bar h}^*}[\hat x_k]=\hat{x}_{k}-e_k$, and $z_k = \hat{x}_k - \hat{\gamma}\nabla f(\hat{x}_k)$, the preceding inequality yields \begin{align*}&{\left\langle\hat{x}_{k+1}-\hat{x}_k,
\hat{x}_k-\hat{\gamma}\nabla f(\hat{x}_k)-\hat{x}_{k+1}\right\rangle}+{\left\langle\hat{x}_{k+1}-\hat{x}_k,\delta_k\right\rangle}
 \\
&+
{\left\langle e_k-\delta_k,
\hat{x}_k-\hat{\gamma}\nabla f(\hat{x}_k)-\hat{x}_{k+1}\right\rangle}
+{\left\langle e_k-\delta_k,\delta_k\right\rangle}
\geq 0.\end{align*}
\noindent  Denote these four inner products, from left to right, by Terms 1--4. Next, we expand Term 1. For Term 2, we apply Young's inequality
$\langle u,v\rangle \leq \tfrac{1}{2\theta}\|u\|^2+\tfrac{\theta}{2}\|v\|^2$,
where $\theta>0$ is arbitrary. For Term 4, we use the inequality
$\langle u,v\rangle \leq \tfrac{1}{2}\|u\|^2+\tfrac{1}{2}\|v\|^2$.
Moreover,
\begin{align*}
\mathrm{Term\ 3}
&=\langle e_k-\delta_k,\hat{x}_k-\hat{x}_{k+1}\rangle
+\langle e_k-\delta_k,-\hat{\gamma}\nabla f(\hat{x}_k)\rangle \\
&\leq
\tfrac{1}{2}\|e_k-\delta_k\|^2
+\tfrac{1}{2}\|\hat{x}_k-\hat{x}_{k+1}\|^2
+\hat{\gamma}\|e_k-\delta_k\|\|\nabla f(\hat{x}_k)\|,
\end{align*}
where the first two terms are bounded using Young's inequality and the last term is bounded using the Cauchy--Schwarz inequality. Therefore, we obtain
\begin{align*}
&\underbrace{
\langle -\hat{\gamma}\nabla f(\hat{x}_k),  (\hat{x}_{k+1}-\hat{x}_k)\rangle-\|\hat{x}_{k+1}-\hat{x}_k\|^2 
}_{=\text{Term 1}}
+\underbrace{
\tfrac{1}{2\theta}\|\hat{x}_{k+1}-\hat{x}_k\|^2+\tfrac{\theta}{2}\|\delta_k\|^2
}_{\text{Term 2}\le}
\\\notag &+\underbrace{
\tfrac{1}{2}\|e_k -\delta_k\|^2 + \tfrac{1}{2}\| \hat{x}_k-\hat{x}_{k+1}\|^2 +  \fyr{\hat{\gamma}}\|e_k -\delta_k\| \|\nabla f(\hat{x}_{k}) \|
}_{\text{Term 3}\le}\\\notag &+\underbrace{
\tfrac{1}{2}\|e_k-\delta_k\|^2+\tfrac{1}{2}\|\delta_k\|^2
}_{\text{Term 4}\le}\geq 0.
\end{align*} 
Choosing $\theta:=\tfrac{2}{\hat{\gamma}L_f}$ yields}
\begin{align*}
\langle\nabla f(\hat{x}_k),\hat{x}_{k+1}-\hat{x}_k\rangle
&\leq
\tfrac{L_f\hat{\gamma}-2}{4\hat{\gamma}}
\|\hat{x}_{k+1}-\hat{x}_k\|^2
+\tfrac{2\|\delta_k\|^2}{L_f\hat{\gamma}^2}
+\|e_k-\delta_k\|\|\nabla f(\hat{x}_k)\|
+\tfrac{\|e_k\|^2}{\hat{\gamma}},
\end{align*}
where we used \fyr{$\theta\geq4$ and $\tfrac{\theta+3}{2}\leq\theta$, in view of $\hat{\gamma}\leq\tfrac{1}{2L_f}$}.}
\ser{From the $L_f$-smoothness property of $f$, we have $f(\hat{x}_{k+1})
\leq f(\hat{x}_k)
+ \langle \nabla f(\hat{x}_k), \hat{x}_{k+1} - \hat{x}_k\rangle
+ \tfrac{L_f}{2}\|\hat{x}_{k+1} - \hat{x}_k\|^2.$
From the two preceding relations, and noting  that
$\tfrac{1}{8\hat{\gamma}}
\le(\tfrac{1}{2\hat{\gamma}}-\tfrac{3L_f}{4})$
in view of the assumption $\hat{\gamma}\le \tfrac{1}{2L_f}$, }
\ser{and  adding $\tfrac{1}{8\hat{\gamma}}\|\delta_k\|^2$ to both sides, 
\begin{align*}
\tfrac{\|\hat{x}_{k+1}-\hat{x}_k\|^2+\|\delta_k\|^2}{8\hat{\gamma}}
\leq
f(\hat{x}_k)-f(\hat{x}_{k+1})
+\fyr{\tfrac{16+L_f\hat{\gamma}}
{8L_f\hat{\gamma}^2}}\|\delta_k\|^2
+\|e_k-\delta_k\|\|\nabla f(\hat{x}_k)\|
+\tfrac{\|e_k\|^2}{\hat{\gamma}}.
\end{align*}
\fyr{Using \eqref{eq:new_g_bound} and $\tfrac{1}{2\hat{\gamma}} \leq \tfrac{1}{L_f\hat{\gamma}^2}$, from the preceding inequality,} we obtain the result.}
\end{proof} 
In the next proposition, we \fyr{derive} upper bounds for $\|e_k\|$ and ${\|\delta_k\|}^2$.
\begin{proposition}\label{prop:nonconvexedeltabound}\em
 \ser{\fyr{Let Assumption~\ref{assump:nonconvex_main} hold and let} $\{\hat{x}_k \}$ be generated by Algorithm~\ref{alg:IPR-VFISTA}.} \ser{Suppose $ X^*_{\bar{h}}$ is $\alpha$-weak sharp minima of order $m\geq 1$.  \fyr{Then, for} all $J_k$ satisfying {$16+16\tfrac{L_h}{\bar{\eta}}  \leq (\tfrac{ {J_k}}{{\ln(J_k)}}  )^2$}, \fyr{we have}
 }
 \begin{align*}
 & \|e_k\|   \leq     \sqrt[\ser{m}]{{\alpha^{-1}} \left(\newchr{\tfrac{(\ln{{(J_{k-1} ))^2}}}{ {J_{k-1}}^2}} 16 \ ({T_{1,k-1}}+T_{3,k-1}) + \tfrac{T_{2,k-1}}{J_{k-1}^4}\right)}  \\                                
   \fyr{\text{and} \quad } &{\|\delta_k\|}^2
\leq  \tfrac{T_{4,k}}{J_k^{4}}
+\tfrac{T_{5,k}}{J_k^{2} \ln(J_k)} 
+ 2 \ \alpha^{\ser{\tfrac{-2}{m}}} \left( \left( \tfrac{16 (T_{1,k-1} + T_{3,k-1}) (\ln(J_{k-1}))^2}{J_{k-1}^2} + \tfrac{T_{2,k-1}}{J_{k-1}^4} \right) \right. \\
 &  \left. \times \left( \tfrac{ 16 \ T_{1,k} (\ln(J_k))^2}{J_k^6} 
  + \tfrac{ T_{2,k}}{J_k^4} 
  + \tfrac{16 \ T_{3,k} (\ln(J_k))^2}{J_k^2} \right) \right)^{\tfrac{1}{\ser{m}}}    \\\notag & + \tfrac{2\ \hat{\gamma} \| \nabla f(\hat{x}_k) \|}{\sqrt[\ser{m}]{ \alpha}} \sqrt[\ser{m}]{ \tfrac{16 \ T_{1,k} (\ln(J_k))^2}{J_k^6} 
  + \tfrac{ T_{2,k}}{J_k^4} 
  + \tfrac{16 \ T_{3,k} (\ln(J_k))^2}{J_k^2} }.
\end{align*}
 \end{proposition}
 \begin{proof}
\fyr{From Algorithm~\ref{alg:IPR-VFISTA}, we have $\hat{x}_{k} = x_{k-1,J_{k-1}}$ for $k \geq 1$.  Invoking  \fyr{Corollary}~\ref{Cor:R-VFISTA} (1.ii) with $p=3$ and from the} nonnegativity of  $f_{s,k}(x) $, we obtain
  \begin{align*}
  & 0 \leq \bar{h}(\hat{x}_k) -{\bar{h}^*}  \leq  \tfrac{\newch{  16 
  T_{1,k-1}(\ln({J_{k-1}}))^2}}{{J_{k-1}}^{6}}
  +\tfrac{ T_{2,k-1}}{{J_{k-1} }^{4}} 
  +   \tfrac{16 T_{3,k-1}(\ln{(J_{k-1})})^2}{{{J_{k-1}}}^2}.
 \end{align*}
\fyr{The bound on $\|e_k\|$ is implied by using} $\|e_k\| =  \operatorname{dist}(\hat{x}_k , X^*_{\bar{h}})  \leq \sqrt[\ser{m}]{\alpha^{-1} {(\bar{h}(\hat{x}_k) - \bar{h}^*)}}$. \fyr{Next, we show the second inequality in (1).} Consider  Definition~\ref{def:eanddelta}. In view of Algorithm~\ref{alg:IPR-VFISTA}, we have $\hat{x}_{k+1}=x_{k,J_k}$ and also, recall that $x^*_{{f_{\fyr{s,k}}}}=\Pi_{{X}^*_{\bar{h}}}[z_k]$. Thus, we have $\delta_k =x_{k, J_k} - x^*_{{f_{s,k}}}$.
\fyr{Invoking Corollary}~\ref{Cor:R-VFISTA}~(1.iii) with $p=3$ and \ser{$m\geq 1$}, 
  \begin{align}\label{eq:delta_2}
  { {\|\delta_k\|}^2}
&\leq  \tfrac{T_{4,k}}{\newchr{{{J_k}}^{4}}}
+\ch{\tfrac{T_{5,k}}{\newchr{{{J_k}}^{2}}\ln({J_k})}} 
 +  \ch{\tfrac{2 \|\nabla {f_{\fyr{s,k}}}(x^*_{{f_{\fyr{s,k}}}})\|}{\sqrt[\ser{m}]{\alpha}} }  {\sqrt[\ser{m}]{\tfrac{ 16\ T_{1,k}\newch{(\ln({J_{k}}))^2}}{{J_{k}}^{6}} 
  +\tfrac{ T_{2,k}}{{J_{k} }^{4}} 
  +  \newch{ \tfrac{16\ T_{3,k}(\ln{(J_{k})})^2}{{{J_{k}}}^2}} }}.
  \end{align}
Now, consider \eqref{prob:proj}. We have $\|\nabla {f_{s,k}} (x) \| = \| z_k - x\|$. Then, by substituting $x=x^*_{{f_{s,k}}}=\Pi_{{X}^*_{\bar{h}}}[z_k]$, we get  $\|\nabla {f_{s,k}} (x^*_{{f_{s,k}}}) \|= \| z_k - \Pi_{{X}^*_{\bar{h}}}[z_k]\|$. In view of Algorithm~\ref{alg:IPR-VFISTA}, we have $z_k = \hat{x}_k - \hat{\gamma}\nabla f(\hat{x}_k)$. Invoking the triangle inequality, we obtain
\begin{align}\label{eq:nablaH}
&\|\nabla {f_{s,k}} (x^*_{{f_{s,k}}}) \| 
 \leq  \| {\hat{x}_k- \Pi_{{X}^*_{\bar{h}}}[\hat{x}_k]}\| +\hat{\gamma} \| \nabla f(\hat{x}_k)\| 
= \|e_k\| + \hat{\gamma} \| \nabla f(\hat{x}_k)\|.
\end{align}
Using \eqref{eq:nablaH}, \eqref{eq:delta_2}, and the bound on $\|e_k\|$, we obtain the second inequality.
 \end{proof}
\fyr{In the next result, under a \ssg{weak sharp minimality}, we show that generated iterate by Algorithm~\ref{alg:IPR-VFISTA} is bounded.} 
\ser{\begin{lemma}\label{lemma:defien_r_k_combined}\em
Suppose Assumption~\ref{assump:nonconvex_main} holds. \fyr{Assume further that} $X^*_{\bar h}$ is bounded and $\alpha$-weakly sharp of order $m=2$. Let $\{\hat{x}_k\}$ be generated by Algorithm~\ref{alg:IPR-VFISTA} \fyr{and let $J_k:=\lceil (k+1)^a\rceil$ for any $k\geq 0$ where $a\ge2$.} Then, \fyr{the sequence} $\{\hat{x}_k\}$ is bounded.
\end{lemma}}
\begin{proof}
    \noindent  In view of Definition~\ref{def:eanddelta}, we have
 $\|e_k\|= \|\hat{x}_k - \Pi_{{X}^*_{\bar{h}}}[\hat{x}_k]\|. $
\fyr{Invoking} the triangle inequality, we obtain
$ \|\hat{x}_k\| \leq D_{X^*_{\bar{h}}} + \|e_k\|.
$ Next,  squaring \ser{both sides} and invoking Proposition~\ref{prop:nonconvexedeltabound} \ser{yield ${\|\hat{x}_k\|}^2  \leq 2  D_{X^*_{\bar {h}}}^2 + \tfrac{2}{\alpha} (\newchr{\tfrac{(\ln{{(J_{k-1} ))^2}} 16  ({T_{1,k-1}}+T_{3,k-1})}{ {J_{k-1}}^2}}  + \tfrac{T_{2,k-1}}{J_{k-1}^4}).$} \fyr{Next, invoking Lemma~\ref{lemma:bounded_T} (i) and substituting the expression for $J_k$, for all $k\ge1$,}
      \begin{align*} \ser{ 
      {\|\hat{x}_k\|}^2  \leq 2D_{X_{\bar{h}}^*}^{\ssg{2}} +  \tfrac{32 a^2 }{\alpha } \left( c_1+c_3\right)+\tfrac{2 B_{\bar{h}}}{\alpha} +\tfrac{64 a^2 c_2 }{\alpha k^{2a-2}} {\| \hat{x}_{k-1}\|}^2  \leq
   c_4 + \tfrac{c_5}{k^{2a-2}} {\| \hat{x}_{k-1}\|}^2,}
   \end{align*}
where \ser{$c_1, c_2, \ \text{and} \ c_3$ are given in Lemma~\ref{lemma:bounded_T},} $c_4 :=  2D_{X_{\bar{h}}^*}^{\ssg{2}} +  \tfrac{32 a^2 }{\alpha } \left( c_1+c_3\right)+\tfrac{2 B_{\bar{h}}}{\alpha}$\ser{,} and {$\ser{c_5:=}\tfrac{64 a^2 c_2 }{\alpha} $}. \fyr{The result follows by Lemma~\ref{lemma:general_boundedness}.}
\end{proof}
\fyr{The following theorem establishes the asymptotic convergence of Algorithm~\ref{alg:IPR-VFISTA}.
\begin{theorem}[asymptotic stationarity for $\text{IPR-VFISTA}_\text{nc}$] \label{theorem:asym_non_cons}\em
    \ser{Consider the problem in \eqref{prob:centr-nonconvex} under Assumption~\ref{assump:nonconvex_main}. Let $\{\hat{x}_k\}$ be generated by Algorithm~\ref{alg:IPR-VFISTA}.} Assume that one of the following holds.\\
      (1) \fyr{$\omega_h=\mathbb{I}_X$, where $X\subset\mathbb{R}^n$ is a nonempty compact convex set. Also, $\{J_k\}$ is an increasing sequence of positive integers satisfying $J_k\to\infty$.}\\
(2) $X^*_{\bar h}$ is bounded and $\alpha$-weakly sharp of order $2$. Also, $J_k:=\lceil (k+1)^a\rceil$ for any $k\geq 0$\ssg{,} where $a\ge2$.\\
 Then, $\liminf_{k\to\infty} \|G_{1/\hat{\gamma}}(\hat{x}_k)\| = 0$.
\end{theorem}}
\begin{proof}
 \fyr{From \eqref{eq:bound_g_new_forasym}, we obtain, for any $k\geq 0$, 
\begin{align}
\label{eq:bound_g_new_forasym_re}
\|G_{1/\hat{\gamma}}(\hat{x}_k)\|^2
&\le \tfrac{16 (f(\hat{x}_k)-f(\hat{x}_{k+1}))}{\hat{\gamma}} 
+{ \tfrac{36}{L_f\hat{\gamma}^3}}\|\delta_k\|^2 +\tfrac{16}{\hat{\gamma}}(\|e_k\|+\|\delta_k\|)
\hat{B}_f
+\tfrac{16}{\hat{\gamma}^2}\|e_k\|^2.
\end{align}
Suppose, by contradiction, 
$\liminf_{k\to\infty}\|G_{1/\hat{\gamma}}(\hat{x}_k)\| > 0$.
Then there exist $\epsilon>0$ and $K_{\fyr{\epsilon}}\geq 0$ such that
$\|G_{1/\hat{\gamma}}(\hat{x}_k)\|^2 \geq \epsilon$
for all $k\geq K_{\fyr{\epsilon}}$. Next, we consider the two settings. First, let conditions in (1) hold. By Lemma~\ref{lemma:limit_dist_nonconvex_cons},  $\lim_{k \to \infty} \|e_k\|=\lim_{k \to \infty} \|\delta_k\| = 0.$ Thus, there exists $\hat{K}_{\epsilon} \geq K_{\epsilon}$ such that for any $k\geq \hat{K}_{\epsilon}$, we have ${ \tfrac{36}{L_f\hat{\gamma}^3}}\|\delta_k\|^2 +\tfrac{16}{\hat{\gamma}}(\|e_k\|+\|\delta_k\|)
\hat{B}_f
+\tfrac{16}{\hat{\gamma}^2}\|e_k\|^2 \leq \tfrac{\epsilon}{2}$. From \eqref{eq:bound_g_new_forasym_re}, we obtain 
$\epsilon \le \tfrac{16 (f(\hat{x}_k)-f(\hat{x}_{k+1}))}{\hat{\gamma}}  +\tfrac{\epsilon}{2}$, for all $k\geq \hat{K}_{\epsilon}$. Summing over $k=\hat{K}_{\epsilon},\ldots, N-1$, where $N>\hat{K}_{\epsilon} $, we obtain $f(\hat{x}_{N}) \leq f(\hat{x}_{\hat{K}_{\epsilon}}) - \tfrac{(N-\hat{K}_{\epsilon}) \hat{\gamma} \epsilon}{32}.$ Taking limits on both sides, we obtain $f(\hat{x}_{N}) \to -\infty$, contradicting Assumption~\ref{assump:nonconvex_main}. Hence, if the conditions in (1) hold, then $\liminf_{k\to\infty} \|G_{1/\hat{\gamma}}(\hat{x}_k)\| = 0$. Next, let conditions in (2) hold. From \ssg{Lemma}~\ref{lemma:defien_r_k_combined}, the sequence $\{\hat{x}_k\}$ is   bounded. Applying Lemma~\ref{lemma:bounded_T}, together with Proposition~\ref{prop:nonconvexedeltabound}, $\lim_{k \to \infty} \|e_k\|=\lim_{k \to \infty} \|\delta_k\| = 0.$ The proof can be completed by following analogous steps as in setting (1).}
\end{proof}
\subsection{\texorpdfstring{\ser{\fyr{Nonasymptotic} convergence rates for \ssg{$\text{IPR-VFISTA}_{\text{nc}}$}}}{Nonasymptotic convergence rates for IPR-VFISTA-nc}}\label{subsection:4.3}
\fyr{In the following, we derive a nonrecursive error bound for the residual mapping.}
\ser{\begin{proposition}\em\label{pro:bound_g_new_fort_theorem}
 Consider problem~\eqref{prob:centr-nonconvex} under Assumption~\ref{assump:nonconvex_main}, and let $\{\hat{x}_k\}$ be generated by Algorithm~\ref{alg:IPR-VFISTA}. \fyr{Define $k^*$ such that} \ssg{$\|G_{1/\hat{\gamma}}(\hat{x}_{k^*})\| := \allowbreak \min_{k=\lfloor K/2\rfloor,\ldots,K-1} \allowbreak \|G_{1/\hat{\gamma}}(\hat{x}_k)\|.$} Then, \fyr{for any $K\geq 1$,}
\begin{align}\label{eq:sumGfirst_withoutgamma}
\|G_{1/\hat{\gamma}}(\fyr{\hat{x}_{k^*}})\|^2
&\leq
\tfrac{32\big(\fyr{f(\hat{x}_{\lfloor K/2 \rfloor})-f(\hat{x}_K)}\big)}
{K\hat{\gamma}}
+\fyr{\ssg{\tfrac{1}{K}
\textstyle\sum_{k=\lfloor K/2 \rfloor}^{K-1}
\left(
\tfrac{32\hat{B}_f}{\hat{\gamma}}
(\|e_k\|+\|\delta_k\|)
\right.}}
\\
&
\fyr{\ssg{\left.
+\tfrac{32\|e_k\|^2}{\hat{\gamma}^2}
+\tfrac{72\|\delta_k\|^2}{L_f\hat{\gamma}^3}
\right)}}.
\notag
\end{align}
\end{proposition}
\begin{proof}
\fyr{Summing} \ser{both sides} of \eqref{eq:bound_g_new_forasym} over $k = \left\lfloor {K}/{2} \right\rfloor, \ldots, K-1$, we obtain
\begin{align}\label{eq:G_first_new}
\textstyle\sum_{k=\lfloor K/2 \rfloor}^{K-1} \|G_{1/\hat{\gamma}}(\hat{x}_k)\|^2
&\leq \tfrac{16}{\hat{\gamma}} \big( f(\hat{x}_{\lfloor K/2 \rfloor}) - f(\hat{x}_K) \big)
+ \tfrac{36}{L_f \hat{\gamma}^3} \textstyle\sum_{k=\lfloor K/2 \rfloor}^{K-1} \|\delta_k\|^2  \\\notag
& + \tfrac{16\hat{B}_f}{\hat{\gamma}} \textstyle\sum_{k=\lfloor K/2 \rfloor}^{K-1} \|e_k - \delta_k\|
+ \tfrac{16}{\hat{\gamma}^2} \textstyle\sum_{k=\lfloor K/2 \rfloor}^{K-1} \|e_k\|^2.
\end{align}
From the definition of $\fyr{\hat{x}_{k^*}}$ and noting that $\tfrac{K}{2} \leq K - \lfloor K/2 \rfloor$, we have
$\tfrac{K}{2} \|G_{1/\hat{\gamma}}(\fyr{\hat{x}_{k^*}})\|^2
\leq \textstyle\sum_{k=\lfloor K/2 \rfloor}^{K-1} \|G_{1/\hat{\gamma}}(\hat{x}_k)\|^2.$ \fyr{Hence, the result follows from \eqref{eq:G_first_new}}.
\end{proof}}
\ser{In the next theorem, we provide convergence rate statements for $\text{IPR-VFISTA}_\text{nc}$.} \fyr{We will utilize the following terms.
\begin{definition}\label{def:thm_nonasymp_nc}\em     Let  us define \ssg{the following terms}, for $K\geq 1$ and $a>1$, $c_{e,K}:=    16 \alpha^{-1}   \left(\ln(\lceil K^a\rceil )\right)^2   (M_1 + M_2+M_3) $ and $c_{\delta,K} := ( \tfrac{M_4 + M_5}{a} ) +  2 {c_{e,K}} +2\hat{\gamma}B_f\sqrt{{c_{e,K}}}$, where $M_i$, for $i \in \{1,2,3, 4,5 \}$, are given by {Lemma~\ref{lemma:bounded_T}}. \fyr{We denote  by $\bar{K}$ the smallest integer $K\geq 2$ satisfying $  \left(\frac{\lceil\left(\left\lfloor {K}/{2} \right\rfloor+1\right)^a\rceil}{\ln\left(\lceil\left(\left\lfloor {K}/{2} \right\rfloor+1\right)^a\rceil\right)}\right)^2 \geq 16+16\frac{L_h}{\bar{\eta}}$, where $\bar{\eta}>0$ is an arbitrary scalar.}
\end{definition}}
\begin{theorem}[error bounds for $\text{IPR-VFISTA}_\text{nc}$]\label{theorem:nonc_all}\em
     \fyr{Let $\{\hat{x}_k\}$ be generated by Algorithm~\ref{alg:IPR-VFISTA} and let Assumption~\ref{assump:nonconvex_main} hold. \ssg{ Assume that $J_k:=\lceil (k+1)^a\rceil$ for any $k\geq 0$ where $a>1$.}  Define $\|G_{{1}/{\hat{\gamma}}}(\fyr{\hat{x}_{k^*}})\| = \min_{k = \left\lfloor {K}/{2} \right\rfloor, \ldots, K-1}\|G_{{1}/{\hat{\gamma}}}(\hat{x}_k)\|$. \ssg{Consider Definition~\ref{def:thm_nonasymp_nc}.} Consider the following settings.}

\medskip

\noindent \fyr{(1) Suppose $X^*_{\bar{h}}$ is $\alpha$-weakly sharp of order $m \in (1,2)$ \fyr{and let $a> m$. Moreover, suppose $\omega_h=\mathbb{I}_X$ for some nonempty compact convex set $X\subset\mathbb{R}^n$}. Then, the following results hold.} 

\noindent \fyr{{(1.i)}  [infeasibility bound]  $
  \operatorname{dist}(\hat{x}_K,X^*_{\bar{h}})  \leq      4\sqrt[m]{ \alpha^{-1}(M_1+M_2+M_3)} \fyr{\tfrac{\ln{{(\lceil K^a\rceil )}}}{ \lceil K^a\rceil}}$ \fyr{for any $K\geq \bar{K}$}.\\
\noindent{{(1.ii)}} [residual bound]  \fyr{For any $K\geq \bar{K}$,}
 \begin{align*}
\|G_{1/\hat{\gamma}}(\hat{x}_{\fyr{k^*}})\|^2
&\le
\left(\tfrac{32(\hat{D}_f-\hat{C}_f)}{\hat{\gamma}}\right)\frac{1}{K}
+\ssg{\left(\tfrac{(144)(3^{{2a}/{m}})c_{\delta,K}^{2/m}}
{3L_f\hat{\gamma}^3}+
\tfrac{64\hat{B}_f(3^{{2a}/{m}})c_{e,K}^{1/m}}
{3\hat{\gamma}} \right)
\frac{1}{K^{\tfrac{2a}{m}}}} \\
&
+ \ssg{\left( 
\tfrac{a(32\hat{B}_f)(3^{{a}/{m}})}
{3m\hat{\gamma}( {a}/{m}-1)}
c_{\delta,K}^{1/m}\right)
\frac{1}{K^{\tfrac{a}{m}}}}
+\left( 
\tfrac{(128) (3^{ {4a}/{m}})c_{e,K}^{2/m}}
{9\hat{\gamma}^2}\right)
\frac{1}{K^{\tfrac{4a}{m}}}.
\end{align*}
\noindent (1.iii) [total complexity] The total iteration complexity is 
$\mathcal{O}(\epsilon^{-(a+1)})$, where $a > m$ and $\epsilon>0$ is an arbitrary scalar such that ${\|G_{{1}/{\hat{\gamma}}} (\hat{x}_{k^*})\|}^2 \leq \epsilon$.}

\medskip 

\noindent \fyr{(2)  Suppose $X^*_{\bar{h}}$ is a bounded set and $\alpha$-weakly sharp of order $m=2$ and let $a> 2$. Then, (1.i), (1.ii), and (1.iii) hold.}
\end{theorem}
 \begin{proof}
 \fyr{(1.i)  Since $X$ is compact and $\hat{x}_{k}:=x_{k-1,J_{k-1}}\in X$ for all $k\geq 1$, the sequence $\{\hat{x}_k\}$ is bounded. From Definition~\ref{def:thm_nonasymp_nc} and $J_k:=\lceil (k+1)^a\rceil$, the conditions of Proposition~\ref{prop:nonconvexedeltabound} hold for any $K\geq \bar{K}$. By Lemma~\ref{lemma:bounded_T} (ii), together with Proposition~\ref{prop:nonconvexedeltabound}, we obtain $$
  \operatorname{dist}(\hat{x}_K,X^*_{\bar{h}}) = \|e_K\| \leq     4\sqrt[m]{ \alpha^{-1}(M_1+M_2+M_3)}{\tfrac{\ln{{(J_{K-1} )}}}{ {J_{K-1}}}} .$$ Substituting $J_K = \lceil(K+1)^a\rceil$, we obtain the result in~(1.i).}\\
 \noindent\fyr{(1.ii) From Proposition~\ref{prop:nonconvexedeltabound} and Definition~\ref{def:thm_nonasymp_nc}, for any $ \left\lfloor {K}/{2} \right\rfloor \leq k \leq K$ and $K\geq \bar{K}$, we obtain ${\|e_k\|}^2 \leq (\tfrac{{c_{e,K}}}{ k^{2a}})^{\tfrac{2}{m}}$ and $ {\|\delta_k\|}^2 \leq (\tfrac{{c_{\delta,K}}}{{k}^{a}})^{\tfrac{2}{m}}.$ This implies that
 \begin{align*}
\textstyle\sum_{k=\lfloor \ssg{\tfrac{K}{2}} \rfloor}^{K-1}\|e_k\|^2
\le (c_{e,K})^{\tfrac{2}{m}}\textstyle\sum_{k=\lfloor \ssg{\tfrac{K}{2}} \rfloor}^{K-1}\tfrac{1}{k^{\ssg{\tfrac{4a}{m}}}}\ssg{,} \
\textstyle\sum_{k=\lfloor \ssg{\tfrac{K}{2}} \rfloor}^{K-1}\|\delta_k\|^2
\le \ssg{(c_{\delta,K}\ssg)}^{\ssg{\tfrac{2}{m}}}\textstyle\sum_{k=\lfloor \ssg{\tfrac{K}{2} }\rfloor}^{K-1}\tfrac{1}{k^{\ssg{\tfrac{2a}{m}}}}.
\end{align*}
Invoking~\cite[Lemma 9]{yousefian2017smoothing} and using $a>m$ and $\lfloor K/2\rfloor \ge K/3$ for $K\ge 2$, we obtain 
$$\textstyle
\sum_{k=\lfloor K/2\rfloor}^{K-1}\frac{1}{k^{b}}
\leq
\left(\frac{b}{b-1}\right)3^{b-1}
\frac{1}{K^{b-1}},
\qquad \hbox{where } b:=\frac{pa}{m} \hbox{ and } p\geq 1 .$$
Thus, applying the preceding bound with $p :=1,2,4$, we obtain
 \begin{align*}
&\tfrac{1}{K}\textstyle\sum_{k=\lfloor K/2 \rfloor}^{K-1}\left(\tfrac{32\hat{B}_f}{\hat{\gamma}}(\|e_k\|+\|\delta_k\|)
+\tfrac{32\|e_k\|^2}{\hat{\gamma}^2}+
\tfrac{72\|\delta_k\|^2}{L_f\hat{\gamma}^3}\right)\\
&\leq
\tfrac{1}{K}\textstyle\sum_{k=\lfloor K/2\rfloor}^{K-1}\left(\left(\tfrac{32\hat{B}_f c_{e,K}^{1/m}}{\hat{\gamma}}\right)\frac{1}{k^{2a/m}}
+\left(\tfrac{32\hat{B}_f c_{\delta,K}^{1/m}}{\hat{\gamma}}\right) \frac{1}{k^{a/m}}
+\left(\tfrac{32c_{e,K}^{2/m}}{\hat{\gamma}^2}\right)
\frac{1}{k^{4a/m}} \right.
\\
&
\left.+\ssg{\left(\tfrac{72c_{\delta,K}^{2/m}}{L_f\hat{\gamma}^3}\right)\frac{1}{k^{2a/m}}}\right)\leq
\left( 
\tfrac{a(32\hat{B}_f)(3^{{a}/{m}})}
{3m\hat{\gamma}( {a}/{m}-1)}
c_{\delta,K}^{1/m}\right)
\frac{1}{K^{\tfrac{a}{m}}} 
+\left( 
\tfrac{(128) (3^{ {4a}/{m}})c_{e,K}^{2/m}}
{9\hat{\gamma}^2}\right)
\frac{1}{K^{\tfrac{4a}{m}}}
\\\notag &\ssg{+\left(\tfrac{(144)(3^{{2a}/{m}})c_{\delta,K}^{2/m}}
{3L_f\hat{\gamma}^3}+
\tfrac{64\hat{B}_f(3^{{2a}/{m}})c_{e,K}^{1/m}}
{3\hat{\gamma}} \right)
\frac{1}{K^{\tfrac{2a}{m}}}}.
\end{align*}
Substituting the preceding bound into~\eqref{eq:sumGfirst_withoutgamma}, we obtain the result in (1.ii).}\\
\noindent \fyr{(1.iii) From part (1.ii) and the assumption $a>m$, we have
\ser{$\|G_{{1}/{\hat{\gamma}}}(\fyr{\hat{x}_{k^*}})\|^2=\mathcal{O}(1/K)$}.
Thus, an $\epsilon$-stationary point is obtained with $K=\mathcal{O}(\epsilon^{-1})$ number of outer iterations. Since $J_k=\lceil(k+1)^a\rceil$, the total iteration complexity is
\ser{$\textstyle\sum_{k=0}^{\mathcal{O}(\epsilon^{-1})}\lceil(k+1)^a\rceil=\mathcal{O}(\epsilon^{-(a+1)})$.}}\\
\noindent \fyr{(2) By Lemma~\ref{lemma:defien_r_k_combined}, the sequence $\{\hat{x}_k\}$ is bounded. The remainder of the proof follows analogously to the proof of part (1).}
 \end{proof}

\begin{remark}\label{rem:nonconvex_comparison1}
\fyr{(i) Theorem~\ref{theorem:asym_non_cons} establishes an asymptotic stationarity guarantee for the problem in~\eqref{prob:centr-nonconvex} with nonconvex $f$, without assuming the \ssg{weak sharp minimality}. (ii) \ser{Theorem~\ref{theorem:nonc_all} complements Theorem~\ref{theorem:asym_non_cons}} by also providing explicit stationarity and infeasibility guarantees. (iii) For the SBO subclass of the nonconvex optimization problem with variational inequality (VI) constraints studied in \cite{samadi2025improved}, \ssg{Theorem~\ref{theorem:nonc_all}} improves the corresponding complexity bound from $\mathcal{O}(\epsilon^{-8})$ to nearly $\mathcal{O}(\epsilon^{-3})$ under the \ssg{weak sharp minimality} with $m=2$, by utilizing acceleration in the inner loop together with a constant outer-loop stepsize in $\mathrm{IPR\text{-}VFISTA}_{\mathrm{nc}}$.} 
\end{remark}
\begin{remark}\label{rem:nonconvex_comparison2}
\fyr{We next compare the proposed method with conditional gradient schemes in~\cite{jiang2023conditional,giang2026conditional}. CG-BiO~\cite{jiang2023conditional} requires an initial point $x_0$ satisfying $h(x_0)-h^*\leq\epsilon_h/2$, whose computation requires $\mathcal{O}(\epsilon_h^{-1})$ conditional-gradient iterations. Corollary~1~(ii) of \cite{jiang2023conditional} under \fyr{weak sharp minimality} of $X^*_h$ of order $m\geq 1$ establishes $|\mathcal{G}(x_{k^*})|\leq\epsilon_f$ together with $h(x_{k^*})-h^*\leq\epsilon_h$, where $\epsilon_h=\Theta(\epsilon_f^m)$ and $\mathcal{G}(x):=\max_{s\in X_h^*} \langle\nabla f(x),x-s\rangle$ is the Frank--Wolfe \fyr{stationarity} measure \fyr{(cf.~\cite{jaggi2013revisiting})}. The total complexity of \fyr{CG-BiO} is $\mathcal{O}(\epsilon_h^{-1}) +\mathcal{O}(\epsilon_f^{-(m+1)})$.  For $\hat{\gamma}\leq 1/L_f$ and $x\in X_{\bar{h}}^*$, $\|G_{1/\hat{\gamma}}(x)\|^2\leq\mathcal{G}(x)/\hat{\gamma}$ (cf.~\cite[(37)]{bomze2021frank}). Theorem~\ref{theorem:nonc_all} gives $\mathcal{O}(\epsilon^{-(m+1+\delta)})$ for $m\in(1,2]$ for an arbitrary $\delta>0$. Thus, the iteration complexity of $\mathrm{IPR\text{-}VFISTA}_{\mathrm{nc}}$ nearly matches that of CG-BiO, while additionally providing an asymptotic stationarity guarantee. 
For $m=2$, Corollary~1 (ii) uses
$\epsilon_h=\mathcal{O}(\epsilon^2)$, which implies
$\mbox{dist}^2(x_{k^*},X_h^*)=\mathcal{O}(\epsilon^2)$. Selecting
$a:=2+\delta$ in Theorem~\ref{theorem:nonc_all} yields, at the end of
the $K$th inner loop, $\operatorname{dist}^2(\hat{x}_K,X_{\bar h}^*)
=\mathcal{O}(\epsilon^{4+2\delta})$. Therefore, \fyr{$\mathrm{IPR\text{-}VFISTA}_{\mathrm{nc}}$ has} a sharper infeasibility guarantee at completed outer iterations, whereas \fyr{CG-BiO} maintains its prescribed lower-level
objective-gap bound throughout the main phase. The IR-FSCG method in~\cite{giang2026conditional}, for a finite-sum setting, achieves an iteration complexity of $\mathcal{O}(\epsilon^{-4})$ for obtaining an $\epsilon$-accurate Frank--Wolfe stationarity measure and an $\epsilon$-accurate lower-level objective-gap guarantee. Moreover, it is shown that the liminf of the stationarity measure is zero and that the lower-level objective gap converges to zero. In comparison, our IPR-VFISTA$_{\mathrm{nc}}$ achieves lower total iteration complexity together with a sharper infeasibility guarantee.}
\end{remark}

\section{Numerical experiments}\label{sec: 6}
\ser{We evaluate the proposed algorithms in strongly convex and nonconvex upper-level settings using three ill-posed inverse problems and an MNIST classifier-selection problem.}
\subsection{Strongly convex upper-level setting}\label{sub:num_sc}
\ser{We evaluate IR-ISTA$_{\mathrm{s}}$ and IR-VFISTA$_{\mathrm{s}}$ and compare them with FBi-PG~\cite{merchav2024dynamic}, IR-IG~\cite{amini2019iterative}, BiG-SAM~\cite{sabach2017first}, Bi-SG~\cite{merchav2023convex}, and FIBA~\cite{helou2017}. \fyr{In \eqref{prob:num_c_sc}, we consider minimizing an elastic-net type regularizer in addressing} ill-posed \fyr{linear} inverse problems,} \ser{where $\mu_f>0$, $\mathbf{A}\in\mathbb{R}^{n\times n}$, and $b\in\mathbb{R}^n$. For the MNIST classifier-selection problem, using the same upper-level objective in \eqref{prob:num_c_sc}, we consider \eqref{eq:mnist_sbo_sc}. Here, $N$ is the sample size and $(u_i,v_i)\in\mathbb{R}^n\times\{-1,1\}$ denotes the $i$th labeled image.} 
\begin{align}
\min_{x}\quad
\bar{f}(x)
:=
\tfrac{\mu_f}{2}\|x\|^2+\|x\|_1,
\quad
\text{s.t.}\quad
x\in
\arg\min_{x\in\mathbb{R}^n}
\bar{h}(x)
:=
\tfrac{1}{2}\|\mathbf{A}x-b\|^2.
\label{prob:num_c_sc}
\\[0ex]
\ser{
\min_x\
\bar f(x),
\quad
\text{s.t.}\quad
x\in
\arg\min_{\|x\|_{\infty}\leq 0.5}
\bar h(x)
:=
\tfrac{1}{N}
\textstyle\sum_{i=1}^{N}
\log\left(1+\exp(-v_i \fyr{\langle u_i, x\rangle})\right)
}.
\label{eq:mnist_sbo_sc}
\end{align}
\subsubsection{Experiments and setup}
Similar to~\cite{beck2014first,amini2019iterative,sabach2017first}, we consider the Baart, Phillips, and Foxgood inverse problems from the Regularization Tools package\footnote{Available at \url{https://www2.imm.dtu.dk/~pcha/Regutools/}}. These problems differ in the construction of $\mathbf{A}\in\mathbb{R}^{n\times n}$ and $b\in\mathbb{R}^n$. We apply \ser{Algorithms~\ref{alg:IR-ISTA-s} and~\ref{alg:IR-R-VF_ne}} to \eqref{prob:num_c_sc} with $n=1000$ and $x_0=1_{n\times 1}$. \ser{For MNIST, we consider binary classification and initialize all methods at $x_0=\mathbf{1}_{784\times1}$.}
\subsubsection{Results and insights}
\begin{figure}[ht]
\centering
\setlength{\tabcolsep}{2pt}
\begin{tabular}{cccc}

\textbf{Baart} &
\textbf{Phillips} &
\textbf{Foxgood} &
\textbf{MNIST}
\\
\hline

\adjustbox{valign=m,height=1.55cm}{%
\includegraphics[width=0.22\textwidth]{
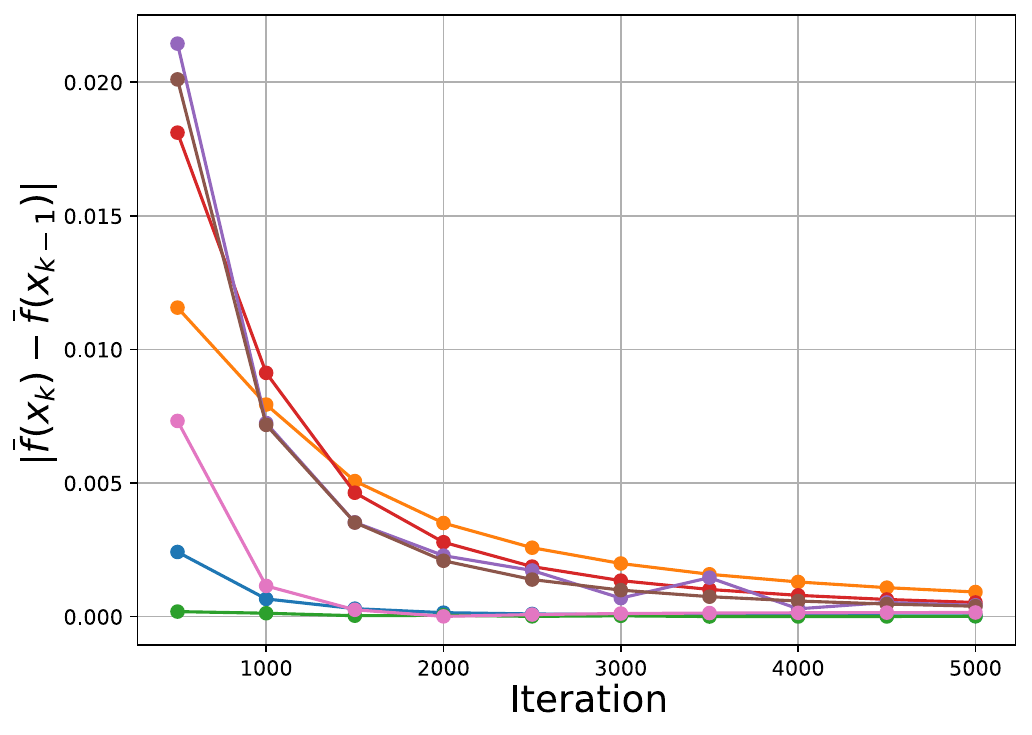}}
&
\adjustbox{valign=m,height=1.55cm}{%
\includegraphics[width=0.22\textwidth]{
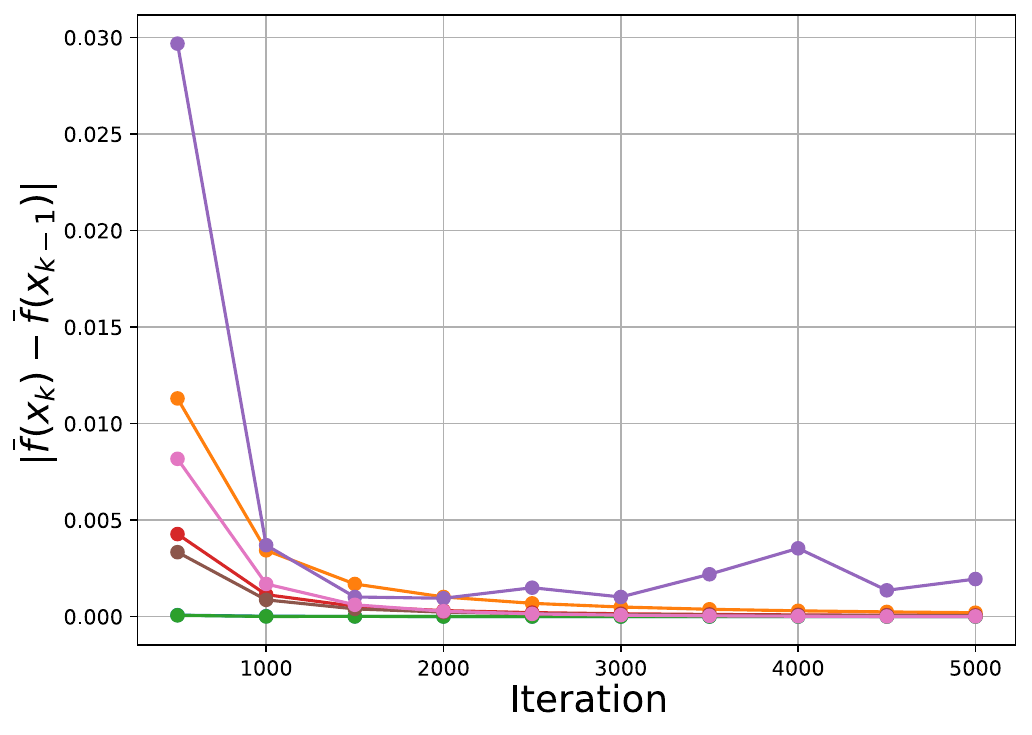}}
&
\adjustbox{valign=m,height=1.55cm}{%
\includegraphics[width=0.22\textwidth]{
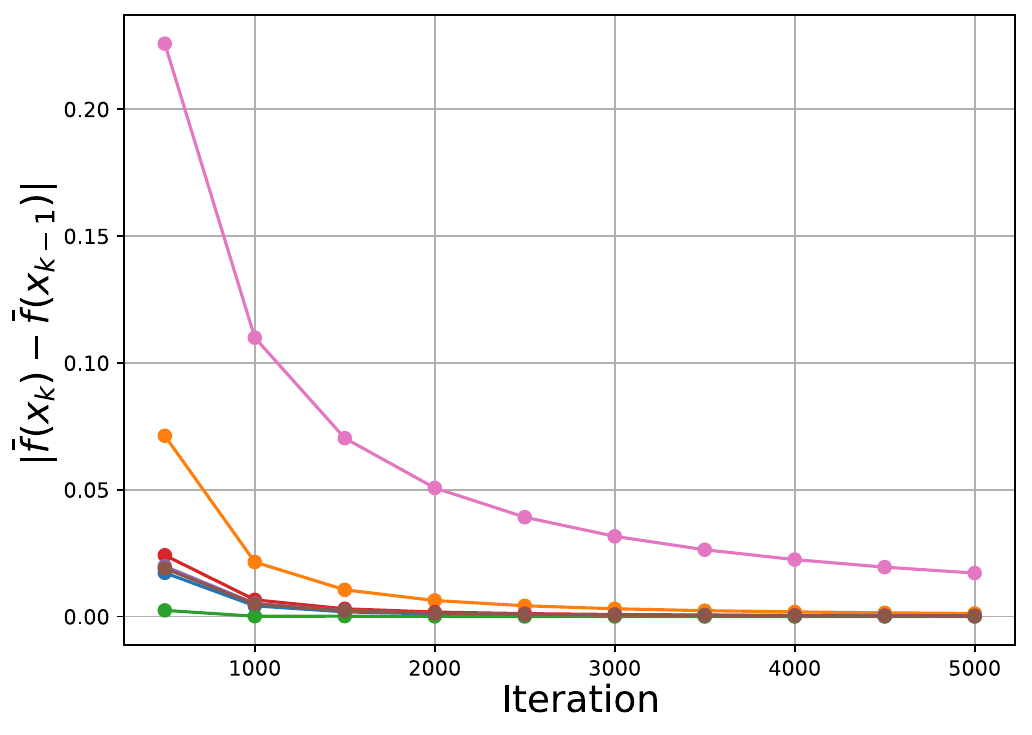}}
&
\adjustbox{valign=m,height=1.55cm}{%
\includegraphics[width=0.22\textwidth]{
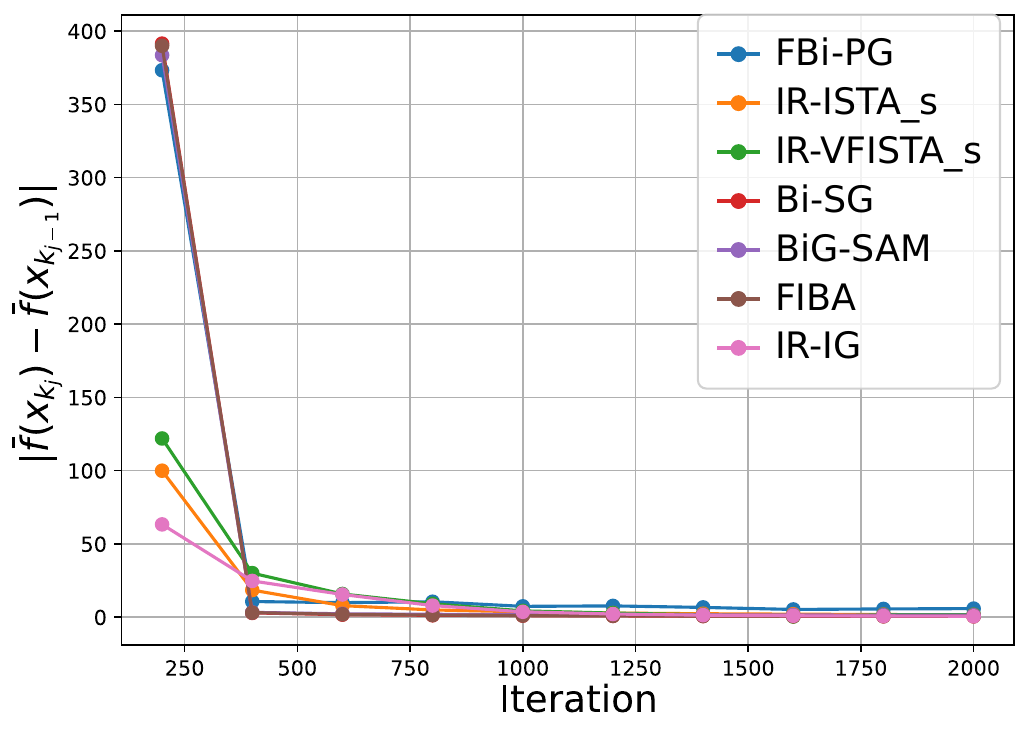}}
\\[2pt]

\adjustbox{valign=m,height=1.55cm}{%
\includegraphics[width=0.22\textwidth]{
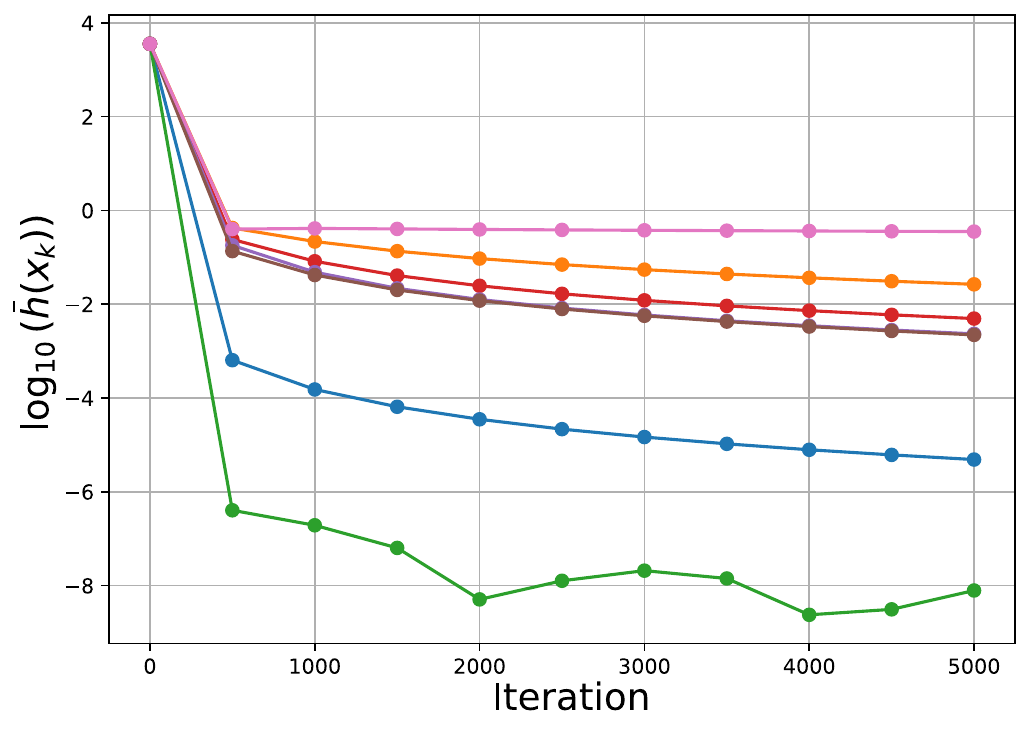}}
&
\adjustbox{valign=m,height=1.55cm}{%
\includegraphics[width=0.22\textwidth]{
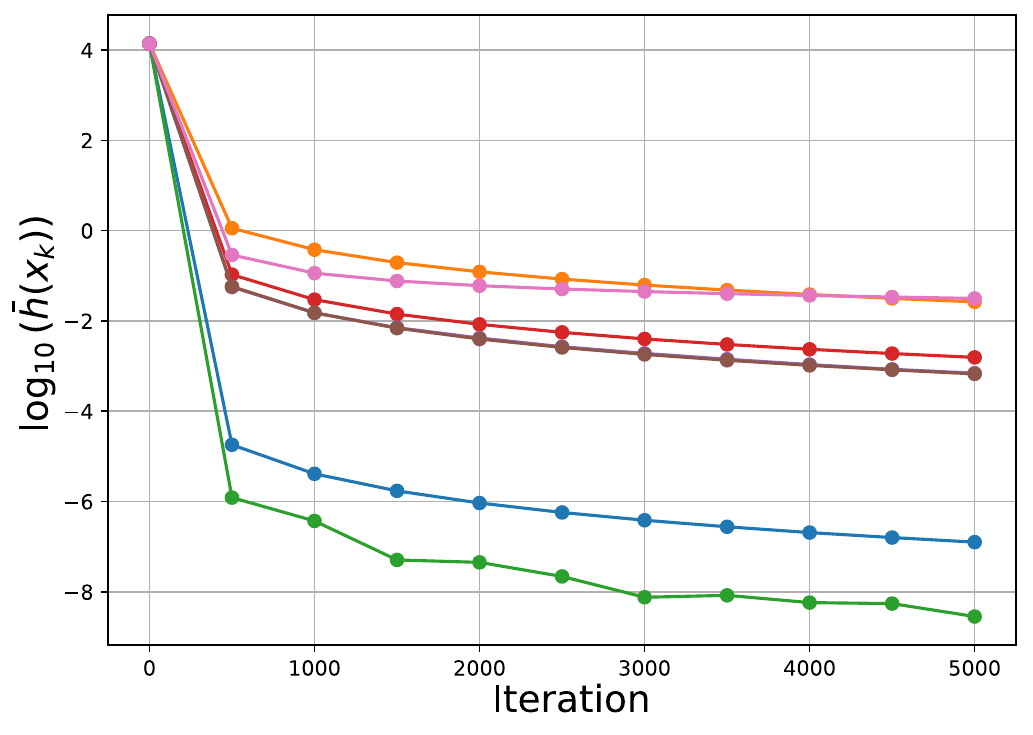}}
&
\adjustbox{valign=m,height=1.55cm}{%
\includegraphics[width=0.22\textwidth]{
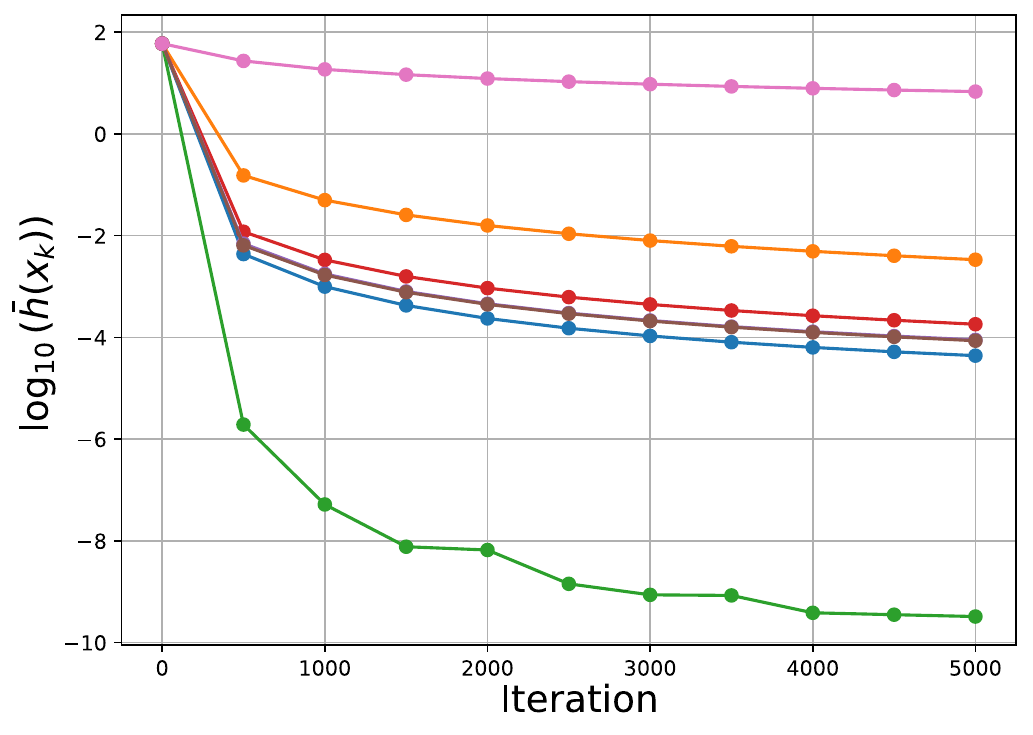}}
&
\adjustbox{valign=m,height=1.55cm}{%
\includegraphics[width=0.22\textwidth]{
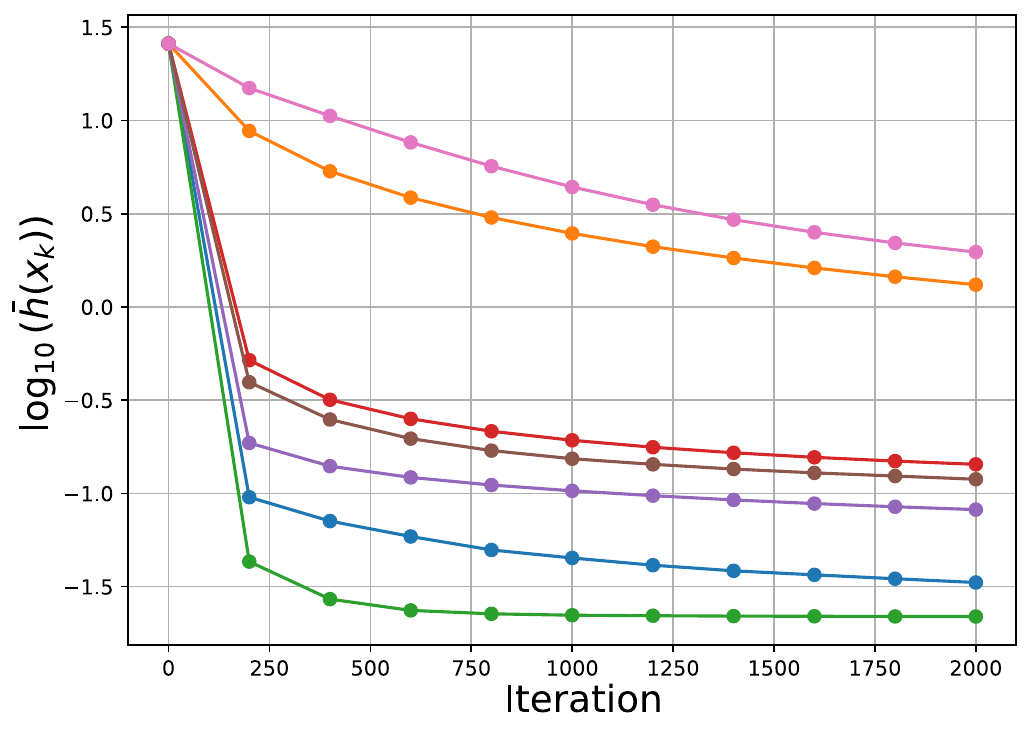}}
\\

\hline

\end{tabular}

\caption{\ser{\small
Comparison of FBi-PG~\cite{merchav2024dynamic}, IR-IG~\cite{amini2019iterative},
BiG-SAM~\cite{sabach2017first}, Bi-SG~\cite{merchav2023convex}, and
FIBA~\cite{helou2017} with our proposed IR-ISTA$_{\mathrm{s}}$ and
IR-VFISTA$_{\mathrm{s}}$ methods addressing
\eqref{prob:num_c_sc} and \eqref{eq:mnist_sbo_sc}.}}
\label{fig:comparison_1000}

\end{figure}
\ser{Figure~\ref{fig:comparison_1000} reports the successive upper-level objective change and the lower-level objective value for the Baart, Phillips, Foxgood, and MNIST problems. Across the four instances, the upper-level objective change generally decreases and stabilizes. For the inverse problems, IR-VFISTA$_{\mathrm{s}}$ reduces the lower-level objective faster than IR-ISTA$_{\mathrm{s}}$ and attains smaller final values. For MNIST, IR-VFISTA$_{\mathrm{s}}$ attains \ssg{the} smallest lower-level \fyr{objective} among the  methods.}
\subsection{Nonconvex upper-level setting}\label{sub:num_nc}
\ser{We next evaluate IPR-VFISTA$_{\mathrm{nc}}$ for a smooth nonconvex upper-level objective and compare it with CG-BiO~\cite{jiang2023conditional}. For both the inverse and MNIST classifier-selection problems, we use a smooth approximation of the log-sum penalty.}
\ser{Specifically, let}
$l(x) := \textstyle\sum_{i=1}^n \log (1 + {|x_i|}{\varepsilon}^{-1}),$
where $\varepsilon>0$. Although $l$ is nonsmooth and nonconvex~\cite{prater2022proximity}, its Moreau envelope with smoothing parameter $\delta>0$ is
$M_{l}^{\delta}(x)=\ch{l\left(\operatorname{prox}_{\delta l}[x]\right)}+\tfrac{1}{2\delta}{\|x-\operatorname{prox}_{\delta l}[x]\|}^2$
~\cite[Definition 6.52]{beck2017first}.
Moreover, the proximal operator is \fyr{separable}, i.e.,
\ch{$
\operatorname{prox}_{\delta l}[x]
=
\left(\operatorname{prox}_{\delta l_i}[x_i]\right)_{i=1}^n,
$}
where
$l_i(x_i)=\log (1+{|x_i|}{\varepsilon}^{-1}).$
If $\sqrt{\delta}\leq\varepsilon$, then~\cite[Proposition 1]{prater2022proximity}
$\operatorname{prox}_{\delta l_i}[x_i]=0$ for
$|x_i|\leq\tfrac{\delta}{\varepsilon}$, and
$\operatorname{prox}_{\delta l_i}[x_i]
=
0.5\operatorname{sign}(x_i)
({|x_i|-\varepsilon+\sqrt{(|x_i|+\varepsilon)^2-4\delta}})$
otherwise.
Consequently, $M_l^\delta$ is $\tfrac{1}{\delta}$-smooth, with
\ch{$
\nabla M_{l}^{\delta}(x)
=
\tfrac{1}{\delta}\left(x-\operatorname{prox}_{\delta l}[x]\right).
$}
\ser{We consider the following ill-posed inverse and classifier-selection problems, respectively\ssg{.}}
\begin{align}
\min_{x}\quad
\bar{f}(x)
:=
M_l^\delta(x),
\quad
\text{s.t.}\quad
x\in
\arg\min_{\|x\|\leq \ser{R}}
\bar{h}(x)
:=
\tfrac{1}{2}\|\mathbf{A}x-b\|^2.
\label{prob:num_ncx}
\\[0ex]
\ser{
\min_x\
\bar f(x),
\quad
\text{s.t.} \ 
x\in
\arg\min_{\|x\|_{\infty}\leq 0.5}
\bar h(x)
:=
\tfrac{1}{N}
\textstyle\sum_{i=1}^{N}
\log\left(1+\exp(-v_i \fyr{\langle u_i, x\rangle})\right)
}.
\label{eq:mnist_sbo_nc}
\end{align}
\ser{As in \eqref{eq:mnist_sbo_sc}, the lower-level problem minimizes the empirical logistic loss, while the upper-level objective \fyr{is a smooth approximation of a sparsity-promoting regularizer}.}
  
\subsubsection{Experimental setup}
\ser{We compare IPR-VFISTA$_{\mathrm{nc}}$ and CG-BiO~\cite{jiang2023conditional} on the Baart, Phillips, and Foxgood problems in~\eqref{prob:num_ncx} and the MNIST binary classification problem in~\eqref{eq:mnist_sbo_nc}. For the inverse problems, we set $n=1000$, $x_0=\mathbf{1}/\sqrt{n}$, and $R=1.05\|x^{\mathrm{true}}\|_2$, \fyr{where $x^{\mathrm{true}}$ denotes the solution of the corresponding inverse problem provided by the Regularization Tools package}. For MNIST, both methods start from the projection of the all-ones vector onto $\{x:\|x\|_\infty\leq0.5\}$.}
\ser{The horizontal axis reports the total computational budget.
For CG-BiO, this budget includes the Frank--Wolfe warm-up and main iterations, with ``pre'' denoting the warm-up iterations. For IPR-VFISTA$_{\mathrm{nc}}$, it denotes the cumulative number of inner accelerated proximal-gradient iterations.}
\ser{For CG-BiO, we use the values of
$\epsilon_{\bar f}=\epsilon_{\bar h}$ reported in the legends
and the constant stepsize prescribed in
\cite[Theorem~2]{jiang2023conditional}. 
We terminate the warm-up when the
lower-level Frank--Wolfe gap is at most $\epsilon_{\bar h}/2$ and solve
the linear minimization oracle  over the
feasible-set--cutting-plane intersection.}
\ser{For IPR-VFISTA$_{\mathrm{nc}}$, we use the $(\bar\eta,\hat{\gamma})$ pairs reported in the legends, and choose $\hat{\gamma}$, $J_k$ and $\eta_k$ according to Theorem~\ref{theorem:nonc_all}~(1). 
}
\begin{figure}[ht]
\centering
\setlength{\tabcolsep}{2pt}
\begin{tabular}{cccc}
\textbf{Baart} & \textbf{Phillips} & \textbf{Foxgood} & \textbf{MNIST} \\
\hline

\adjustbox{valign=m, height=1.55cm}{%
\includegraphics[width=0.22\textwidth]{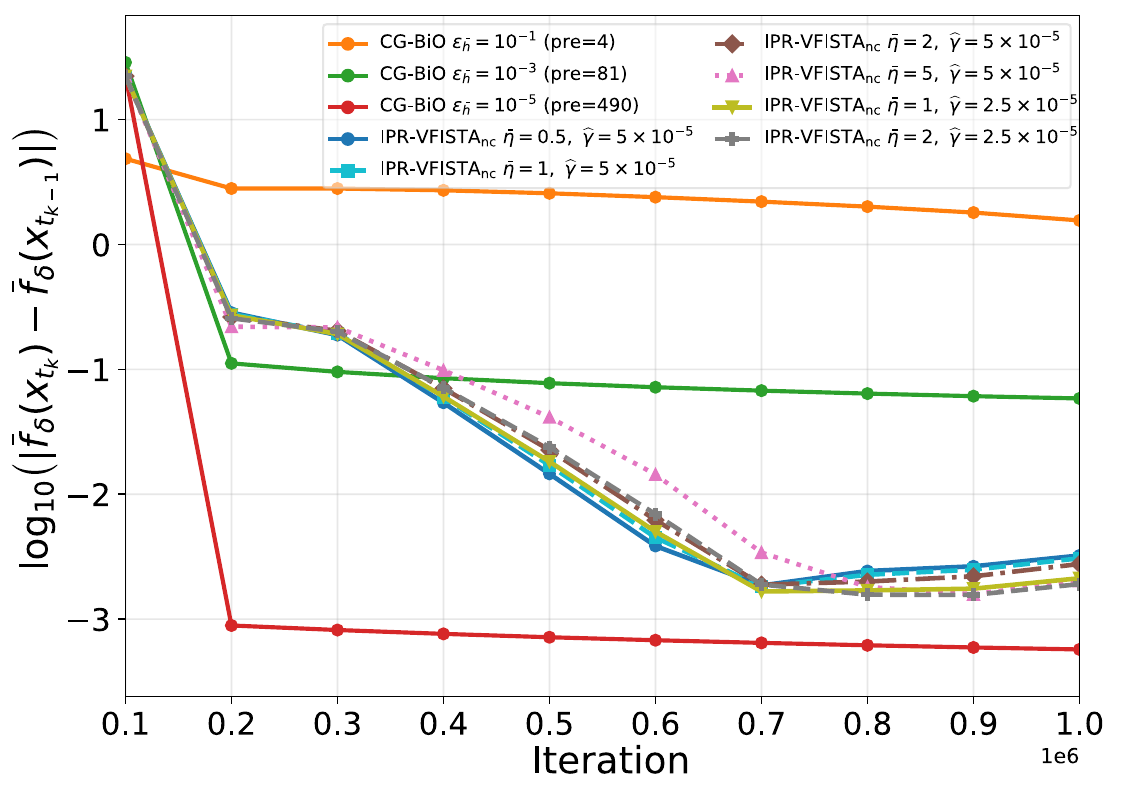}}
&
\adjustbox{valign=m, height=1.55cm}{%
\includegraphics[width=0.22\textwidth]{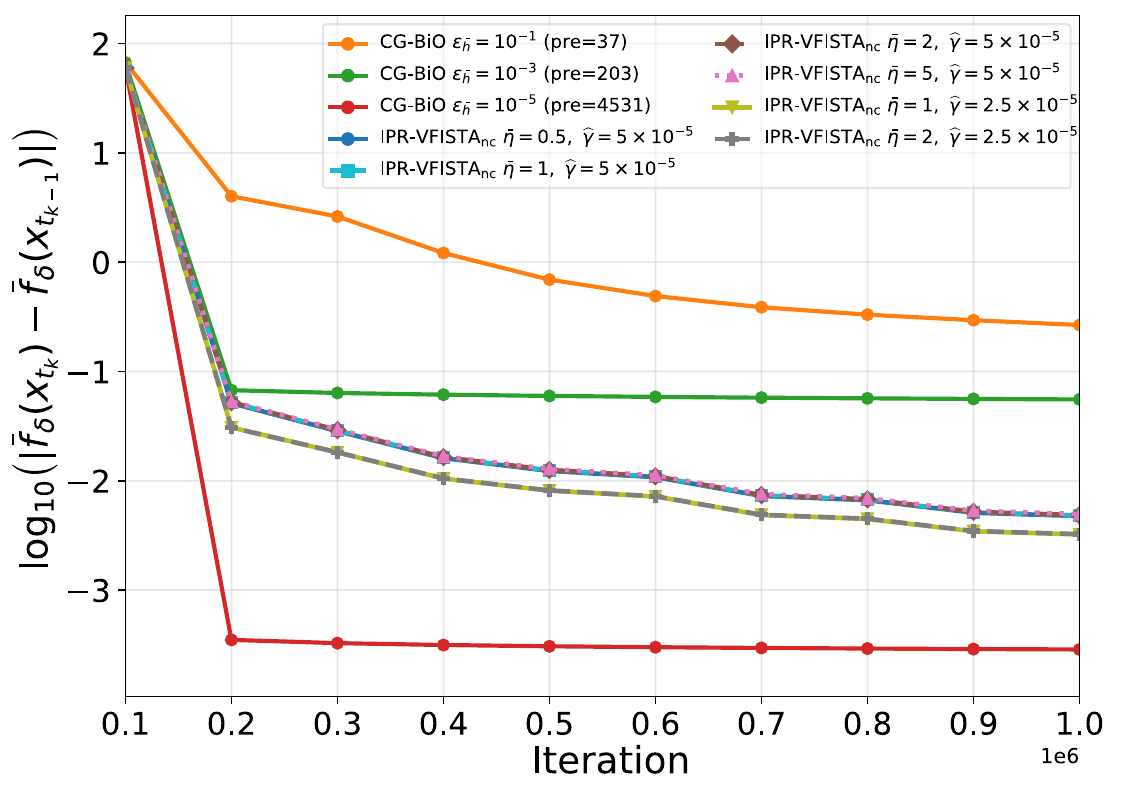}}
&
\adjustbox{valign=m, height=1.55cm}{%
\includegraphics[width=0.22\textwidth]{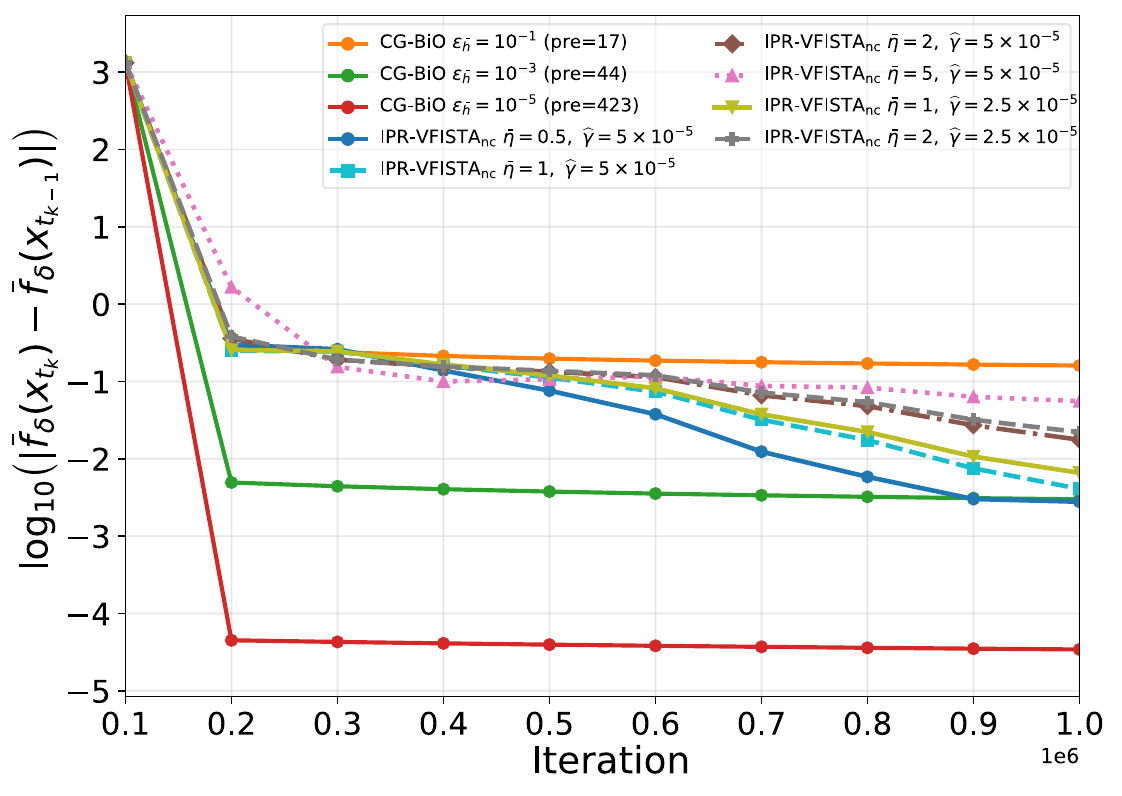}}
&
\adjustbox{valign=m, height=1.55cm}{%
\includegraphics[width=0.22\textwidth]{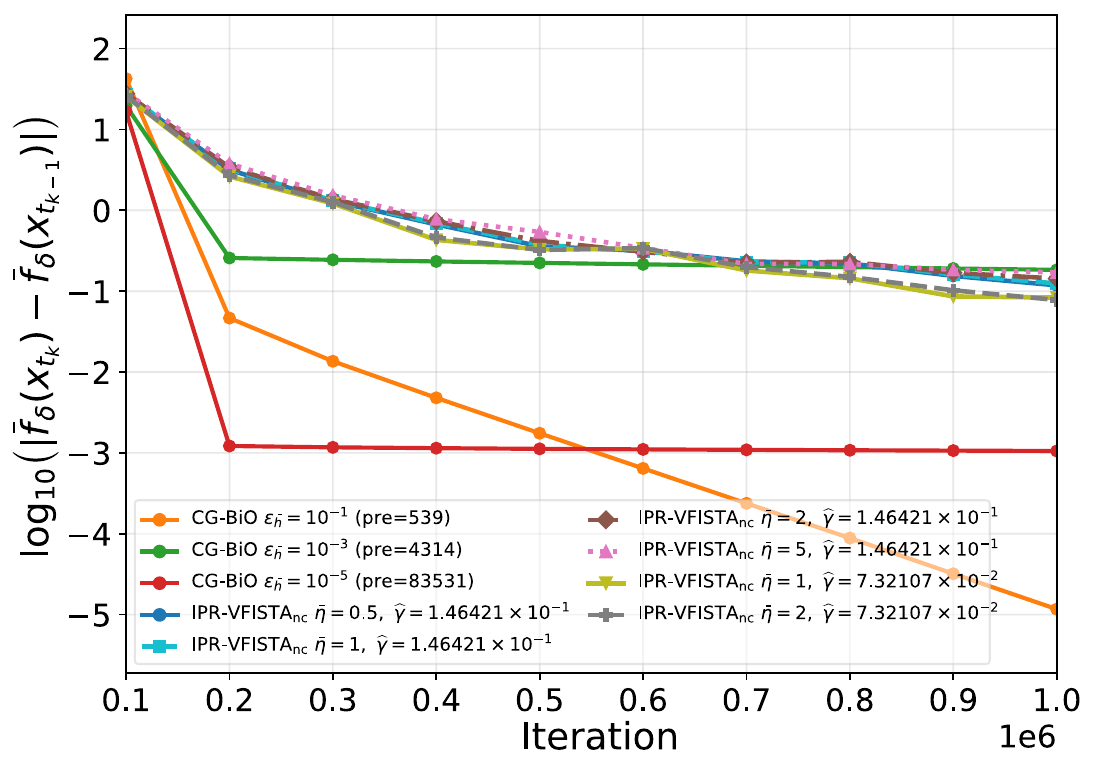}}
\\[2pt]

\adjustbox{valign=m, height=1.55cm}{%
\includegraphics[width=0.22\textwidth]{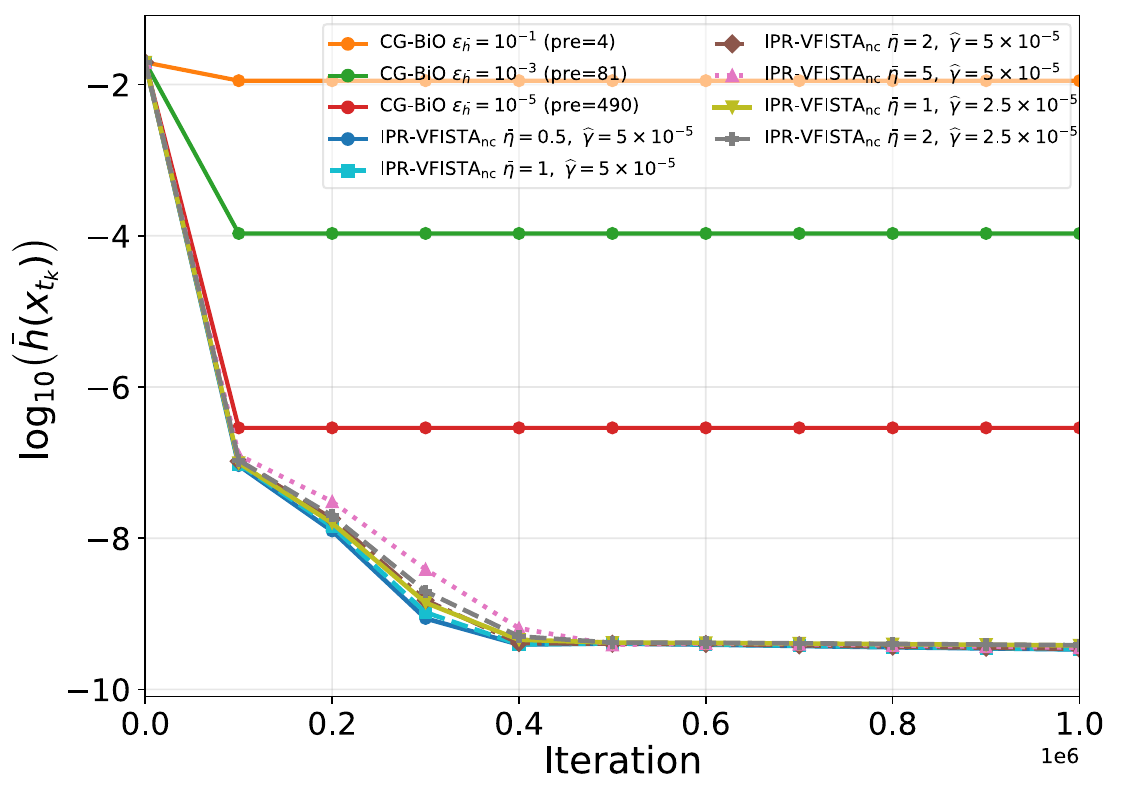}}
&
\adjustbox{valign=m, height=1.55cm}{%
\includegraphics[width=0.22\textwidth]{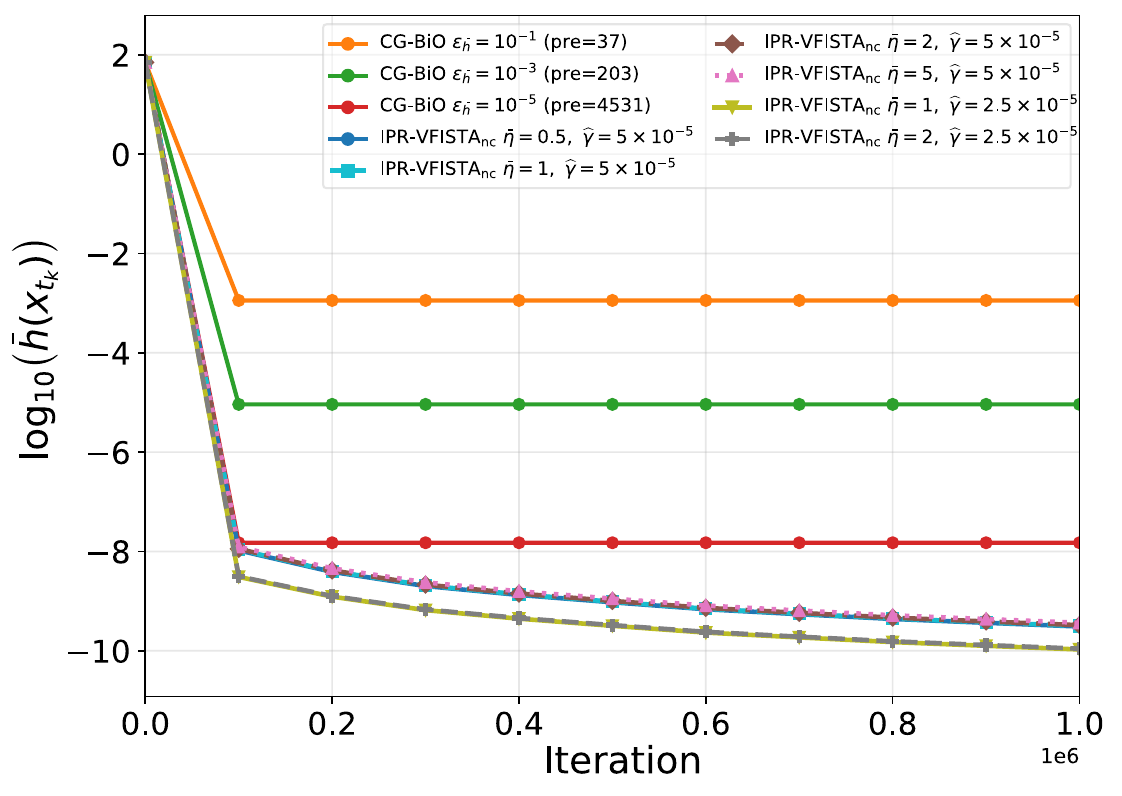}}
&
\adjustbox{valign=m, height=1.55cm}{%
\includegraphics[width=0.22\textwidth]{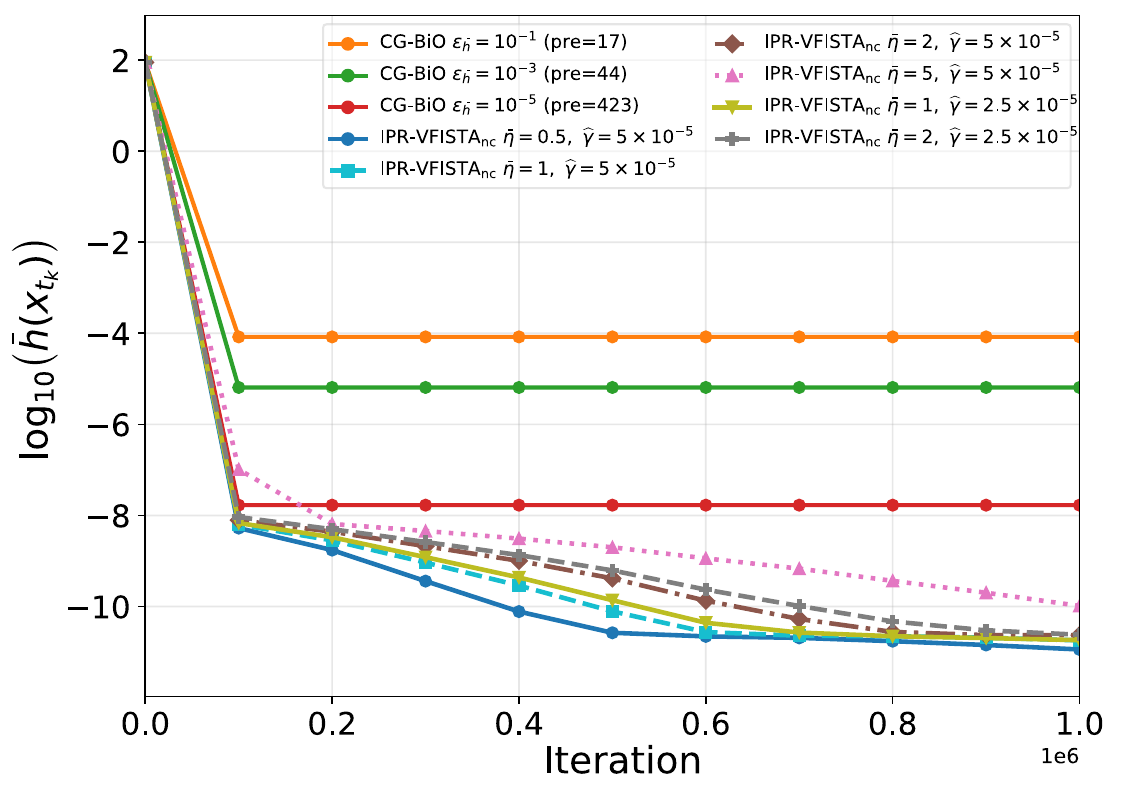}}
&
\adjustbox{valign=m, height=1.55cm}{%
\includegraphics[width=0.22\textwidth]{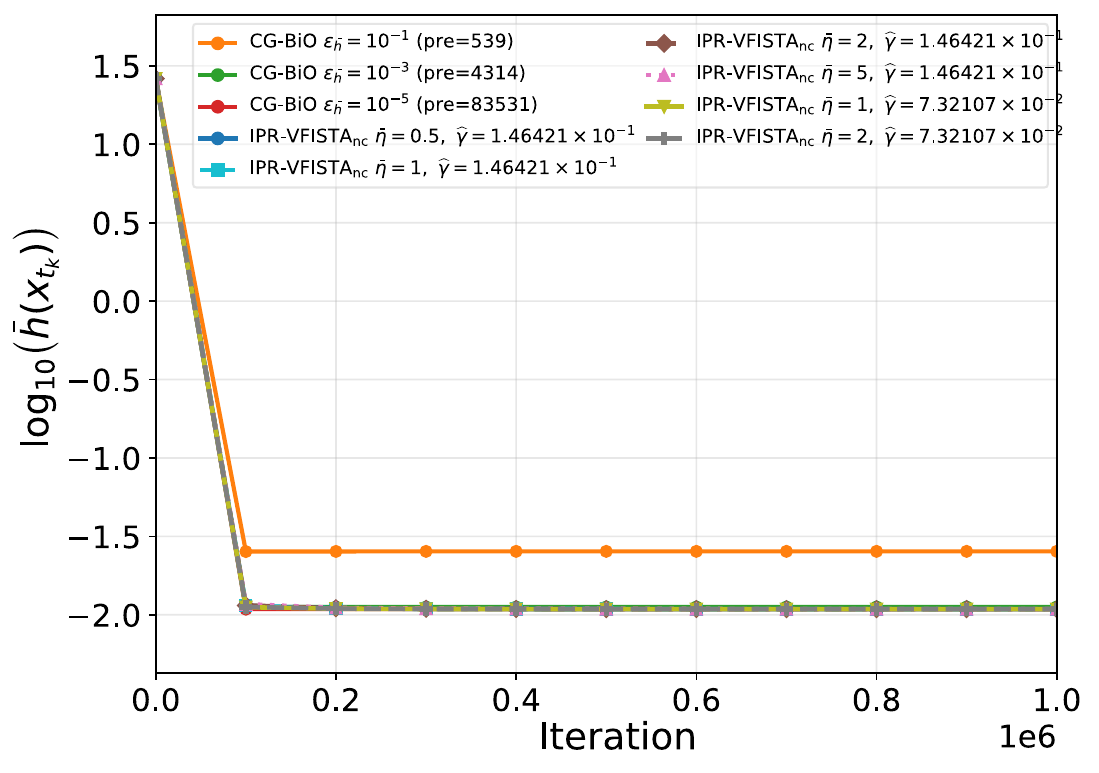}}
\\

\end{tabular}

\caption{\small
\ser{Comparison of $\text{IPR-VFISTA}_{\mathrm{nc}}$ and
CG-BiO~\cite{jiang2023conditional} addressing
\eqref{prob:num_ncx} and \eqref{eq:mnist_sbo_nc}.}}
\label{fig:num_result_ncv_n}
\end{figure}
\subsubsection{Results and insights}
\ser{Figure~\ref{fig:num_result_ncv_n} compares the methods under equal computational budgets. The first and second rows report successive upper-level changes and lower-level objective values, respectively.}
\ser{For the inverse problems, IPR-VFISTA$_{\mathrm{nc}}$ consistently decreases the upper-level changes and attains smaller final lower-level values than CG-BiO across all tested parameters. Its final upper-level changes are generally smaller than those of CG-BiO with $\bar{\epsilon}_h\in\{10^{-1},10^{-3}\}$, while CG-BiO with $\bar{\epsilon}_h=10^{-5}$ attains smaller changes at \fyr{a higher} preprocessing cost. For MNIST, the IPR-VFISTA$_{\mathrm{nc}}$ variants exhibit comparable upper-level changes and attain lower-level values comparable to CG-BiO with $\bar{\epsilon}_h\in\{10^{-3},10^{-5}\}$. Overall, IPR-VFISTA$_{\mathrm{nc}}$ consistently reduces the lower-level objective and admits asymptotic convergence, whereas CG-BiO~\cite{jiang2023conditional} provides nonasymptotic guarantees \ssg{but does not establish the asymptotic convergence result}.}
\section{Concluding remarks}\label{sec:7}
\fyr{In this work, we propose three iteratively regularized (IR) first-order methods for SBO problems. For SBO problems with strongly convex upper-level and convex lower-level objectives, we introduce IR proximal gradient and its accelerated variant and establish asymptotic convergence together with simultaneous nonasymptotic convergence rates for infeasibility and suboptimality measures. To the best of our knowledge, these guarantees provide the best-known simultaneous convergence rates for this class of composite SBO problems. For SBO problems with smooth nonconvex upper-level objectives, we propose an inexactly projected IR \fyr{accelerated} gradient method and establish an asymptotic stationarity guarantee without requiring \ssg{lower-level weak sharp minimality}. \ssg{Moreover, the total iteration complexity of the proposed method matches that of existing methods while providing a sharper lower-level infeasibility guarantee.}}

\bibliographystyle{siamplain}
\bibliography{references_fy}
\section{\texorpdfstring{\ssg{Appendix}}{Appendix}}
\subsection{\texorpdfstring{\ssg{Proof for Lemma~\ref{lemma:weaksharp_examples}}}{Proof for Lemma}}\label{subsection:proof_weaksharp_examples}
\begin{proof}
    \ser{(i) If $\min_{x \in \mathbb{R}^n} h(x) < h^* $, then \ser{the problem $\min_{\|x\|\leq 1} h(x):=\tfrac{1}{2}\|\mathbf{A}x-b\|^2$} admits the $\alpha$-weak sharp minimality of order $m=2$ in view of \cite[Lemma 3.6]{jiang2022holderian}. If $\min_{x \in \mathbb{R}^n}h(x) = h^*   $, this condition still holds under either of the following conditions \ser{~\cite[Lemma 3.7]{jiang2022holderian}}: (a) if $\lambda_{\text{min}} > 0$ in which the property holds with $\alpha = \sqrt{\tfrac{1}{\lambda_{\text{min}}}}$. (b) if $\lambda_{\text{min}} = 0$ and $\|{{\bf{Q}}^{\dagger}q} \|<1$.  Thus, there exists $\alpha > 0$ such that $
\alpha\ \operatorname{dist}(x,X_{\bar{h}^*})^2 \leq  \bar{h} (x) - \bar{h}^*,$  for all $x \in \mathbb{R}^n,
$ where $\bar{h}(x) = h(x) + \mathbb{I}_S(x)$, $\bar{h}^*$ is the optimal objective function of \ser{the problem}, and $S = \{x: \|x\| \leq 1 \}$.}

\ser{(ii)} \fyr{Define $X^*  := \arg\min_{x \in X} f(x)$ and $
f^*  := \min_{x \in X} f(x)$.   Let $x_{\mathcal R} :=  \Pi_{\mathcal R}[x]$ and $r := x - \Pi_{\mathcal R}[x]$. By the projection theorem,
$\langle r, z - \Pi_{\mathcal R}[x]\rangle \leq 0$ for all $z \in \mathcal R$. Since $\mathcal R$ is a subspace, both $z$ and $-z$ belong to $\mathcal R$. Applying the inequality to $z$ and $-z$ gives
$\langle r, z\rangle \leq \langle r, \Pi_{\mathcal R}[x]\rangle$ and  $\langle -r, z \rangle \leq \langle r ,\Pi_{\mathcal R}[x]\rangle $. Taking $z := 2\Pi_{\mathcal R}[x]$ and $z := -2\Pi_{\mathcal R}[x]$ yields $\langle r, \Pi_{\mathcal R}[x]\rangle \leq 0$ and $\langle r, \Pi_{\mathcal R}[x] \rangle \geq 0$,  implying that $\langle r,\Pi_{\mathcal R}[x]\rangle  = 0$. Combining the preceding relation with $\langle r, z - \Pi_{\mathcal R}[x]\rangle \leq 0$ yields $\langle r,z\rangle  \leq  0$ for all $z \in \mathcal R$. Since $-z \in \mathcal R$, we also have $\langle r, z\rangle  \geq  0$. Hence, $\langle r , z \rangle= 0$ for all $z \in \mathcal{R}$. From the definition of $r$, we conclude that any vector $x \in \mathbb{R}^n$ can be characterized uniquely as $
x = x_\mathcal{R} + x_\mathcal{N}$,  where  $x_\mathcal{R} \in \mathcal{R}$ and $ x_\mathcal{N} \in \mathcal{N}$ (cf. the orthogonal decomposition theorem). 
Since each $u_i \in \mathcal{R}$, we have $\langle u_i, x_\mathcal{N}\rangle =0$, hence
$
f(x)=f(x_\mathcal{R}).$ Note that  $
   \{  \Pi_\mathcal{R}[x] \mid x \in X\} = X_\mathcal{R}$. Define $g(z):= f(z)$, for $ z\in X_\mathcal{R}$. Then $
\min_{x\in X} f(x) = \min_{z\in X_\mathcal{R}} g(z).$ Next, we show that $f$ is strongly convex on $X_\mathcal{R}$. Let $\sigma(t) := 1/( 1 + \exp(-t))$ denote the logistic sigmoid function. Then, the Hessian of the logistic loss is $$
\nabla^2 f(x)
= \tfrac{1}{N}\textstyle\sum_{i=1}^N \sigma(v_i \langle u_i, x\rangle )\bigl(1-\sigma(v_i \langle u_i , x\rangle )\bigr)\, u_i u_i^\top.
$$
Since $X$ is bounded, there exists $B>0$ such that $|\langle u_i, x\rangle | \leq B$ for all $x \in X$ and all $i$. The function $
\phi(t) := \sigma(t)(1-\sigma(t))
$
is continuous and strictly positive on $[-B,B]$, and thus attains a strictly positive minimum
$
c_X := \min_{|t|\leq B} \sigma(t)(1-\sigma(t)) > 0.
$
Therefore,
$
\sigma(v_i \langle u_i, x\rangle )(1-\sigma(v_i \langle u_i, x\rangle)) \geq c_X
$
for all $x \in X$, which implies that 
$
\nabla^2 f(x)
\succeq c_X \left(\tfrac{1}{N}\sum_{i=1}^N u_i u_i^\top\right)
= c_X \Sigma,$ 
where $\Sigma := \tfrac{1}{N}\sum_{i=1}^N u_i u_i^\top.
$ Since $\mathcal{R} \neq \{0\}$, we have $\Sigma \succ 0$ on $\mathcal{R}$. Thus, there exists $\lambda_{\min}^\mathcal{R}(\Sigma) > 0$ such that
$
h^\top \Sigma h \geq \lambda_{\min}^\mathcal{R}(\Sigma)\|h\|^2$ for all  $h \in \mathcal{R}$.
Hence,  $g$ is $\mu$-strongly convex on $X_{\mathcal{R}}$ with
$
\mu := c_X \lambda_{\min}^\mathcal{R}(\Sigma).
$
Let $x_{\mathcal{R}}^* $ be the unique minimizer of $g$ over $X_{\mathcal{R}}$. Then, we have 
$
g(x_{\mathcal{R}}) - g(x_{\mathcal{R}}^* )
\geq \tfrac{\mu}{2}\|x_{\mathcal{R}}-x_{\mathcal{R}}^* \|^2.
$
Next, we characterize the optimal solution set $X^*$ in terms of $x^*_{\mathcal{R}}$. 
Since $f(x)=f(x_\mathcal{R}) = g(x_\mathcal{R})$ for all $x \in X$, any minimizer $x^* \in X^*$ must satisfy
$
g(\Pi_\mathcal{R}[x^*]) = f(x^*)=f^*.
$
Hence $\Pi_\mathcal{R}[x^*]$ is a minimizer of $g$. Since $g$ is strongly convex on $X_\mathcal{R}$, its minimizer is unique, and therefore
$
\Pi_\mathcal{R}[x^*] = x_\mathcal{R}^*.
$ This implies that
$$
X^*  = \{x \in X  \mid \Pi_\mathcal{R}[x]=x_\mathcal{R}^* \}.$$
Fix $x\in X$ and write $x=x_{\mathcal{R}} + x_{\mathcal{N}}$. Define
$z^* := x_{\mathcal{R}}^* + x_{\mathcal{N}}$.  Since $x_{\mathcal{R}}^* \in X_{\mathcal{R}} $ and
$x_{\mathcal{N}} \in X_{\mathcal{N}}$, we have $z^* \in X$. Moreover, since
$x_{\mathcal{N}} \in \mathcal{N}$, we have $f(z^*)=f(x_{\mathcal{R}}^*)=f^*$, 
which implies that $z^* \in X^*$. Therefore, 
$X^*=\{x_{\mathcal{R}}^* + y_{\mathcal{N}} : y_{\mathcal{N}} \in X_{\mathcal{N}}\}$. Using the orthogonality of $\mathcal{R}$ and $\mathcal{N}$, we obtain
\begin{align*}
\operatorname{dist}^2(x,X^*)
&=
\min_{y_{\mathcal{N}}\in X_{\mathcal{N}}}
\left\|
(x_{\mathcal{R}} - x_{\mathcal{R}}^*)
+ (x_{\mathcal{N}} - y_{\mathcal{N}})
\right\|^2 \\
&=\min_{y_{\mathcal{N}}\in X_{\mathcal{N}}}
\left( \|x_{\mathcal{R}} - x_{\mathcal{R}}^*\|^2
+ \|x_{\mathcal{N}} - y_{\mathcal{N}}\|^2 \right) \\
&= \|x_{\mathcal{R}} - x_{\mathcal{R}}^*\|^2,
\end{align*}
where the last equality follows by choosing $y_{\mathcal{N}}=x_{\mathcal{N}}$. 
Hence, by the strong convexity of $g$,
\begin{align*}
f(x)-f^*
&= g(x_{\mathcal{R}})-g(x_{\mathcal{R}}^*) \geq
\frac{\mu}{2}
\|x_{\mathcal{R}}-x_{\mathcal{R}}^*\|^2 =\frac{\mu}{2}\operatorname{dist}^2(x,X^*).
\end{align*}
Therefore, the problem satisfies weak sharp minimality with order $2$ and $\alpha:=\frac{\mu}{2}$.}
\end{proof} 

\subsection{\texorpdfstring{\ssg{Proof for Corollary~\ref{Cor:R-VFISTA}}}{Proof for Corollary}}\label{appendix:proof_cor18}
\begin{proof}
{{(1.i)}} The lower bound follows  from \eqref{eq:lower_f_strongly}. Next, we show the upper bound. \ch{From the condition \ch{$(\newchr{\tfrac{K}{\ln(K)}})^2 \geq \tfrac{ (L_h + \bar{\eta}L_f) (p+1)^2}{\mu_f \bar{\eta}}$}, we \ch{have} \ch{$\bar{\eta} \geq  \tfrac{(L_h + \bar{\eta}L_f)}{\mu_f}  (\tfrac{(p+1) \ln({K})}{ K})^2.$}} This implies that $ \eta \leq \bar{\eta} $ \ch{and thus $(1 - \tfrac{1}{\sqrt{\kappa}} )^K \leq (1 - \tfrac{1}{\sqrt{\bar\kappa}} )^K$, where $\bar\kappa = \tfrac{L_h + \bar{\eta} L_f}{\eta{\mu_f}}$. } Now, by invoking the preceding inequality in \eqref{prop:R-VFISTA}, for any $k=K$ and any $x \in \mathbb{R}^n$, 
 \begin{align*}
 &{\eta}\bar{f}(x_K) - {\eta}\bar{f}(x) + \bar{h}(x_K) - \bar{h}(x) 
  \leq (1 - \tfrac{1}{\sqrt{\bar\kappa}} )^K   ( \bar{g}_{\eta}(x_0) - \bar{g}_{\eta}(x)  + \tfrac{\eta {\mu_f}}{2} {\|x_0 - x\|}^2 ),
 \end{align*}
By substituting  $x=x^*$,  using $0 \leq \bar{h}(x_K) - \bar{h}^* $, and dividing both sides by $\eta$, we obtain 
   \begin{align}\label{in:up_bound_e}
 & \bar{f}(x_K) - \bar{f}^*
  \leq   ( \bar{f}(x_0) - \bar{f}^* + \tfrac{\bar{h}(x_0) -\bar{h}^*}{\eta}  + \tfrac{{\mu_f}}{2} {\|x_0-x^*\|}^2 ) (1 - \tfrac{1}{\sqrt{\bar\kappa}}  )^K.
 \end{align}
From \ch{$\eta = \tfrac{(L_h + \bar{\eta} L_f)}{\mu_f} (\tfrac{(p+1) \ln(K)}{K})^2$} and that $1 - x \leq \exp(-x)$ for any $x \in \mathbb{R}$, we have
 \begin{align}\label{in:1-e}
 (1 - \tfrac{1}{\sqrt{\bar\kappa}}  )^{K} \leq \exp( \tfrac{-K}{\sqrt{\bar\kappa}})  
=\exp\left({- (p+1) \ln({K}) }\right)\leq \ch{\tfrac{1}{K^{(p+1)}}}.
\end{align}
The upper bound in (1.i) is obtained by substituting the value of $\eta$ in $\tfrac{\bar{h}(x_0) -\bar{h}^*}{\eta}$ in \eqref{in:up_bound_e} and then, by invoking \eqref{in:1-e}. 

\noindent {{(1.ii)}} Note that $\bar{h}^*:= \inf_{x\in \mathbb{R}^n} \bar{h}(x)$ implies that the lower bound for infeasibility is zero.  Now,  consider \eqref{prop:R-VFISTA} for $k=K$ and $x:=\Pi_{X^*_{\bar{h}}}[x_0]$. Then, by recalling $\bar{h}(\Pi_{X^*_{\bar{h}}}[x_0])=\bar{h}^*$ and  ${\|x_0 - \Pi_{X^*_{\bar{h}}}[x_0]\|}^2=\operatorname{dist}^2(x_0, X^*_{\bar{h}})$ and also by considering $\hat{C}_{\bar{f}} \leq \bar{f}(x_K)$ and $\hat{C}_{\bar{f}} \leq \bar{f}(\Pi_{X^*_{\bar{h}}}[x_0])$,  we obtain  
    \begin{align*}
  \bar{h}(x_K) -{\bar{h}^*}  &\leq 
    \eta (1 - \tfrac{1}{\sqrt{\bar\kappa}} )^K    (  \bar{f}(x_0) - \hat{C}_{\bar{f}}      + \tfrac{{\mu_f}}{2} \operatorname{dist}^2(x_0, X^*_{\bar{h}})  )\\
  & +(1 - \tfrac{1}{\sqrt{\bar\kappa}} )^K(\bar{h}(x_0) - \bar{h}^*) + \eta  {( \bar{f}(\Pi_{X^*_{\bar{h}}}[x_0])-\hat{C}_{\bar{f}})} ,
 \end{align*}
 where we used $\bar{\eta} > \eta$. By substituting \ch{$\eta=\tfrac{(L_h + \bar{\eta}L_f)}{\mu_f} (\tfrac{(p+1) \newchr{\ln({K})}}{ K})^2$}  and applying \eqref{in:1-e} to the preceding inequality, we obtain the  result. 
 
\ser{\noindent {{(1.iii)}}}  \fyr{Applying} the weak sharp minimality of $ X^*_{\bar{h}} $ as defined in Definition~\ref{def:weaksharp} to \ser{\eqref{prob:distance-weak-strongly}}, \fyr{together with the bounds established in (1.i) and (1.ii)}, we obtain the result.

\ser{\noindent {{(1.iv)}}}  \fyr{Consider} \eqref{eq:lower_f_strongly}. \fyr{Dropping} the nonnegative term from the \ser{left}-hand side and \ser{using} Definition~\ref{def:weaksharp}, \fyr{together with the bounds established in (1.ii)}, we obtain the result.

\noindent {{(2.i)}}  
Note that $\bar{h}^* := \inf_{x\in \mathbb{R}^n} \bar{h}(x)$ implies that infeasibility is bounded below by zero.
Consider   \eqref{prop:R-VFISTA} for $k:=K$ and \fyr{$x:=x^*$}. Then, from \eqref{eq:lower_f_strongly}, we obtain
\begin{align*}
 & \bar{h} (x_K)-\bar{h}^*-\eta{\|\tilde{\nabla} \bar{f}(x^*)\| \operatorname{dist}(x_K,X^*_{\bar{h}})} + \tfrac{\eta \mu_f}{2} {\|x_0 - x^*\|}^2  
 \\\notag &\leq (1 - \tfrac{1}{\sqrt{\kappa}} )^K   ( \bar{g}_{\eta}\left(x_0\right) - \bar{g}_{\eta}(x^*) + \tfrac{\eta {\mu_f}}{2} {\|x_0-x^*\|}^2 ).
\end{align*}
\noindent From the $\alpha$-weak sharp minimality of $X^*_{\bar{h}}$ and dropping the nonnegative term from the left-hand side, 
\ser{and in} view of $\|\tilde{\nabla} \bar{f}(x^*)\| \leq \tfrac{\alpha}{2 \eta}$, we obtain 
\begin{align*} 
    &{ \bar h(x_K) -\bar h^*}   \leq (1 - \tfrac{1}{{\sqrt{\kappa}}} )^K    ( 2(\bar{g}_{\eta}(x_0) - \bar{g}_{\eta}(x^*)) + {\eta {\mu_f}} {\|x_0-x^*\|}^2  ).
\end{align*}
\noindent Therefore, in view of Definition~\ref{def:uni_funcs}, we obtain the  result.

\noindent{{(2.ii)}} This result follows by recalling  Definition~\ref{def:weaksharp} and using \fyr{Corollary}~\ref{Cor:R-VFISTA}~(2.i).

\noindent{{(2.iii)}}
\noindent Now, consider \eqref{prop:R-VFISTA} for $k:=K$  \fyr{and $x:=x^*$. We obtain}
\begin{align}\label{eq:f_s_c_w_known}
   \bar{f}(x_K)-{\bar{f}^*}& \leq   ( \bar{f}(x_0) - \bar{f}^*+\tfrac{1}{\eta} (\bar{h}(x_0) - \bar{h}^*)  + \tfrac{ {\mu_f}}{2} {\|x_0 - x^*\|}^2  ) (1 - \tfrac{1}{\sqrt{\kappa}} )^K.
\end{align}
To obtain the lower bound, 
consider  \eqref{eq:lower_f_strongly}. By dropping $\tfrac{\mu_f}{2} {\|{x}_{k+1} - x^*\|}^2$ from the left-hand side, we obtain
$
-\|\tilde{\nabla} \bar{f}(x^*)\|  \operatorname{dist}({x}_{k+1},X^*_{\bar{h}}) \leq     \bar{f}({x}_{k+1}) -\bar{f}^*.$
Then, from \fyr{Corollary}~\ref{Cor:R-VFISTA}~(2.ii) and $ \|\tilde{\nabla} \bar{f}(x^*)\| \leq \tfrac{\alpha}{2\eta}$ in the preceding inequality, we obtain the result.

\noindent {{(2.iv)}} Consider \eqref{eq:lower_f_strongly}. From $\|\tilde{\nabla} \bar{f}(x^*)\| \leq \tfrac{\alpha}{2 \eta}$ and  \fyr{Corollary}~\ref{Cor:R-VFISTA} (2.ii),  we obtain
\begin{align*} 
\tfrac{\mu_f}{2} {\|x_{\fyr{K}} - x^*\|}^2 &\leq  \bar{f}(x_K) - \bar{f}^*   +       (   \bar{f}(x_0)-\bar{f}^* +\tfrac{1}{\eta}(\bar{h}(x_0)- \bar{h}^*) \ch{+} \tfrac{ {\mu_f  }}{2}{\|x_0-x^*\|}^2  ) (1 - \tfrac{1}{\sqrt{\kappa}} )^K.
\end{align*}
By using the upper bound for $ \bar{f}(x_K)-{\bar{f}^*}$ in \eqref{eq:f_s_c_w_known}, we arrive at the  result.
\end{proof}
\end{document}